\documentclass[reqno]{amsart}
\usepackage{amssymb, hyperref}

\AtBeginDocument{
\vspace{9mm}}

\begin{document}
\title[
]
{SCATTERING FOR THE RADIAL FOCUSING INLS EQUATION IN HIGHER DIMENSIONS } 
\author[L.G FARAH, C.M.GUZM\'AN]
{LUIZ GUSTAVO FARAH AND CARLOS M. GUZM\'AN}  

\address{LUIZ G. FARAH \hfill\break
Department of Mathematics, Federal University of Minas Gerais, BRAZIL}
\email{lgfarah@gmail.com}

\address{CARLOS M. GUZM\'AN \hfill\break
Department of Mathematics, Federal University of Minas Gerais, BRAZIL}
\email{carlos.guz.j@gmail.com}


\begin{abstract}
We consider the inhomogeneous nonlinear Schr\"odinger equation 
$$
i u_t +\Delta u+|x|^{-b}|u|^\alpha u = 0, 
$$
where $\frac{4-2b}{N}<\alpha<\frac{4-2b}{N-2}$ (when $N=2$, $\frac{4-2b}{N}<\alpha<\infty$) and $0<b<\min\{N/3,1\}$. For a radial initial data $u_0\in H^1(\mathbb{R}^N)$ under a certain smallness condition we prove that the corresponding solution is global and scatters. The smallness condition is related to the ground state solution of $-Q+\Delta Q+ |x|^{-b}|Q|^{\alpha}Q=0$ and the critical Sobolev index $s_c=\frac{N}{2}-\frac{2-b}{\alpha}$. This is an extension of the recent work \cite{paper2} by the same authors, where they consider the case $N=3$ and $\alpha=2$. The proof is inspired by the concentration-compactness/rigidity method developed by Kenig-Merle \cite{KENIG} to study $H^1(\mathbb{R}^N)$-critical problem and also Holmer-Roudenko \cite{HOLROU} in the case of $H^1(\mathbb{R}^N)$-subcritical equations.
\end{abstract}
\maketitle
\numberwithin{equation}{section}
\newtheorem{theorem}{Theorem}[section]
\newtheorem{proposition}[theorem]{Proposition}
\newtheorem{lemma}[theorem]{Lemma}
\newtheorem{corollary}[theorem]{Corollary}
\newtheorem{remark}[theorem]{Remark}
\newtheorem{definition}[theorem]{Definition}
\section{Introduction}

\ This paper is concerned with the initial value problem (IVP) for the focusing inhomogenous nonlinear Schr\"odinger (INLS) equation
\begin{equation}\label{INLS}
\left\{\begin{array}{cl}
i\partial_tu +\Delta u + |x|^{-b} |u|^\alpha u =0, & \;\;\;t\in \mathbb{R} ,\;x\in \mathbb{R}^N,\\
u(0,x)=u_0(x), &
\end{array}\right.
\end{equation}
where $u = u(t, x)$ is a complex-valued function in space-time $\mathbb{R}\times \mathbb{R}^N$ and $b>0$.

\ Note that when $b=0$ the above equation is the classical nonlinear Schr\"odinger equation (NLS) which appears in the description of nonlinear waves for various physical phenomena. On the other hand, in the end of the last century, it was suggested that in some situations laser beam propagation can be modeled by the inhomogeneous nonlinear Schr\"odinger equation in the following form 
\begin{equation}\label{INLS*}
i\partial_tu +\Delta u + V(x) |u|^\alpha u =0.
\end{equation}
We refer the reader to Gill \cite{GILL} and Liu-Tripathi \cite{LIU} for more physical details. From the mathematical point of view, the INLS model \eqref{INLS*} has been investigated by several authors over the last two decades. For instance, Merle \cite{MERLE} and Rapha\"el-Szeftel \cite{RAFAEL}, assuming $k_1<V(x)<k_2$ with $k_1,k_2>0$, study the problem of existence/nonexistence of minimal mass blow-up solutions. Fibich-Wang \cite{FIBI}, for $V(\epsilon |x|)$ with $\epsilon$ small and $V\in C^4(\mathbb{R}^N)\cap  L^{\infty}(\mathbb{R}^N)$, consider the stability of solitary waves. We should point out that in all these works the authors assume that $V(x)$ is a bounded function, so the well-posedness theory for the NLS equation can be directly applied also in this case. However, such assumption does not hold for the INLS equation \eqref{INLS} and several challenging technical difficulties arise in its study.

We briefly review some existence results available in the literature. Let us first introduce the following number:
\begin{equation}\label{def2*}
2^*:=\left\{\begin{array}{cl}
\frac{4-2b}{N-2}&\;\;\;\;N\geq 3,\\
\infty &\;\;\;\;\;N=2.
\end{array}\right.
\end{equation}
Genoud and Stuart \cite{GENOUD}-\cite{GENSTU}, using the abstract theory developed by Cazenave \cite{CAZEWE}, showed that \eqref{INLS} is locally well-posed in $H^1(\mathbb{R}^N)$ if $0<\alpha <2^*$ and globally if $\frac{4-2b}{N}<\alpha <2^{\ast}$ for small initial data. Recently, the second author in \cite{CARLOS} gave an alternative proof of these results, using the contraction mapping principle based on the Strichartz estimates satisfied by the linear flow. This new approach will be very important to carry out the analysis in the present study. For other recent works about the INLS model we refer the reader to Hong \cite{HONG}, Killip-Murphy-Visan-Zheng \cite{KMVZ2016} and Combet-Genoud\cite{COMGEN}.

We focus on the $L^2$-supercritical and $H^1$-subcritical case. Let us briefly explain this terminology. For a fixed $\delta >0$, the rescaled function  $u_\delta(t,x)=\delta^{\frac{2-b}{\alpha}}u(\delta^2 t,\delta x)$ is solution of \eqref{INLS} if only if $u(t,x)$ is. This scaling property gives rise to a scale-invariant norm. Indeed, computing the homogeneus Sobolev norm of $u_\delta(0,x)$ we have
$$
\|u_\delta(0,.)\|_{\dot{H^s}}=\delta^{s-\frac{N}{2}+\frac{2-b}{\alpha}}\|u_0\|_{\dot{H^s}}.
$$
Thus, the Sobolev space $H^{s_c}(\mathbb{R}^N)$, with $s_c=\frac{N}{2}-\frac{2-b}{\alpha}$, is invariant under the above scaling. The number $s_c$ is commonly referred  as the critical Sobolev index. Now, $L^2$-supercritical and $H^1$-subcritical equations refer to the case where $0<s_c<1$. A simple computation shows that the last relation is equivalent to $\frac{4-2b}{N}<\alpha <2^*$.

\ It is well-known that the INLS equation \eqref{INLS} has the following conserved quantities
\begin{equation}\label{mass}
M[u(t)]=\int_{\mathbb{R}^N}|u(t,x)|^2dx
\end{equation}
and
\begin{equation}\label{energy}
E[u(t)]=\frac{1}{2}\int_{\mathbb{R}^N}| \nabla u(t,x)|^2dx-\frac{1}{\alpha +2}\left\| |x|^{-b}|u|^{\alpha +2}\right\|_{L^1_x},
\end{equation}
which are called Mass and Energy, respectively. Furthermore, since 
\begin{equation}\label{QI1}
\|u_\delta\|_{L^2_x}=\delta^{-s_c}\|u\|_{L^2_x},\;\;\;\;\|\nabla u_\delta\|_{L^2_x}=\delta^{1-s_c}\|\nabla u\|_{L^2_x}
\end{equation}
and
$$
\left \| |x|^{-b}|u_\delta|^{\alpha+2} \right\|_{L^1_x}=\delta^{2(1-s_c)}\left \| |x|^{-b}|u|^{\alpha+2} \right\|_{L^1_x},
$$
it is easy to see that the following quantities are scale invariant
\begin{equation}\label{QI2}
E[u_\delta]^{s_c}M[u_\delta]^{1-s_c}=E[u]^{s_c}M[u]^{1-s_c},\;\;\|\nabla u_\delta\|^{s_c}_{L^2_x}\|u_\delta\|^{1-s_c}_{L^2_x}=\|\nabla u\|^{s_c}_{L^2_x}\|u\|^{1-s_c}_{L^2_x}.
\end{equation}
These quantities were introduced in Holmer-Roudenko \cite{HOLROU} (see also Duyckaerts-Holmer-Roudenko \cite{DUCHOLROU}) in order to describe the dichotomy between blowup/global regularity for the $3D$ cubic nonlinear Schr\"odinger equation (NLS). Here, these quantities also play an important role in our analysis.

The main goal is to extend our result in \cite{paper2} to general dimensions $N\geq 2$. More precisely, we want to obtain sufficient conditions on the initial data $u_0\in H^1(\mathbb{R}^N)$ such that the corresponding solution is global and scatters according to the next definition. 

\begin{definition}\label{defscattering}
A global solution $u(t)$ to the Cauchy problem \eqref{INLS} scatters forward in time
in $H^1(\mathbb{R}^N)$, if there exists  $\phi^+\in H^1(\mathbb{R}^N)$ such that
$$
\lim_{t\rightarrow +\infty}\|u(t)-U(t)\phi^+\|_{H^1}=0.
$$
Also, we say that  $u(t)$ scatters backward in time if there exists $\phi^-\in H^1(\mathbb{R}^N)$ such that 
$$
\lim_{t\rightarrow -\infty}\|u(t)-U(t)\phi^-\|_{H^1}=0.
$$
Here, $U(t)$ denotes unitary group associated with the linear equation $i\partial_tu +\Delta u=0$, with initial data $u_0$.
\end{definition}

\ The global theory for the $L^2$-supercritical and $H^1$-subcritical INLS equation \eqref{INLS} was already investigated by the first author in \cite{LG}, where he proved the following result.

\begin{theorem}\label{TG}
	Let $\frac{4-2b}{N}<\alpha <2^*$ (or equivaly $0<s_c<1$) and $0<b<\min\{2,N\}$. Suppose that $u(t)$ is the solution of \eqref{INLS} with initial data $u_0\in H^1(\mathbb{R}^N)$ satisfying
\begin{equation}\label{EMC} 
E[u_0]^{s_c}M[u_0]^{1-s_c}<E[Q]^{s_c}M[Q]^{1-s_c} 
\end{equation}
and 
\begin{equation}\label{GFC}
\|  \nabla u_0 \|_{L^2}^{s_c} \|u_0\|_{L^2}^{1-s_c}<\|\nabla Q \|_{L^2}^{s_c} \|Q\|_{L^2}^{1-s_c},
\end{equation}
then $u(t)$ is a global solution in $H^1(\mathbb{R}^N)$. Furthermore, for any $t\in \mathbb{R}$ we have
\begin{equation}\label{GR}
\|  \nabla u(t) \|_{L^2}^{s_c} \|u(t)\|_{L^2}^{1-s_c}<\|\nabla Q \|_{L^2}^{s_c} \|Q\|_{L^2}^{1-s_c},
\end{equation}
where $Q$ is unique positive solution of the elliptic equation 
\begin{equation}\label{GSE} 
-Q+\Delta Q+ |x|^{-b}|Q|^{\alpha}Q=0.
\end{equation}
\end{theorem}

\ In this work we prove, for radial initial data, that the global solution obtained in the above theorem also scatters, under some extra restrictions on the parameters $b$ and $\alpha$. These restrictions are probably technical and are direct consequence of the approach used to estimate the nonlinear part $|x|^{-b} |u|^\alpha u$ (see Lemmas \ref{LG1}, \ref{LG2}, \ref{LG3} and Proposition \ref{ECS}). The method of the proof is based on the concentration-compactness and rigidity technique developed by Kenig-Merle \cite{KENIG} and Holmer-Roudenko \cite{HOLROU} (see also Fang-Xie-Cazenave \cite{JIANCAZENAVE} and Guevara \cite{GUEVARA}) for the NLS equation. Our main theorem reads as follows.
\begin{theorem}\label{SCATTERING} Let
\begin{equation}\label{def2_*}
2_*:=\left\{\begin{array}{cl}
\frac{4-2b}{N-2}&\;\;\;\;N\geq 4,\\
3-2b&\;\;\;\;N= 3,\\
\infty &\;\;\;\;\;N=2.
\end{array}\right.
\end{equation} 
Assume that $\frac{4-2b}{N}<\alpha<2_*$ and $0<b<\min\{\frac{N}{3},1\}$. If $u_0\in H^1(\mathbb{R}^N)$ is radial and \eqref{EMC}-\eqref{GFC} are satisfied, then the corresponding solution $u(t)$ of \eqref{INLS} is global in $H^1(\mathbb{R}^N)$ and scatters both forward and backward in time. 
\end{theorem} 

\ For $N=3$, we impose an extra assumption namely $\frac{4-2b}{3}<\alpha<3-2b$. So in 3D it is still an open problem to prove scattering for the global solutions given by Theorem \ref{TG}, when $3-2b\leq \alpha<4-2b$. However, the cubic INLS equation in 3D ($\alpha=2$ and $N=3$) is included in the assumptions of Theorem \ref{SCATTERING} and this is exactly the case considered by the authors in \cite{paper2}. So, the present article can be viewed as an extension of this study to all spacial dimensions $N\geq 2$. In particular, when $N=2$ or $N\geq 4$ the above theorem asserts scattering for all range of $L^2$-supercritical and $H^1$-subcritical INLS equations \eqref{INLS} (recall \eqref{def2*}), assuming that the initial data is radial and satisfies the assumptions \eqref{EMC}-\eqref{GFC}.

\ Similarly as in the NLS model (also $3D$ cubic INLS), to establish scattering we use the following criteria\footnote{The proof will be given after  Proposition \ref{GWPH1}.}.

\begin{proposition}\label{SCATTERSH1}{\bf ($H^1$ scattering)} Let $u(t)$ be  a global solution of \eqref{INLS} with initial data $u_0 \in  H^1(\mathbb{R}^N)$. If $\|u\|_{S(\dot{H}^{s_c})}< +\infty$ and $\sup\limits_{t\in \mathbb{R}}\|u(t)\|_{H^1_x}\leq B$. Then $u(t)$ scatters in $H^1(\mathbb{R}^N)$ as $t \rightarrow \pm\infty$.
\end{proposition}

Our paper is organized as follows. In the next section we introduce some notation and estimates. In Section $3$, we outline the proof of our main result (Theorem \ref{SCATTERING}), assuming all the technical points. In Section $4$, we recall some properties of the ground state and we collect many preliminary results of the Cauchy problem \eqref{INLS}. Next in Section $5$, we establish a profile decomposition result and an
Energy Pythagorean expansion for such a decomposition.  In Section $6$, we construct a critical solution denoted by $u_c$ and we show some of its properties. Finally, Section $7$ is devoted to a rigidity theorem. 
\section{Notation and preliminares}

\ In this section, we introduce some general notations and give basic results that will be used along the work.
\subsection{Some notation}
Given a set $A\subset \mathbb{R}^N$ then $A^C=\mathbb{R}^N \backslash A$ denotes the complement of $A$. We use $c$ to denote various constants that may vary line by line. 
Given any positive numbers $a$ and $b$, the notation $a\lesssim b$ means that there exists a positive constant $c$ that $a\leq cb$. $C_{p,q}$ denotes a constant depending on $p$ and $q$. Given $x,y \in \mathbb{R}^N$ then $x\cdot y$ denotes the inner product of $x$ and $y$ on $\mathbb{R}^N$. We denote by $a^{\pm}=a\pm \varepsilon$ with $\varepsilon>0$ small enough.

\ For $s\in \mathbb{R}$, $J^s$ and $D^s$ denote the Bessel and the Riesz potentials of order $s$, given via Fourier transform by the formulas
$$
\widehat{J^s f}=(1+|y|^2)^{\frac{s}{2}}\widehat{f}\;\;\;\textnormal{and} \;\;\;\;\widehat{D^sf}=|y|^s\widehat{f},
$$
where the Fourier transform of $f(x)$ is given by
$$
\widehat{f}(y)=\int_{\mathbb{R}^N}e^{ix.y}f(x)dx.
$$
We also denote the support of a function $f$, by 
$$supp(f)=\overline{\{f:\mathbb{R}^N\rightarrow \mathbb{C}: f(x)\neq 0\}}.
$$ 
\subsection{Functional spaces} We start with $C^\infty_0(\mathbb{R}^N)$ denoting the space of functions with continuous derivatives of all orders and compact support in $\mathbb{R}^N$. 

\ We use $\|.\|_{L^p}$ to denote the $L^p(\mathbb{R}^N)$ norm with $p\geq 1$. If necessary, we use subscript to inform which variable we are concerned with. The mixed norms in the spaces $L^q_tL^r_x$ and $L^q_TL^r_x$ of $f(x,t)$ are defined, respectively, as 
$$
\|f\|_{L^q_tL^r_x}=\left(\int_{\mathbb{R}}\|f(t,.)\|^q_{L^r_x}dt\right)^{\frac{1}{q}}
$$
and
$$
\|f\|_{L^q_TL^r_x}=\left(\int_T^\infty\|f(t,.)\|^q_{L^r_x}dt\right)^{\frac{1}{q}}
$$
with the usual modifications when $q=\infty$ or $r=\infty$.

\ We also define the norm of the Sobolev spaces $H^{s,r}(\mathbb{R}^N)$ and $\dot{H}^{s,r}(\mathbb{R}^N)$, respectively, by
$$
	\|f\|_{H^{s,r}}:=\|J^sf\|_{L^r}\;\;\;\;\textnormal{and}\;\;\;\;\|f\|_{\dot{H}^{s,r}}:=\|D^sf\|_{L^r}.
$$
If $r=2$ we denote $H^{s,2}=H^s$ and $\dot{H}^{s,2}=\dot{H}^s$.

\ Next we recall some Strichartz norms. We begin with the following definitions:
\begin{definition} The pair $(q,r)$ is called $L^2$-admissible if it satisfies the condition
\begin{equation*}
\frac{2}{q}=\frac{N}{2}-\frac{N}{r},
\end{equation*}
where 
\begin{equation}\label{L2Admissivel}
\left\{\begin{array}{cl}
2\leq & r  \leq \frac{2N}{N-2}\hspace{0.5cm}\textnormal{if}\;\;\;  N\geq 3,\\
2 \leq  & r < +\infty\;  \hspace{0.5cm}\textnormal{if}\;\; \;N=2,\\
2 \leq & r \leq + \infty\;  \hspace{0.5cm}\textnormal{if}\;\;\;N=1.
\end{array}\right.
\end{equation}
\end{definition}				
\begin{remark}
We included in the above definition the improvement, due
to M. Keel and T. Tao \cite{tao:keel}, to the limiting case for Strichartz’s inequalities.
\end{remark}
\begin{definition}\label{Hsdefinition}
 We say the pair $(q,r)$ is $\dot{H}^s$-admissible if\footnote{It is worth mentioning that the pair $\left(\infty,\frac{2N}{N-2s}\right)$ also satisfies the relation \eqref{CPA1}, however, in our work we will not make use of this pair when we estimate the nonlinearity $|x|^{-b} |u|^\alpha u$.}
\begin{equation}\label{CPA1}
\frac{2}{q}=\frac{N}{2}-\frac{N}{r}-s,
\end{equation}
where 
\begin{equation}\label{CPA2}
\left\{\begin{array}{cl}
\frac{2N}{N-2s}< & r  \leq\left(\frac{2N}{N-2}\right)^{-}\;\;\hspace{0.4cm}\textnormal{if}\;\;  N\geq 3,\\
\frac{2}{1-s} <  & r \leq \left((\frac{2}{1-s})^+\right)'\;  \hspace{0.2cm}\textnormal{if}\;\; \;N=2,\\
\frac{2}{1-2s} < & r \leq + \infty\; \; \hspace{1.2cm}\textnormal{if}\;\;\;N=1.
\end{array}\right.
\end{equation}
Moreover, $(a^+)'$ is the number such that 
\begin{equation}\label{a^+}
\frac{1}{a}=\frac{1}{(a^+)'}+\frac{1}{a^+},
\end{equation}
that is $(a^+)':=\frac{a^+.a}{a^+-a}$. Finally we say that $(q,r)$ is $\dot{H}^{-s}$-admissible if 
$$
\frac{2}{q}=\frac{N}{2}-\frac{N}{r}+s,
$$
where
\begin{equation}\label{H-s}
\left\{\begin{array}{cl}
\left(\frac{2N}{N-2s}\right)^{+}\leq & r  \leq\left(\frac{2N}{N-2}\right)^{-}\;\;\hspace{0.4cm}\textnormal{if}\;\;  N\geq 3,\\
\left(\frac{2}{1-s}\right)^{+} \leq  & r \leq \left((\frac{2}{1+s})^+\right)'\;  \hspace{0.2cm}\textnormal{if}\;\; \;N=2,\\
\left(\frac{2}{1-2s}\right)^{+} \leq & r \leq + \infty\; \; \hspace{1.2cm}\textnormal{if}\;\;\;N=1.
\end{array}\right.
\end{equation}
\end{definition}				
\  Given $s\in \mathbb{R}$, let $\mathcal{A}_s=\{(q,r);\; (q,r)\; \textnormal{is} \;\dot{H}^s-\textnormal{admissible}\}$ and $(q',r')$ is such that $\frac{1}{q}+\frac{1}{q'}=1$ and $\frac{1}{r}+\frac{1}{r'}=1$ for $(q,r)\in \mathcal{A}_s$. We define the following Strichartz norm
$$
\|u\|_{S(\dot{H}^{s})}=\sup_{(q,r)\in \mathcal{A}_{s}}\|u\|_{L^q_tL^r_x} 
$$
and the dual Strichartz norm
$$
\|u\|_{S'(\dot{H}^{-s})}=\inf_{(q,r)\in \mathcal{A}_{-s}}\|u\|_{L^{q'}_tL^{r'}_x}.
$$
\begin{remark}
Note that, if $s=0$ then $\mathcal{A}_0$ is the set of all $L^2$-admissible pairs. Moreover, if $s=0$, $S(\dot{H}^0)=S(L^2)$ and $S'(\dot{H}^{0})=S'(L^2)$. We just write $S(\dot{H}^s)$ or $S'(\dot{H}^{-s})$ if the mixed norm is evaluated over $\mathbb{R}\times\mathbb{R}^N$. To indicate a restriction to a time interval $I\subset (-\infty,\infty)$ and a subset $A$ of $\mathbb{R}^N$, we will consider the notations $S(\dot{H}^s(A);I)$ and $S'(\dot{H}^{-s}(A);I)$. 
\end{remark}

\subsection{Basic estimates} We start with two important remarks (the second one provides a condition for the integrability of $|x|^{-b}$ on $B$ and $B^C$).		
\begin{remark}\label{RB} 
Let $B=B(0,1)=\{ x\in \mathbb{R}^N;|x|\leq 1\}$ and $b>0$. If $x\in B^C$ then $|x|^{-b}<1$ and so
$$ 
\left	\||x|^{-b}f \right\|_{L^r_x}\leq \|f\|_{L_x^r(B^C)}+\left\||x|^{-b}f\right\|_{L_x^r(B)}. 
$$
\end{remark}

\begin{remark}\label{RIxb} 
Note that if $\frac{N}{\gamma}-b>0$ then $\||x|^{-b}\|_{L^\gamma(B)}<+\infty$. Indeed 
\begin{equation*}
\int_{B}|x|^{-\gamma b}dx=c\int_{0}^{1}r^{-\gamma b}r^{N-1}dr=c_1 \left. r^{N-\gamma b} \right |_0^1<+\infty\;\;\textnormal{if}\;\;\frac{N}{\gamma} - b>0.
\end{equation*}
Similarly, we have that $\||x|^{-b}\|_{L^\gamma(B^C)}$ is finite if $\frac{N}{\gamma}- b<0$.
\end{remark}

\ Now, we list (without proving) some well known estimates associated to the linear Schr\"odinger operator.
  
\begin{lemma}\label{ILE}  
If $t\neq 0$, $\frac{1}{p}+\frac{1}{p'}=1$ and $p'\in[1,2]$, then $U(t):L^{p'}(\mathbb{R}^N)\rightarrow L^p(\mathbb{R}^N)$ is continuous and 
$$
\|U(t)f\|_{L^p_x}\lesssim|t|^{-\frac{N}{2}(\frac{1}{p'}-\frac{1}{p})}\|f\|_{L^{p'}}.
$$
\begin{proof} See Linares-Ponce \cite[Lemma $4.1$]{FELGUS}.
\end{proof}
\end{lemma} 
\begin{lemma}\label{Lemma-Str} The following statements hold.
 \begin{itemize}
\item [(i)] (Linear estimates).
\begin{equation}\label{SE1}
\| U(t)f \|_{S(L^2)} \leq c\|f\|_{L^2},
\end{equation}
\begin{equation}\label{SE2}
\|  U(t)f \|_{S(\dot{H}^s)} \leq c\|f\|_{\dot{H}^s}.
\end{equation}
\item[(ii)] (Inhomogeneous estimates).
\begin{equation}\label{SE3}					 
\left \| \int_{\mathbb{R}} U(t-t')g(.,t') dt' \right\|_{S(L^2)}\;+\; \left \| \int_{0}^t U(t-t')g(.,t') dt' \right \|_{S(L^2) } \leq c\|g\|_{S'(L^2)},
\end{equation}
\begin{equation}\label{SE5}
\left \| \int_{0}^t U(t-t')g(.,t') dt' \right \|_{S(\dot{H}^s) } \leq c\|g\|_{S'(\dot{H}^{-s})}.
\end{equation}
\end{itemize}
\end{lemma} 
\noindent The inequalities of Lemma \ref{Lemma-Str} are the well known Strichartz estimates. The relations \eqref{SE3} and \eqref{SE5} will be very useful to perform estimates on the nonlinearity $|x|^{-b}|u|^\alpha u$. We refer the reader to Linares-Ponce \cite{FELGUS} and Kato \cite{KATO} (see also Holmer-Roudenko \cite{HOLROU} and Guevara \cite{GUEVARA}).

\ We end this section by recalling the Sobolev inequalities and giving a useful remark.
\begin{lemma}\label{SI} Let $s\in (0,+\infty)$ and $1\leq p<+\infty$.
\begin{itemize}
\item [(i)] If $s\in \left(0,\frac{N}{p}\right)$ then $H^{s,p}(\mathbb{R}^N)$ is continuously embedded in $L^r(\mathbb{R}^N)$ where $s=\frac{N}{p}-\frac{N}{r}$. Moreover, 
\begin{equation}\label{SEI} 
\|f\|_{L^r}\leq c\|D^sf\|_{L^{p}}.
\end{equation}
\item [(ii)] If $s=\frac{N}{2}$ then $H^{s}(\mathbb{R}^N)\subset L^r(\mathbb{R}^N)$ for all $r\in[2,+\infty)$. Furthermore,
\begin{equation}\label{SEI1} 
\|f\|_{L^r}\leq c\|f\|_{H^{s}}.
\end{equation}
\end{itemize}
\begin{proof} See Bergh-L\"ofstr\"om \cite[Theorem $6.5.1$]{BERLOF} (see also Linares-Ponce \cite[Theorem $3.3$]{FELGUS} and Demenguel-Demenguel \cite[Proposition 4.18]{DEMENGEL}). 
\end{proof}
\end{lemma}
As a consequence of Lemma \ref{SI} $(i)$ (particular case: $p=2$ and $s\in (0,\frac{N}{2})$) we have that $H^s(\mathbb{R}^N)$ is continuously embedded in $L^r(\mathbb{R}^N)$ and 
\begin{equation}\label{SEI22} 
\|f\|_{L^r}\leq c\|f\|_{H^{s}},
\end{equation}    
where $r\in[2,\frac{2N}{N-2s}]$.
\begin{remark}\label{nonlinerity}
Let $F(x,z)=|x|^{-b}|z|^\alpha z$, and $f(z)=|z|^\alpha z$. The complex derivative of $f$ is
\begin{equation*}\label{nonli1}
f_z(z)=\frac{\alpha+2}{2}|z|^\alpha\;\;\;\;\;\textnormal{and}\;\;\;\; f_{\bar{z}}    (z)=\frac{\alpha}{2}|z|^{\alpha-2}z^2. 
\end{equation*}
For $z,w\in \mathbb{C}$, we get 
\begin{equation*}
 f(z)-f(w)=\int_{0}^{1}\left[f_z(w+\theta(z-w))(z-w)+f_{\bar{z}}(w+\theta(z-w))\overline{(z-w)}\right]d\theta.
\end{equation*}
Hence,
\begin{equation}\label{FEI}
 |F(x,z)-F(x,w)|\lesssim |x|^{-b}\left( |z|^\alpha+ |w|^\alpha \right)|z-w|.
\end{equation}

\ Our interest now is to estimate $\nabla \left( F(x,z)-F(x,w) \right)$. A simple computation gives  
\begin{equation}\label{NONLI11}
\nabla F(x,z)=\nabla(|x|^{-b})f(z)+|x|^{-b} \nabla f(z)
\end{equation}
where $f(z)=f'(z)\nabla z=f_z(z)\nabla z+f_{\bar{z}}(z) \overline{\nabla z}$.\\
We first estimate $|\nabla (f(z)-f(w))|$. Observe that
 \begin{equation}\label{NONLI55}
 \nabla (f(z)-f(w))=f'(z)(\nabla z-\nabla w)+(f'(z)-f'(w))\nabla w.
 \end{equation}
So, since (the proof of the following estimate can be found in Cazenave-Fang-Han \cite[Remark $2.3$]{CAZENAVECONTINUOUS})
\begin{equation*}\label{NONLI66}
|f_z(z)-f_z(w)|\lesssim \left\{\begin{array}{cl}
(|z|^{\alpha-1}+|w|^{\alpha-1})|z-w| &\textnormal{if}\;\alpha> 1,\\
|z-w|^\alpha&\textnormal{if}\;0<\alpha\leq 1
\end{array}\right.     
\end{equation*}
and 
\begin{equation*}\label{NONLI77}
| f_{\bar{z}}(z)-f_{\bar{z}}(w)|\lesssim \left\{\begin{array}{cl}
(|z|^{\alpha-1}+|w|^{\alpha-1})|z-w| & \textnormal{if}\;\alpha> 1,\\
 |z-w|^\alpha& \textnormal{if}\;0<\alpha\leq 1,
\end{array}\right.     
\end{equation*}
 we get that
\begin{equation*}\label{NONLI88}
|\nabla (f(z)-f(w))|\lesssim |z|^\alpha|\nabla (z- w)|+(|z|^{\alpha-1}+|w|^{\alpha-1})|\nabla w||z-w|\;\;\textnormal{if}\;\alpha > 1
\end{equation*}
and  
\begin{equation*}\label{NONLI99}
|\nabla (f(z)-f(w))|\lesssim |z|^\alpha|\nabla (z- w)|+|z-w|^\alpha|\nabla w|\;\;\;\textnormal{if}\;\;0<\alpha \leq 1,
\end{equation*}
where we have used \eqref{NONLI55}. Therefore, by \eqref{NONLI11}, \eqref{FEI} and the two last inequalities we obtain \vspace{0.1cm}
\begin{equation}\label{SECONDEI}
\left|\nabla \left(F(x,z)-F(x,w)\right)\right|\lesssim  |x|^{-b-1}(|z|^{\alpha}+|w|^{\alpha})|z-w|+|x|^{-b}|z|^\alpha|\nabla (z- w)|+E, 
\end{equation}
where 
\begin{eqnarray*}\label{NONLI100}
 E &\lesssim& \left\{\begin{array}{cl}
 |x|^{-b}\left(|z|^{\alpha-1}+|w|^{\alpha-1}\right)|\nabla w||z-w| & \textnormal{if}\;\;\;\alpha > 1 \vspace{0.2cm} \\ 
|x|^{-b}|\nabla w||z-w|^{\alpha} & \textnormal{if}\;\;\;0<\alpha\leq 1.
\end{array}\right.
\end{eqnarray*}
 \end{remark}

	
	
			
	

\section{Outline of the proof of Theorem \ref{SCATTERING}}\label{SPMR}


\ In this short section, we give the proof of Theorem \ref{SCATTERING}, assuming all preliminary results. We start with the following definition.
\begin{definition}
We shall say that SC($u_0$) holds if the solution $u(t)$ with initial data $u_0\in H^1(\mathbb{R}^N)$ is global and \eqref{HsFINITE} holds.
\end{definition}

\ Let $u(t)$ be the corresponding $H^1$ solution for the IVP \eqref{INLS} with radial data $u_0\in H^1(\mathbb{R}^N)$ satisfying \eqref{EMC} and \eqref{GFC}. We already know by Theorem \ref{TG} that the solution is globally defined and $\sup\limits_{t\in \mathbb{R}}\|u(t)\|_{H^1}< \infty$. Furthermore, if 
\begin{equation}\label{HsFINITE} 
\|u\|_{S(\dot{H}^{s_c})}<+\infty
\end{equation}
then $u$ scatters in $H^1(\mathbb{R}^N)$ ( see Proposition \ref{SCATTERSH1}). To achieve the scattering property \eqref{HsFINITE}, we follow the exposition in Holmer-Roudenko \cite{HOLROU}
and Fang-Xie-Cazenave \cite{JIANCAZENAVE} (see also our work \cite{paper2}),
which was based in the ideas introduced by Kenig-Merle \cite{KENIG}. Indeed, define
\begin{definition}
For each $\delta > 0$ define the set $A_\delta$ to be the collection of all initial data in $H^1(\mathbb{R}^N)$ satisfying
\begin{align*}
A_\delta=\{u_0\in H^1:E[u_0]^{s_c}M[u_0]^{1-s_c}<\delta\;\textnormal{and}\;\|\nabla u_0\|^{s_c}_{L^2}\| u_0\|^{1-s_c}_{L^2}<\|\nabla Q\|^{s_c}_{L^2}\| Q\|^{1-s_c}_{L^2} \}
\end{align*}
and define
\begin{equation}\label{deltac}
\delta_c=\sup \{\; \delta>0:\; u_0\; \in A_\delta\;  \Longrightarrow SC(u_0)\; \textnormal{holds} \}=\sup_{\delta>0} B_\delta.
\end{equation}
\end{definition}
\noindent Note that there always exists a $\delta>0$ such that the above statement is true, i.e., $B_\delta \neq \emptyset$ (see the proof at the end of this section). 
 
 \ To prove Theorem \ref{SCATTERING} we have two cases to consider. If $\delta_c\geq E[Q]^{s_c}M[Q]^{1-s_c}$ then we are done. Assume now, by contradiction, that $\delta_c<E[Q]^{s_c}M[Q]^{1-s_c}$. There exists a sequence of radial solutions $u_n$ to \eqref{INLS} with $H^1$ initial data $u_{n,0}$ (rescale all of them to have $\|u_{n,0}\|_{L^2} = 1$ for all $n$) such that\footnote{We can rescale $u_{n,0}$ such that $\|u_{n,0}\|_{L^2} = 1$. Indeed, if $u^\lambda_{n,0}(x)=\lambda^{\frac{2-b}{\alpha}}u_{n,0}(\lambda x)$ then by \eqref{QI2} we have $E[u^\lambda_{n,0}]^{s_c}M[u^\lambda_{n,0}]^{1-s_c}<E[Q]^{s_c}M[Q]^{1-s_c}$ and $\|\nabla u^\lambda_{n,0}\|^{s_c}_{L^2}\| u^\lambda_{n,0}\|^{1-s_c}_{L^2}<\|\nabla Q\|^{s_c}_{L^2}\| Q\|^{1-s_c}_{L^2} $. Moreover, since $\|u^\lambda_{n,0}\|_{L^2} = \lambda^{-s_c}\|u_{n,0}\|_{L^2}$ by \eqref{QI1}, setting $\lambda^{s_c}=\|u_{n,0}\|_{L^2}$ we have $\|u^\lambda_{n,0}\|_{L^2} = 1$.}   
\begin{equation}\label{CC0}   
\|\nabla u_{n,0}\|^{s_c}_{L^2} <  \|\nabla Q\|^{s_c}_{L^2}\|Q\|^{1-s_c}_{L^2}
\end{equation}
and
\begin{equation*}
E[u_n]^{s_c} \searrow \delta_c\; \textnormal{as}\; n \rightarrow +\infty,
\end{equation*} 	
for which SC($u_{n,0}$) does not hold for any $n\in \mathbb{N}$, that is $\|u_n\|_{S(\dot{H}^{s_c})}=+\infty$, since we get by Theorem \ref{TG} that $u_n$ is globally defined. Thus using a profile decomposition result (see Proposition \ref{LPD}) on the sequence $\{u_{n,0}\}_{n\in \mathbb{N}}$ we can construct a critical solution of \eqref{INLS}, denoted by $u_c$, that lies exactly at the threshold $\delta_c$, satisfies \eqref{CC0} (it implies that $u_c$ is globally defined again by Theorem \ref{TG}) and $\|u_c\|_{S(\dot{H}^{s_c})}=+\infty$ (see Proposition \ref{ECS}). Moreover, we show that the critical solution $u_c$ has the property that $K=\{u_c(t):t\in[0,+\infty)\}$ is precompact in $H^1(\mathbb{R}^N)$ (see Proposition \ref{PSC}). Finally, the rigidity theorem (Theorem \ref{RT}) will imply that $u_c$ (critical solution) is identically zero, which contradicts the fact that $\|u_c\|_{S(\dot{H}^{s_c})}=+\infty$. 

\ To complete the proof it remains to establish $B_\delta \neq \emptyset$. Indeed, the Strichartz estimate \eqref{SE2}, interpolation and Lemma \ref{LGS} (i) imply that
  \begin{eqnarray*}
  \|U(t)u_0\|_{S(\dot{H}^{s_c})}&\leq& c\|u_0\|_{\dot{H}^{s_c}}\leq c\|\nabla u_0\|^{s_c}_{L^2}\|u_0\|^{1-s_c}_{L^2}\\
  &\leq& c\left(\frac{N\alpha+2b}{\alpha s_c}\right)^{\frac{s_c}{2}}E[u_0]^{\frac{s_c}{2}}M[u_0]^{\frac{1-s_c}{2}}.
\end{eqnarray*}
So if $u_0\in A_\delta$ then $
E[u_0]^{s_c}M[u_0]^{1-s_c}<\left(\frac{\alpha s_c}{N\alpha+2b}\right)^{s_c}\delta'^2,
$ which implies $\|U(t)u_0\|_{S(\dot{H}^{s_c})}\leq c\delta'$.
Therefore, by the small data theory (Proposition \ref{GWPH1}) we obtain that $SC(u_0)$ holds for $\delta'>0$ small enough.

\section{Energy bounds and the Cauchy problem}

\ We divide this section  in two parts. First, we recall some properties that are related to our problem and we provide important estimates. Subsequently, we show the basic results concerning the IVP \eqref{INLS} that will help us in the proof of Theorem \ref{SCATTERING}.

\ We start with the following Gagliardo-Nirenberg inequality (it was obtained by the first
author in \cite{LG})
\begin{equation}\label{GNI} 
\left\||x|^{-b}|u|^{\alpha+2} \right\|_{L^1_x}\leq C_{GN}\|\nabla u\|^{\frac{N\alpha+2b}{2}}_{L^2_x}\|u\|^{\frac{4-2b-\alpha(N-2)}{2}}_{L^2_x},
\end{equation}
with the sharp constant 
\begin{equation}\label{GNI1}
C_{GN}=\frac{2(\alpha +2)}{N\alpha +2b}\left(\frac{4-2b-\alpha(N-2)}{N\alpha +2b}\right)^{\alpha s_c/2}\frac{1}{\|Q\|^{\alpha}_{L^2}}
\end{equation}
where $Q$ is the ground state solution of \eqref{GSE}. Furthermore, $Q$ satisfies  

\begin{equation}\label{GS1}
\|\nabla Q\|^2_{L^2}=\frac{N\alpha+2b}{4-2b-\alpha(N-2)}\|Q\|^2_{L^2}
\end{equation}
and 
\begin{equation}\label{GS2}
\left\||x|^{-b}|Q|^{\alpha+2} \right\|_{L^1}=\frac{2(\alpha+2)}{N\alpha+2b}\|\nabla Q\|^2_{L^2}.
\end{equation}
Combining the relations \eqref{GNI1}, \eqref{GS1} and \eqref{GS2} we deduce (recalling $s_c=\frac{N}{2}-\frac{2-b}{\alpha}$)
\begin{equation}\label{GNI2}
C_{GN}=\frac{2(\alpha+2)}{(N\alpha+2b)\|\nabla Q\|^{\alpha s_c}_{L^2}\|Q\|^{\alpha(1-s_c)}_{L^2}}.
\end{equation}
Also, we get 
\begin{equation}\label{EGS} 
E[Q]=\frac{1}{2}\|\nabla Q\|^2_{L^2}-\frac{1}{\alpha+2}\left\||x|^{-b}|Q|^{\alpha+2}\right\|_{L^1}=\frac{\alpha s_c}{N\alpha+2b}\|\nabla Q\|^2_{L^2}.
\end{equation}

\ In the sequel, we show the radial Sobolev Gagliardo-Nirenberg inequality in $N$ dimension. The proof follows the ideas introduced by Strauss \cite{Strauss}.
\begin{lemma}\label{RSGN} Let $N\geq 2$, $R>0$ and $f\in H^1(\mathbb{R}^N)$ a radial function. Then the following inequality holds
\begin{equation}\label{RSGN1}
\sup_{|x|\geq R}|f(x)|\leq \frac{1}{R^{\frac{N-1}{2}}}\|f\|^{\frac{1}{2}}_{L^2}\|\nabla f\|^{\frac{1}{2}}_{L^2}.
\end{equation}
\begin{proof}
Since $f$ is radial we deduce
\begin{eqnarray*}
\sup_{|x|\geq R}|f(x)|^2&=&\sup_{|x|\geq R}\frac{1}{2}\int_{|x|}^{+\infty} \partial_r(f^2)dr\\
&\leq & \int_{R}^{+\infty} f \partial_rfdr\\
&\leq & \left(\int_{R}^{+\infty} |f|^2dr \right)^{\frac{1}{2}}\left(\int_{R}^{+\infty} |\partial_rf|^2dr \right)^{\frac{1}{2}},
\end{eqnarray*}
where we have used that $f$ has to vanish at infinite and the Cauchy-Schwarz inequality.  On the other hand, the fact that $|x|\geq R$ (or $r\geq R$) implies $1 \leq \frac{r}{R}$ so
\begin{eqnarray*}
\sup_{|x|\geq R}|f(x)|^2&\leq & \left(\int_{R}^{+\infty} |f|^2 \left(\frac{r}{R}\right)^{N-1} \right)^{\frac{1}{2}}\left(\int_{R}^{+\infty} |\partial_rf|^2\left(\frac{r}{R}\right)^{N-1}dr \right)^{\frac{1}{2}}\\
&\leq & \frac{1}{R^{\frac{N-1}{2}}} \left(\int_{R}^{+\infty} |f|^2 r^{2(N-1)} \right)^{\frac{1}{2}}\frac{1}{R^{\frac{N-1}{2}}}\left(\int_{R}^{+\infty} |\partial_rf|^2 r^{2(N-1)}dr \right)^{\frac{1}{2}}\\
&=&  \frac{1}{R^{N-1}} \left(\int_{R}^{+\infty} |f|^2 dx \right)^{\frac{1}{2}}\left(\int_{R}^{+\infty} |\nabla f|^2 dx \right)^{\frac{1}{2}}\\
&\leq& \frac{1}{R^{N-1}}\|f\|_{L^2}\|\nabla f\|_{L^2},
\end{eqnarray*}
where in the third line we have used the fact that $|\partial_r f|=|\nabla f|$ for radial functions. We finish the proof taking the square root on both sides.
\end{proof}
\end{lemma}

\indent We now provide some useful energy inequalities.   

\begin{lemma}\label{LGS} 
Let $v \in H^1(\mathbb{R}^N)$ such that 
\begin{equation}\label{LGS0}
E[v]^{s_c}M[v]^{1-s_c}<E[Q]^{s_c}M[Q]^{1-s_c}
\end{equation}
and
\begin{equation}\label{LGS1}
\|\nabla v\|^{s_c}_{L^2}\|v\|_{L^2}^{1-s_c}\leq \|\nabla Q\|^{s_c}_{L^2}\|Q\|_{L^2}^{1-s_c}.
 \end{equation}	
 Then, the following statements hold
\begin{itemize}
\item [(i)] $\frac{\alpha s_c}{N\alpha +2b}\|\nabla v\|^2_{L^2}\leq E(v)\leq \frac{1}{2} \|\nabla v\|^{2}_{L^2}$,
\item [(ii)] $\|\nabla v\|^{s_c}_{L^2}\|v\|^{1-s_c}_{L^2}\leq w^{\frac{1}{2}} \|\nabla Q\|^{s_c}_{L^2}\|Q\|^{1-s_c}_{L^2}$, 
\item [(iii)] $16A E[v]\leq 8A \|\nabla v\|_{L^2}^2\leq 8 \|\nabla v\|^2_{L^2}-\frac{4(N\alpha+2b)}{\alpha+2}\left\||x|^{-b}|v|^{\alpha+2}\right\|_{L^1}$,
\end{itemize}
where\footnote{Note that, the relation \eqref{LGS0} implies that $w<1$ and $A>0$.} $w=\frac{E[v]^{s_c}M[v]^{1-s_c}}{E[Q]^{s_c}M[Q]^{1-s_c}}$ and $A=(1-w^{\frac{\alpha}{2}})$.
\begin{proof}
 (i) The definition of Energy \eqref{energy} yields the second inequality. The first one is obtained by observing that (using \eqref{GNI}, \eqref{GNI2} and \eqref{LGS1})
\begin{eqnarray*}
E[v] &\geq& \frac{1}{2}\|\nabla v\|^2_{L^2} -  \frac{C_{GN}}{\alpha+2}\|\nabla v\|^{\frac{N\alpha+2b}{2}}_{L^2}\|v\|^{\frac{4-2b-\alpha(N-2)}{2}}_{L^2}\\
&=&\frac{1}{2}\|\nabla v\|^2_{L^2}\left(1- \frac{2C_{GN}}{\alpha+2} \|\nabla v\|^{\alpha s_c}_{L^2}\|v\|^{\alpha(1-s_c)}_{L^2}  \right)\\
&\geq& \frac{1}{2}\|\nabla v\|^2_{L^2}\left(1- \frac{2C_{GN}}{\alpha+2} \|\nabla Q\|^{\alpha s_c}_{L^2}\|Q\|^{\alpha(1-s_c)}_{L^2}  \right)\\
&=&\frac{N\alpha -(4-2b)}{2(N\alpha+2b)}\|\nabla v\|^2_{L^2}\\
&=&\frac{\alpha s_c}{N\alpha+2b}\|\nabla v\|^2_{L^2}.
\end{eqnarray*}	  
	
 (ii) The first inequality in (i) implies $\|\nabla v\|^2_{L^2}\leq \frac{N\alpha +2b}{\alpha s_c}E(v)$, multiplying it by $M[v]^\sigma=\|v\|_{L^2}^{2\sigma}$, where $\sigma=\frac{1-s_c}{s_c}$, we get
\begin{eqnarray*}
\|\nabla v\|^2_{L^2}\| v\|^{2\sigma}_{L^2}&\leq& \frac{N\alpha +2b}{\alpha s_c}E[v]M[v]^\sigma\\
 & = & \frac{N\alpha +2b}{\alpha s_c}\frac{E[v]M[v]^\sigma}{E[Q]M[Q]^\sigma} E[Q]M[Q]^\sigma\\
&=& w \|\nabla Q\|^2\|Q\|^{2\sigma}_{L^2},
\end{eqnarray*}
where we have used \eqref{EGS}.
		
(iii) Let $B=8 \|\nabla v\|^2_{L^2}-\frac{4(N\alpha+2b)}{\alpha+2}\left\||x|^{-b}|v|^{\alpha+2}\right\|_{L^1}$. Applying the Gagliardo-Niremberg inequality \eqref{GNI} and item (ii) we have  
\begin{eqnarray*}
B&\geq& 8\|\nabla v\|^2_{L^2}- \frac{4(N\alpha+2b)C_{GN}}{\alpha+2}\|\nabla v\|^{\frac{N\alpha+2b}{2}}_{L^2}\|v\|^{\frac{4-2b-\alpha(N-2)}{2}}_{L^2}   \\      
&\geq& \|\nabla v\|^2_{L^2}\left(8- \frac{4(N\alpha+2b)}{\alpha+2}C_{GN}w^\frac{\alpha}{2} \|\nabla Q\|^{\alpha s_c}_{L^2}\|Q\|^{\alpha(1-s_c)}_{L^2}  \right)\\
&=& \|\nabla v\|^2_{L^2}8(1-w^{\frac{\alpha}{2}}),
\end{eqnarray*} 
where in the last equality, we have used \eqref{GNI2}. Next, the first inequality obviously holds.	
\end{proof}	
\end{lemma}

\ We end the first part of this section with a important lemma. Define
 \begin{equation}\label{defr*}
\widetilde{p}=\frac{2N}{N-2}\;\; \textnormal{if}\;\; N\geq 3 \;\;\textnormal{and}\;\;\widetilde{p}=\infty\;\; \textnormal{if} \;\;N=2.
\end{equation}
\begin{lemma}\label{LEWO}
Let $\frac{4-2b}{N}<\alpha<2^*$ and $0<b<\min\{2,N\}$. If $f$ and $g\in H^1(\mathbb{R}^N)$ then 
 \begin{itemize}
\item [(i)] $\left\|  |x|^{-b} |f|^{\alpha+1}g \right\|_{L^1} 
 \leq c \|f\|^{\alpha+1}_{L^{\alpha+2}} \|g\|_{L^{\alpha+2}} + c\|f\|^{\alpha+1}_{L^{r}}\|g\|_{L^{r}}$
\item [(ii)] $\left\|  |x|^{-b} |f|^{\alpha+1}g \right\|_{L^1} 
 \leq c \|f\|^{\alpha+1}_{{H^1}} \|g\|_{H^1}$
\item [(iii)] $ \lim\limits_{|t|\rightarrow +\infty}  \left\|  |x|^{-b} |U(t) f|^{\alpha+1}g \right\|_{L^1_x}=0.$ 
 \end{itemize}
where\footnote{Since $\alpha<\widetilde{p}$ we have $\frac{N(\alpha+2)}{N-b} <\widetilde{p}$, which implies that there exists $r$ such that $\frac{N(\alpha+2)}{N-b} <r<\widetilde{p}$ . } $2<\frac{N(\alpha+2)}{N-b}<r<\widetilde{p}$.
\begin{proof} (i) We divide the estimate in the regions $B^C$ and $B$. Indeed, it follows from Remark \ref{RB} and the H\"older inequality (since $1=\frac{\alpha+1}{\alpha+2}+\frac{1}{\alpha+2}$) that
\begin{eqnarray}\label{LEWO1}
\left\|  |x|^{-b} |f|^{\alpha+1}g \right\|_{L^1} 
&\leq &\left\|   |f|^{\alpha+1}g \right\|_{L^1(B^C)} 
+\left\|  |x|^{-b} |f|^{\alpha+1}g \right\|_{L^1(B)} \nonumber\\
&\leq& \|f\|^{\alpha+1}_{L^{\alpha+2}}\|g\|_{L^{\alpha+2}}+\|x|^{-b}|\|_{L^\gamma(B)}\|f\|^{\alpha+1}_{L^{(\alpha+1)\beta}}\|g\|_{L^r} \nonumber\\
 &=& \|f\|^{\alpha+1}_{L^{\alpha+2}}\|g\|_{L^{\alpha+2}}+\|x|^{-b}|\|_{L^\gamma(B)} \|f\|^{\alpha+1}_{L^{r}} \|g\|_{L^r},
\end{eqnarray}
where
 \begin{equation}\label{LEWO2}
1=\frac{1}{\gamma}+\frac{1}{\beta}+\frac{1}{r}\;\;\;\;\;\textnormal{and}\;\;\;\;\;\;r=(\alpha+1)\beta.
\end{equation}
Since $r>\frac{N(\alpha+2)}{N-b}$ and using \eqref{LEWO2}, we obtain
$$
\frac{N}{\gamma}=N-\frac{N(\alpha+2)}{r}>b,
$$
which implies that $\||x|^{-b}\|_{L^\gamma(B)}$ is bounded (see Remark \ref{RIxb}). This completes the proof of item (i).

\ (ii) Applying the Sobolev inequality \eqref{SEI1} (for $N=2$ and $s=1$) and \eqref{SEI22} (for $N\geq 3$ and $s=1$), it is easy to see that $H^1 \hookrightarrow  L^{\alpha+2} $ and $H^1 \hookrightarrow L^{r}$ (where $2<\frac{N(\alpha+2)}{N-b}<r<\widetilde{p}$), then by \eqref{LEWO1} we get (ii).  

\ (iii) We also have (using the same argument as (i) and (ii))
 \begin{eqnarray}\label{lematecnico}
 \left\|  |x|^{-b} |U(t) f|^{\alpha+1}g \right\|_{L^1_x}\leq c\|U(t)f\|^{\alpha+1}_{L^{\alpha+2}} \|g\|_{H^1}+ c\|U(t)f\|^{\alpha+1}_{L^{r}}\|g\|_{H^1},
\end{eqnarray}
for $2<\frac{N(\alpha+2)}{N-b}<r<\widetilde{p}$.\\
To complete the proof we show that $\|U(t)f\|_{L^{r}_x}$ and $\|U(t)f\|_{L^{\alpha+2}_x}$ $\rightarrow 0$ as 
$|t|\rightarrow +\infty$. Indeed, it suffices to show  (since $r$ and $\alpha+2$ belong to $(2,\widetilde{p})$)
\begin{equation}\label{LEWO3}
\lim\limits_{|t|\rightarrow +\infty}\|U(t)f\|_{L^p_x}=0, 
\end{equation}
where $2<p<\widetilde{p}$. Let $\widetilde{f}\in H^1\cap L^{p'}$, the Sobolev embedding \eqref{SEI1} if $N=2$ or \eqref{SEI22} if $N\geq 3$ and Lemma \ref{ILE} yield
 $$
 \|U(t)f\|_{L^{p}_x}\leq c \|f-\widetilde{f}\|_{{H}^1}+c|t|^{-\frac{N(p-2)}{2p}}\|\widetilde{f}\|_{L^{p'}}.
 $$
Note that the exponent of $|t|$ is negative (since $p>2$), then approximating $f$ by $\widetilde{f}\in C^\infty_0$ in $H^1$, we obtain \eqref{LEWO3}.
\end{proof}
\end{lemma}

\ Our interested now is to show a miscellaneous of results for the Cauchy problem \eqref{INLS}. We begin by recalling the small data global theory in $H^1$ (it was obtained by the second author in \cite{CARLOS}). After that, we prove the $H^1$
-scattering criterion, the perturbation theory  and the existence of wave operators. To this end, the heart of the proof is to establish good estimates on the nonlinearity. The next lemmas provide these estimates. 
\begin{lemma}\label{LG1} 
Let $N\geq 2$, $ \frac{4-2b}{N}<\alpha<2^*$ and $0<b<\min\{\frac{N}{3},1\}$. Then there exist $c>0$ and $\theta\in (0,\alpha)$ sufficiently small such that 
\begin{equation}\label{ELG1} 
\left \||x|^{-b}|u|^\alpha v \right\|_{S'(\dot{H}^{-s_c})} \lesssim \| u\|^{\theta}_{L^\infty_tH^1_x}\|u\|^{\alpha-\theta}_{S(\dot{H}^{s_c})} \|v\|_{S(\dot{H}^{s_c})}.
\end{equation}
\begin{proof} See \cite[Lemma $4.1$, with $s=1$]{CARLOS}.
\end{proof} 
\end{lemma} 
 
\begin{lemma}\label{LG2} 
Let $N\geq 2$, $ \frac{4-2b}{N}<\alpha<2^*$ and $0<b<\min\{\frac{N}{3},1\}$. Then there exist $c>0$ and $\theta\in (0,\alpha)$ sufficiently small such that 
\begin{equation}\label{ELG2} 
\left\||x|^{-b}|u|^\alpha v \right\|_{S'(L^2)}\lesssim \| u\|^{\theta}_{L^\infty_tH^1_x}\|u\|^{\alpha-\theta}_{S(\dot{H}^{s_c})} \| v\|_{S(L^2)}.
\end{equation}
\begin{proof} See \cite[Lemma $4.2$, with $s=1$]{CARLOS}.
\end{proof}
\end{lemma} 
\begin{remark}\label{RGP} As an immediate consequence of Lemma \ref{LG2}, we have
following estimate for $\alpha>1$
$$
\left \||x|^{-b}|u|^{\alpha-1} v  w\right\|_{S'(L^2)} \leq c \| u\|^{\theta}_{L^\infty_tH^1_x}\|u\|^{\alpha-1-\theta}_{S(\dot{H}^{s_c})} \|v\|_{S(\dot{H}^{s_c})}\|  w\|_{S(L^2)},
$$
where $\theta\in(0,\alpha-1)$ is a sufficiently small number. Indeed, we can repeat all the computations replacing $|u|^\alpha v$ by $|u|^{\alpha-1}vw$ or, to be more precise, replacing $|u|^\alpha v=|u|^\theta |u|^{\alpha-\theta }v$ by $|u|^{\alpha-1} vw=|u|^\theta |u|^{\alpha-1-\theta }vw$.

\ The last inequality is important in the perturbation theory. 
\end{remark} 

\ Before stating the next lemma, we define the following numbers:   		
\begin{equation}\label{paradmissivel1}
 k=\frac{4\alpha(\alpha+1-\theta)}{4-2b-\alpha}\hspace{1.5cm}\;p=\frac{6\alpha(\alpha+1-\theta)}{(4-2b)(\alpha-\theta)+\alpha}
\end{equation}
 and
\begin{equation}\label{paradmissivel2}
 l=\frac{4\alpha(\alpha+1-\theta)}{\alpha(3\alpha-2+2b)-\theta(3\alpha-4+2b)},
 \end{equation}	 
where $\theta\in (0,\alpha)$ small enough. It is easy to see that $(l,p)$ is $L^2$-admissible and $(k,p)$ is $\dot{H}^{s_c}$-admissible. In Appendix we verify the conditions of admissible pair.

\begin{lemma}\label{LG3} Let $N\geq 2$, $ \frac{4-2b}{N}<\alpha<2_*$ and $0<b<\min\{\frac{N}{3},1\}$. There exist $c>0$ such that
 \begin{equation*}
 \left\|\nabla (|x|^{-b}|u|^\alpha u)\right\|_{S'(L^2)}\leq c\| u\|^{\theta}_{L^\infty_tH^1_x}\|u\|^{\alpha-\theta}_{S(\dot{H}^{s_c})} \|\nabla u\|_{S(L^2)}+c\| u\|^{1+\theta}_{L^\infty_tH^1_x}\|u\|^{\alpha-\theta}_{S(\dot{H}^{s_c})},
\end{equation*}
where $\theta\in (0,\alpha)$ sufficiently small. 
\begin{proof} 
For $N\geq 4$ and $N=2$, the above inequality was already proved in \cite[Lemmas $4.3$, with $s=1$]{CARLOS} and \cite[Lemmas $4.7$, with $s=1$]{CARLOS}, respectively. Now, we only consider the case $N=3$. We claim that
$$
 \left\|\nabla (|x|^{-b}|u|^\alpha u)\right\|_{S'(L^2)}\leq c\| u\|^{\theta}_{L^\infty_tH^1_x}\|u\|^{\alpha-\theta}_{S(\dot{H}^{s_c})} \|\nabla u\|_{S(L^2)}.
$$

\ Indeed, in view of $(2,6)$ is $L^2$-admissible in $3D$ we deduce (dividing in $B$ e $B^C$) 
\begin{equation*}
 \left\|\nabla\left(|x|^{-b}|u|^\alpha u \right)\right\|_{S'(L^2)}\leq \left\|\nabla \left(|x|^{-b}|u|^\alpha u \right)\right\|_{L_t^{2'}L^{6'}_x(B)}+\left\|\nabla \left(|x|^{-b}|u|^\alpha u \right)\right\|_{L_t^{2'}L^{6'}_x(B^C)}.
\end{equation*}
Let $A\subset \mathbb{R}^N$. The product rule for derivatives and H\"older's inequality imply that
\begin{eqnarray}\label{E3D}
 \left\|  \nabla \left(|x|^{-b}|u|^\alpha u \right) \right \|_{L^{6'}_x(A)}&\leq& \left\|  |x|^{-b}\nabla \left(|u|^\alpha u \right) \right \|_{L^{6'}_x(A)}+\left\|  \nabla \left(|x|^{-b}\right)|u|^\alpha u \right \|_{L^{6'}_x(A)}\nonumber \\ 
 &\leq& M_1(t,A)+M_2(t,A),
\end{eqnarray}  
 where 
 $$
 M_1(t,A)=\left \||x|^{-b} \right \|_{L^\gamma(A)} \left \|\nabla(|u|^\alpha u)\right\|_{L^{\beta}_x}\;\;\;\;\;\;\;\;M_2(t,A)=\left\|\nabla(|x|^{-b})\right\|_{L^d(A)} \left\||u|^\alpha u\right\|_{L^{e}_x}
 $$    
and
\begin{equation}\label{SLG3}
\frac{1}{6'}=\frac{1}{\gamma}+\frac{1}{\beta}=\frac{1}{d}+\frac{1}{e}. 
 \end{equation}

\ First, we estimate $M_1(t,A)$. By H\"older's inequality we deduce
\begin{eqnarray}\label{LG4Hs3}
 M_1(t,A) &\leq& \||x|^{-b}\|_{L^\gamma(A)}  \|u\|^{\theta}_{L^{\theta r_1}_x}   \|u\|^{\alpha-\theta}_{L_x^{(\alpha-\theta)r_2}}  \|\nabla u\|_{ L^p_x}  \nonumber \\
&=&   \||x|^{-b}\|_{L^\gamma(A)}   \|u\|^{\theta}_{L^{\theta r_1}_x} \|u\|^{\alpha-\theta}_{L_x^p}  \|\nabla u\|_{ L^p_x},  
\end{eqnarray}
where 
\begin{equation}\label{LG4Hs31}
 \frac{1}{\beta}=\frac{1}{r_1}+\frac{1}{r_2}+\frac{1}{p}\;\;\;\;\textnormal{and}\;\;\;\; p=(\alpha-\theta)r_2.
\end{equation}  
Combining \eqref{SLG3} and \eqref{LG4Hs31} we obtain
\begin{equation*}
\frac{3}{\gamma}=\frac{5}{2}-\frac{3}{r_1}-\frac{3(\alpha+1-\theta)}{p},
\end{equation*}
 which implies, by \eqref{paradmissivel1}
\begin{equation}\label{LG4Hs32}
\frac{3}{\gamma}-b=\frac{\theta(2-b)}{\alpha}-\frac{3}{r_1}.
\end{equation}
In to order to show that $\||x|^{-b}\|_{L^\gamma(A)}$ is finite we need $\frac{3}{\gamma}-b>0$ if $A=B$ and $\frac{3}{\gamma}-b<0$ if $A=B^C$, by Remark \ref{RIxb}. Indeed if $\theta r_1=\frac{6}{3-2s}$, by \eqref{LG4Hs32} we have 
$$\frac{3}{\gamma}-b=\theta(s-s_c)>0$$ 
and if $\theta r_1=2$ then 
$$\frac{3}{\gamma}-b=-\theta s_c<0.$$
Therefore, the inequality \eqref{LG4Hs3} and the Sobolev embedding \eqref{SEI22} yield
\begin{equation}\label{M1tA}
 M_1(t,A) \leq c \|u\|^{\theta}_{H^1_x}\|u\|^{\alpha-\theta}_{L_x^p}  \|\nabla u\|_{ L^p_x}.  
\end{equation}

\  To estimate $M_2(t,A)$ we use the pairs $(\bar{a},\bar{r})=\left(4(\alpha-2\theta), \frac{6\alpha(\alpha-2\theta)}{\alpha(3-2b)-2\theta(4-2b)}\right)$ $\dot{H}^{s_c}$-admissible and $(q,r)=\left(\frac{4(\alpha-2\theta)}{\alpha-3\theta},\frac{6(\alpha-2\theta)}{2\alpha-3\theta}\right)$ $L^2$-admissible.
Applying the H\"older inequality and the Sobolev embedding \eqref{SEI} we get
\begin{eqnarray}\label{SLG32}
 M_2(t,A)& \leq& \||x|^{-b-1} \|_{L^d(A)} \|u\|^{\theta}_{L^{\theta r_1}_x}   \|u\|^{\alpha-\theta}_{L_x^{(\alpha-\theta)r_2}} \| u\|_{L_x^{ r_3}} \nonumber \\ 
&\leq &   \||x|^{-b-1}\|_{L^d(A)}  \|u\|^{\theta}_{L^{\theta r_1}_x} \| u\|^{\alpha-\theta}_{L_x^{\bar{r}}}  \|\nabla u\|_{L_x^r}
\end{eqnarray}
if
\begin{equation*}
\left\{\begin{array}{cl}
\frac{1}{e}=&\frac{1}{r_1}+\frac{1}{r_2}+\frac{1}{r_3}\\ 
 1=&\frac{3}{r}-\frac{3}{r_3}\\ 
\bar{r}=&(\alpha-\theta)r_2.
\end{array}\right.
\end{equation*}
Note that the second equation is valid since $r<3$. Similarly as before, in order to show that $\||x|^{-b-1}\|_{L^d(A)}$ is bounded, we need $\frac{3}{d}-b-1>0$ when $A$ is the ball $B$ and $\frac{3}{d}-b-1<0$ when $A=B^C$, by Remark \ref{RIxb}. In fact, it follows from \eqref{SLG3}, the previous system and the values of $q$, $r$, $\bar{q}$ and $\bar{r}$ defined above that
\begin{eqnarray}\label{SLG33}
\frac{3}{d}-b-1&=& \frac{5}{2}-b-\frac{3}{r_1}-\frac{3(\alpha-\theta)}{\bar{r}}-\frac{3}{r}  \nonumber\\
&=&\frac{\theta(2-b)}{\alpha}-\frac{3}{r_1}.
\end{eqnarray}
Choosing $r_1$ such that 
$$
\theta r_1>\frac{3\alpha}{2-b} \textrm{ when } A=B \quad \textrm{and} \quad \theta r_1<\frac{3\alpha}{2-b} \textrm{ when } A=B^C,
$$
we obtain $\frac{3}{d}-b-1>0$ and $\frac{3}{d}-b-1<0$, respectively, that is $|x|^{-b-1}\in L^d(A)$. In addition, by the Sobolev embedding \eqref{SEI22} (since $2<\frac{3\alpha}{2-b}<6$) and \eqref{SLG32}, it follows that   
\begin{equation*}
 M_2(t,A)\leq c \|u\|^{\theta}_{H^1_x}\|u\|^{\alpha-\theta}_{L^{\bar{r}}_x}   \|\nabla  u\|_{L_x^r}.
\end{equation*}

\ Therefore, combining \eqref{E3D}, \eqref{M1tA} and the last inequality we obtain
$$
 \left\|  \nabla \left(|x|^{-b}|u|^\alpha u \right) \right \|_{L^{6'}_x} \leq c \|u\|^{\theta}_{H^1_x}\|u\|^{\alpha-\theta}_{L_x^p}  \|\nabla u\|_{ L^p_x} +c \|u\|^{\theta}_{H^1_x}\|u\|^{\alpha-\theta}_{L^{\bar{r}}_x}   \|\nabla  u\|_{L_x^r}.
$$

\ Finally, using the H\"older inequality in the time variable (since $
\frac{1}{2'}=\frac{\alpha-\theta}{k}+\frac{1}{l}=\frac{\alpha-\theta}{\bar{a}}+\frac{1}{q})
$,
we conclude 
\begin{eqnarray*}
 \left\|  \nabla \left(|x|^{-b}|u|^\alpha u \right) \right \|_{L^{2'}_tL_x^{6'}}  &\leq& c \|u\|^{\theta}_{L^\infty_tH^1_x}\|u\|^{\alpha-\theta}_{L^k_tL_x^p}  \|\nabla u\|_{L^l_t L^p_x}+\|u\|^{\theta}_{L^\infty_tH^1_x}\|u\|^{\alpha-\theta}_{L^{\bar{a}}_tL^{\bar{r}}_x}   \|\nabla u\|_{L^q_tL_x^r}.
\end{eqnarray*}
The proof is completed recalling that $(q,r)$ and $(l,p)$ are $L^2$-admissible as well as $(k,p)$ and $(\bar{a},\bar{r})$ are $\dot{H}^{s_c}$-admissible.
\end{proof}
\end{lemma} 
\begin{remark}Note that, in the previous lemma we need the assumption $\alpha<3-2b$. Indeed, to verify that ($\bar{a},\bar{r}$) satisfies the condition of admissible pair (condition \eqref{CPA2} with $N=3$), we have to show  $\frac{3\alpha}{2-b}=\frac{6}{3-2s_c}<\bar{r}<6$. Note that $\bar{r}>\frac{3\alpha}{2-b}$ is equivalent to $2(\alpha-2\theta)(2-b)>\alpha(3-2b)-2\theta(4-2b)\Leftrightarrow \alpha>0$. Also, $\bar{r}<6\Leftrightarrow 2\theta(4-2b-\alpha)<\alpha(3-2b-\alpha)$, which is true if $\alpha<3-2b$ (since $\theta>0$ is a small number).
\end{remark}
\begin{remark}\label{RSglobal} We also have the following estimate (a consequence of the previous lemma)
$$
 \left\| |x|^{-b-1}|u|^\alpha v  \right\|_{S'(L^2)}\lesssim  \| u\|^{\theta}_{L^\infty_tH^1_x}\|u\|^{\alpha-\theta}_{S(\dot{H}^{s_c})} \left( \|\nabla v\|_{S(L^2)}+\| v\|_{L^\infty_tH^1_x} \right).
 $$
\end{remark} 



\ We now state our first result concerning the Cauchy problem \eqref{INLS}.
\begin{proposition}\label{GWPH1}{\bf (Small data global theory in $H^1$)}
Let $N\geq 2$, $\frac{4-2b}{N}<\alpha<2_*$ with $0<b<\min\{\frac{N}{3},1\}$ and $u_0 \in H^1(\mathbb{R}^N)$. Suppose $\|u_0\|_{H^1}\leq A$. There exists $\delta=\delta(A)>0$ such that if $\|U(t)u_0\|_{S(\dot{H}^{s_c})}<\delta$, then there exists a unique global solution $u$ of \eqref{INLS} such that
\begin{equation*}\label{NGWP3}
\|u\|_{S(\dot{H}^{s_c})}\leq 2\|U(t)u_0\|_{S(\dot{H}^{s_c})} 
\end{equation*}
and
\begin{equation*}\label{NGWP4}
\|u\|_{S\left(L^2\right)}+\|\nabla  u\|_{S\left(L^2\right)}\leq 2c\|u_0\|_{H^1}.
\end{equation*}
\begin{proof} See \cite[Theorem $1.9$, with $s=1$]{CARLOS}.
\end{proof}
\end{proposition}
\begin{remark} It is worth mentioning that the previous results were proved in \cite{CARLOS}  under the condition $0<b<\widetilde{2}$ (see definition \eqref{def2s}). Consequently, it is easy to see that they also hold for $0<b<\min\{\frac{N}{3},1\}$.
\begin{equation}\label{def2s}
\widetilde{2}:= \left\{\begin{array}{cl}
\frac{N}{3}&\;\;N=1,2,3,\\
2&\;\;N\geq 4.
\end{array}\right.
\end{equation}

\end{remark}


\ As mentioned in the introduction, Proposition \ref{SCATTERSH1} gives us the criterion to establish scattering. We prove it in the sequel.

\begin{proof}[\bf{Proof of Proposition \ref{SCATTERSH1}}] Note that
 \begin{equation}\label{SCATTER1}
 \|u\|_{S(L^2)}+\|\nabla u\|_{S(L^2)}<+\infty.
 \end{equation}
 
\ Indeed, using the fact that $\|u\|_{S(\dot{H}^{s_c})}<+\infty$, given $\delta>0$ we can decompose $[0,\infty)$ into $n$ intervals $I_j=[t_j,t_{j+1})$ such that $\|u\|_{S(\dot{H}^{s_c};I_j)}<\delta$ for all $j=1,...,n$. Let the integral equation  on the time interval $I_j$
\begin{equation*}\label{SCATTER2}
u(t)=U(t-t_j)u(t_j)+i\int_{t_j}^{t_{j+1}}U(t-s)(|x|^{-b}|u|^\alpha u)(s)ds.
\end{equation*}
Applying the Strichartz estimates \eqref{SE1} and \eqref{SE3}, we obtain  
\begin{equation}\label{SCATTER3} 
\|u\|_{S(L^2;I_j)}\leq c\|u(t_j)\|_{L^2_x}+c\left\||x|^{-b}|u|^\alpha u \right\|_{S'(L^2;I_j)}	
\end{equation}	 
and
\begin{equation}\label{SCATTER4}
\|\nabla u\|_{S(L^2;I_j)}\leq c\|\nabla u(t_j)\|_{L^2_x}+c\left\|\nabla(|x|^{-b}|u|^\alpha u) \right\|_{S'(L^2;I_j)}.
\end{equation}
From Lemmas \ref{LG2} and \ref{LG3} we have
\begin{eqnarray*}
\left\||x|^{-b}|u|^\alpha u \right\|_{S'(L^2;I_j)}  &\leq& c \| u \|^\theta_{L^\infty_{I_j}H^1_x}\| u \|^{\alpha-\theta}_{S(\dot{H}^{s_c};I_j)}\|u\|_{S(L^2;I_j)},
\end{eqnarray*}
 \begin{equation*}
\|\nabla(|x|^{-b}|u|^\alpha u)\|_{S'(L^2;I_j)}\leq c\| u\|^{\theta}_{L^\infty_{I_j}H^1_x}\|u\|^{\alpha-\theta}_{S(\dot{H}^{s_c};I_j)} \left( \|\nabla u\|_{S(L^2;I_j)} + \| u\|_{L^\infty_{I_j}H^1_x} \right).
\end{equation*}
Hence, the relations \eqref{SCATTER3}, \eqref{SCATTER4} and the two last estimates imply	
$$
\|u\|_{S(L^2;I_j)}\leq c B+cB^\theta\delta^{\alpha-\theta}\|u\|_{S(L^2;I_j)}
$$
and
\begin{equation}\label{nablau}
\|\nabla u\|_{S(L^2;I_j)}\leq c B+ cB^{\theta+1}\delta^{\alpha-\theta}+cB^\theta\delta^{\alpha-\theta}\|\nabla u\|_{S(L^2;I_j)},
\end{equation}
where we have used the assumption $\sup\limits_{t\in \mathbb{R}}\|u(t)\|_{H^1}\leq B$. \\
Taking $\delta>0$ such that $c B^\theta\delta^{\alpha-\theta}<\frac{1}{2}$ we deduce
$$
\| u\|_{S(L^2;I_j)}+	\|\nabla u\|_{S(L^2;I_j)} \leq cB,
$$
and by summing over the $n$ intervals, we conclude the proof of \eqref{SCATTER1}.	
	
 \ Returning to the proof of the proposition, let
$$
\phi^+=u_0+i\int\limits_{0}^{+\infty}U(-s)|x|^{-b}(|u|^\alpha u)(s)ds.
$$
It is easy to see that $\phi^+ \in H^1(\mathbb{R}^N)$. Indeed, by the same arguments as before, we have that
\begin{equation*}
 \|\phi^+\|_{L^2}\leq c\|u_0\|_{L^2}+c \| u \|^\theta_{L^\infty_tH^1_x}\| u \|^{\alpha-\theta}_{S(\dot{H}^{s_c})}\|u\|_{S(L^2)}
\end{equation*}	 
and
\begin{eqnarray*}
\|\nabla\phi^+\|_{L^2}&\leq & c\|\nabla u_0\|_{L^2}+c\| u\|^{\theta}_{L^\infty_tH^1_x}\|u\|^{\alpha-\theta}_{S(\dot{H}^{s_c})}\left( \| \nabla u\|_{S(L^2)} +\| u\|_{L^\infty_tH^1_x}\right). \end{eqnarray*}
Therefore, \eqref{SCATTER1} yields $\|\phi^+\|_{H^1}<+\infty$.
  
\ On the other hand, since $u$ is a solution of \eqref{INLS} we get
$$
 u(t)-U(t)\phi^+=-i\int\limits_{t}^{+\infty}U(t-s)|x|^{-b}(|u|^\alpha u)(s)ds.
$$
Moreover, we deduce (again as before)
$$
 \|u(t)-U(t)\phi^+\|_{L^2_x}\leq c \| u \|^\theta_{L^\infty_tH^1_x}\| u \|^{\alpha-\theta}_{S(\dot{H}^{s_c};[t,\infty))}\|u\|_{S(L^2)}
$$
and
\begin{eqnarray*}
\|\nabla(u(t)-U(t)\phi^+)\|_{L^2_x}   &\leq & c  \| u\|^{\theta}_{L^\infty_tH^1_x}\|u\|^{\alpha-\theta}_{S(\dot{H}^{s_c};[t,\infty))}  \left( \| \nabla u\|_{S(L^2)} +\| u\|_{L^\infty_tH^1_x} \right).
\end{eqnarray*}
Since $\|u\|_{S(\dot{H}^{s_c};[t,\infty))}\rightarrow 0$ as $t \rightarrow +\infty$ and using \eqref{SCATTER1}, we conclude that 
$$\|u(t)-U(t)\phi^+\|_{H^1_x}\rightarrow 0 \,\,\textnormal{as}\,\,t\rightarrow +\infty.$$ 

\ In the same way we define 
$$
\phi^-=u_0+i\int_0^{-\infty}U(-s)|x|^{-b}(|u|^\alpha u)(s)ds,
$$
so that we obtain $\phi^-\in H^1$ and 
$$
 u(t)-U(t)\phi^-=i\int\limits_{-\infty}^{t}U(t-s)|x|^{-b}(|u|^\alpha u)(s)ds,
$$
which also satisfies (using the same argument as before)
$$
\|u(t)-U(t)\phi^-\|_{H^1_x}\rightarrow 0 \,\,\textnormal{as}\,\,t\rightarrow -\infty.
$$
\end{proof}	

\ Now, the purpose is to study the perturbation theory for \eqref{INLS}. We begin proving the short-time perturbation result. 
\begin{proposition}\label{STP}{\bf (Short-time perturbation).} 
Let $I\subseteq \mathbb{R}$ be a time interval containing zero and let $\widetilde{u}$ defined on $I\times \mathbb{R}^N$ be a solution to 
\begin{equation*}\label{PE}
i\partial_t \widetilde{u} +\Delta \widetilde{u} + |x|^{-b} |\widetilde{u}|^\alpha \widetilde{u} =e,
\end{equation*}  
with initial data $\widetilde{u}_0\in H^1(\mathbb{R}^N)$, satisfying 
\begin{equation}\label{PC11}  
\sup_{t\in I}  \|\widetilde{u}(t)\|_{H^1_x}\leq M \;\; \textnormal{and}\;\; \|\widetilde{u}\|_{S(\dot{H}^{s_c}; I)}\leq \varepsilon,
\end{equation}
for some positive constant $M$ and some small $\varepsilon>0$.
	
\indent  Let $u_0\in H^1(\mathbb{R}^N)$ such that 
\begin{equation}\label{PC22}
\|u_0-\widetilde{u}_0\|_{H^1}\leq M'\;\; \textnormal{and}\;\; \|U(t)(u_0-\widetilde{u}_0)\|_{S(\dot{H}^{s_c}; I)}\leq \varepsilon,\;\;\textnormal{for }\; M'>0.
\end{equation}
Assume also the following conditions
\begin{equation}\label{PC33}
\|e\|_{S'(L^2; I)}+\|\nabla e\|_{S'(L^2; I)}+  \|e\|_{S'(\dot{H}^{-s_c}; I)}\leq \varepsilon.
\end{equation}
\indent There exists $\varepsilon_0(M,M')>0$ such that if $\varepsilon<\varepsilon_0$, then there is a unique solution $u$ to \eqref{INLS} on $I\times \mathbb{R}^N$ with initial data $u_0$, at the time $t=0$, satisfying 
\begin{equation}\label{C} 
\|u\|_{S(\dot{H}^{s_c}; I)}\lesssim \varepsilon 
\end{equation}
and
\begin{equation}\label{C1}
\|u\|_{S(L^2; I)}+\|\nabla u\|_{S(L^2; I)}\lesssim c(M,M').
\end{equation}
\begin{proof} First, we claim (we will show it later): if $\|\widetilde{u}\|_{S(\dot{H}^{s_c};I)}\leq \varepsilon_0$, for some $\varepsilon_0>0$ enough small, then 
\begin{equation}\label{ST1} 
\|\widetilde{u}\|_{S(L^2;I)} \lesssim M\;\;\;\textnormal{and}\;\;\;\;\|\nabla\widetilde{u}\|_{S(L^2;I)} \lesssim M.
\end{equation}
   
\ Assume, without loss of generality, that $0=\inf I$. First, we prove the existence of a solution $w$ for the following Cauchy problem  
	
\begin{equation}\label{IVPP} 
\left\{\begin{array}{cl}
i\partial_tw +\Delta w + H(x,\widetilde{u},w)+e= 0,&  \\
w(0,x)= u_0(x)-\widetilde{u}_0(x),& 
\end{array}\right.
\end{equation}
where $H(x,\widetilde{u},w)=|x|^{-b} \left(|\widetilde{u}+w|^\alpha (\widetilde{u}+w)-|\widetilde{u}|^\alpha \widetilde{u}\right)$.	
	
\ Indeed, let  
\begin{equation}\label{IEP} 
G (w)(t):=U(t)w_0+i \int_0^t U(t-s)(H(x,\widetilde{u},w)+e)(s)ds
\end{equation}
and define
$$
B_{\rho,K}=\{ w\in C(I;H^1(\mathbb{R}^N)):\;\|w\|_{S(\dot{H}^{s_c};I)}\leq \rho\;\textnormal{and}\;\|w\|_{S(L^2;I)}+\|\nabla w\|_{S(L^2;I)}\leq K    \}.
 $$
We need to show (for a suitable choice of the parameters $\rho>0$ and $K>0$) that $G$ in \eqref{IEP} defines a contraction on $B_{\rho,K}$.
Indeed, we deduce by the Strichartz inequalities (\ref{SE1}), (\ref{SE2}), (\ref{SE3}) and \eqref{SE5} that 
\begin{equation}\label{SP1}
\|G(w)\|_{S(\dot{H}^{s_c};I)}\lesssim \|U(t)w_0\|_{S(\dot{H}^{s_c};I)}+ \| H(\cdot,\widetilde{u},w) \|_{S'(\dot{H}^{-s_c};I)}+\|e \|_{S'(\dot{H}^{-s_c};I)}
\end{equation}
\begin{equation}\label{SP2}
\|G(w)\|_{S(L^2;I)}\lesssim \|w_0\|_{L^2}+ \| H(\cdot,\widetilde{u},w) \|_{S'(L^2;I)}+\|e\|_{S'(L^2;I)}
\end{equation}
and	
\begin{equation}\label{SP3}
\|\nabla G(w)\|_{S(L^2;I)}\lesssim  \|\nabla w_0\|_{L^2}+ \| \nabla H(\cdot,\widetilde{u},w)\|_{S'(L^2;I)}+\|\nabla e\|_{S'(L^2;I)}.
\end{equation}	
On the other hand, since 
\begin{equation}\label{EI} 
\left| |\widetilde{u}+w|^\alpha(\widetilde{u}+w)-|\widetilde{u}|^\alpha\widetilde{u} \right|\lesssim |\widetilde{u}|^\alpha|w|+|w|^{\alpha+1}
\end{equation}
we obtain (using \eqref{FEI})
$$
\|H(\cdot,\widetilde{u},w)\|_{S'(\dot{H}^{-s_c};I)}\leq \||x|^{-b}|\widetilde{u}|^\alpha w\|_{S'(\dot{H}^{-s_c};I)}+\||x|^{-b}|w|^\alpha w\|_{S'(\dot{H}^{-s_c};I)},
$$
which implies by Lemma \ref{LG1} that
\begin{align}\label{SP4}
 \|H(\cdot,\widetilde{u},w)\|_{S'(\dot{H}^{-s_c};I)}\lesssim  \left(\| \widetilde{u} \|^\theta_{L^\infty_tH^1_x}\| \widetilde{u} \|^{\alpha-\theta}_{S(\dot{H}^{s_c};I)}+ \| w \|^\theta_{L^\infty_tH^1_x} \| w\|^{\alpha-\theta}_{S(\dot{H}^{s_c};I)}    \right) \| w \|_{S(\dot{H}^{s_c};I)}.
\end{align}
Using Lemma \ref{LG2} we also have
\begin{align}\label{SP5}
\|H(\cdot,\widetilde{u},w)\|_{S'(L^2;I)}\lesssim \left(\| \widetilde{u} \|^\theta_{L^\infty_tH^1_x}\| \widetilde{u} \|^{\alpha-\theta}_{S(\dot{H}^{s_c};I)} + \| w \|^\theta_{L^\infty_tH^1_x}\| w\|^{\alpha-\theta}_{S(\dot{H}^{s_c};I)} \right)\| w \|_{S(L^2;I)}.
\end{align}
Now we are interested in estimating $\|\nabla H(\cdot,\widetilde{u},w)\|_{S'(L^2;I)}$. The relations \eqref{SECONDEI} and \eqref{EI} imply that
\begin{equation*} 
|\nabla H(x,\widetilde{u},w)| \lesssim |x|^{-b-1}(|\widetilde{u}|^{\alpha}+|w|^{\alpha})|w|+|x|^{-b}(|\widetilde{u}|^\alpha+|w|^\alpha) |\nabla w| +E,
 \end{equation*}
 where
\begin{eqnarray*} 
 E &\lesssim& \left\{\begin{array}{cl}
 |x|^{-b}\left(|\widetilde{u}|^{\alpha-1}+|w|^{\alpha-1}\right)|w||\nabla \widetilde{u}| & \textnormal{if}\;\;\;\alpha> 1 \vspace{0.2cm} \\ 
|x|^{-b}|\nabla \widetilde{u}||w|^{\alpha} & \textnormal{if}\;\;\;\alpha\leq 1.
\end{array}\right.
\end{eqnarray*}
Thus, Lemma \ref{LG2} and Remark \ref{RSglobal} lead to
$$
\|\nabla H(\cdot,\widetilde{u},w)\|_{S'(L^2;I)} \lesssim \left(\| \widetilde{u} \|^\theta_{L^\infty_tH^1_x}\| \widetilde{u} \|^{\alpha-\theta}_{S(\dot{H}^{s_c};I)} + \| w \|^\theta_{L^\infty_tH^1_x}\| w\|^{\alpha-\theta}_{S(\dot{H}^{s_c};I)} \right)\|\nabla w \|_{S(L^2;I)}
$$
$$
\hspace{2.0cm} +\left(\| \widetilde{u} \|^\theta_{L^\infty_tH^1_x}\| \widetilde{u} \|^{\alpha-\theta}_{S(\dot{H}^{s_c};I)} + \| w \|^\theta_{L^\infty_tH^1_x}\| w\|^{\alpha-\theta}_{S(\dot{H}^{s_c};I)} \right)\| w \|_{L^\infty_tH^1_x}
$$
$$
\hspace{2.4cm}+\left(\| \widetilde{u} \|^\theta_{L^\infty_tH^1_x}\| \widetilde{u} \|^{\alpha-\theta}_{S(\dot{H}^{s_c};I)} + \| w \|^\theta_{L^\infty_tH^1_x}\| w\|^{\alpha-\theta}_{S(\dot{H}^{s_c};I)} \right)\|\nabla w \|_{S(L^2;I)} +E_1
$$
$$
\hspace{3.2cm}\lesssim \left(\| \widetilde{u} \|^\theta_{L^\infty_tH^1_x}\| \widetilde{u} \|^{\alpha-\theta}_{S(\dot{H}^{s_c};I)} + \| w \|^\theta_{L^\infty_tH^1_x}\| w\|^{\alpha-\theta}_{S(\dot{H}^{s_c};I)} \right)\|\nabla w \|_{S(L^2;I)}
$$
\begin{equation}\label{SP6}
\hspace{2.8cm} +\left(\| \widetilde{u} \|^\theta_{L^\infty_tH^1_x}\| \widetilde{u} \|^{\alpha-\theta}_{S(\dot{H}^{s_c};I)} + \| w \|^\theta_{L^\infty_tH^1_x}\| w\|^{\alpha-\theta}_{S(\dot{H}^{s_c};I)} \right)\| w \|_{L^\infty_tH^1_x}+E_1,
\end{equation}
where (using Remark \ref{RGP})
\begin{align*}
E_1\lesssim & \left\{\begin{array}{cl}
\left(\| \widetilde{u} \|^\theta_{L^\infty_tH^1_x} \| \widetilde{u} \|^{\alpha-1-\theta}_{S(\dot{H}^{s_c};I)} + \| w \|^\theta_{L^\infty_tH^1_x} \| w \|^{\alpha-1-\theta}_{S(\dot{H}^{s_c};I)} \right) \| w \|_{S(\dot{H}^{s_c};I)}  \|\nabla \widetilde{u} \|_{S(L^2;I)},&\alpha> 1 \vspace{0.2cm} \\
\| w \|^\theta_{L^\infty_tH^1_x}\| w\|^{\alpha-\theta}_{S(\dot{H}^{s_c};I)} \|\nabla \widetilde{u} \|_{S(L^2;I)}\;,\;\;\;\;\;\;\;\;\;\;\alpha\leq 1.
\end{array}\right.
\end{align*}
 
\ Next, combining \eqref{SP4}, \eqref{SP5} and if $u\in B(\rho,K)$, we get
\begin{eqnarray}\label{SP7}
 \|H(\cdot,\widetilde{u},w)\|_{S'(\dot{H}^{-s_c};I)} \lesssim \left(M^\theta\varepsilon^{\alpha-\theta}+K^\theta \rho^{\alpha-\theta}\right)\rho
\end{eqnarray}
and
\begin{eqnarray}\label{SP8}
\|H(\cdot,\widetilde{u},w)\|_{S'(L^2;I)}\lesssim \left(M^\theta\varepsilon^{\alpha-\theta}+K^\theta \rho^{\alpha-\theta}\right)K.
\end{eqnarray}
In addition, \eqref{SP6} and \eqref{ST1} imply
\begin{align}\label{SP9}
 \|\nabla H(\cdot,\widetilde{u},w)\|_{S'(L^2;I)} \lesssim \left(M^\theta\varepsilon^{\alpha-\theta}+K^\theta \rho^{\alpha-\theta}\right)K +E_1
\end{align}
where
\begin{eqnarray*}
E_1&\lesssim & \left\{\begin{array}{cl}
\left( M^\theta \varepsilon^{\alpha-1-\theta} + K^\theta \rho^{\alpha-1-\theta} \right) \rho M&\textnormal{if}\;\;\alpha > 1, \vspace{0.2cm} \\
K^\theta \rho^{\alpha-\theta} M\;\;\;\;\;\; \textnormal{if}\;\;\;\;\;\alpha\leq 1.
\end{array}\right.
\end{eqnarray*}
Hence, it follows from \eqref{SP1}-\eqref{SP2} together with \eqref{SP7}- \eqref{SP8} that 
$$
\|G(w)\|_{S(\dot{H}^{s_c};I)}\leq  c\varepsilon+ cA\rho
$$
and
$$
\|G(w)\|_{S(L^2;I)}\leq cM'+c\varepsilon +cAK,
 $$
 where we also used the hypothesis \eqref{PC22}-\eqref{PC33} and $A=M^\theta\varepsilon^{\alpha-\theta}+K^\theta \rho^{\alpha-\theta}$. We also get, using \eqref{SP3}, \eqref{SP9}, that if $\alpha> 1$  
\begin{equation*}
\|\nabla G(w)\|_{S(L^2;I)}\leq cM'+c\varepsilon +cAK+cB \rho M,
\end{equation*}
where $B= M^\theta \varepsilon^{\alpha-1-\theta} + K^\theta \rho^{\alpha-1-\theta}$, and if $\alpha\leq 1$ 
\begin{equation*}
\|\nabla G(w)\|_{S(L^2;I)}\leq cM'+c\varepsilon + cAK + K^\theta \rho^{\alpha-\theta} M.
\end{equation*}
Choosing $\rho=2c\varepsilon$, $K=3cM'$ and $\varepsilon_0$ sufficiently small such that 
$$
cA<\frac{1}{3}\;\;\;\;\textnormal{and}\;\;\;c(\varepsilon+B \rho M+K^\theta \rho^{\alpha-\theta} M)<\frac{K}{3},
$$ 
we have
\begin{equation*}
\|G(w)\|_{S(\dot{H}^{s_c};I)}\leq \rho\;\;\;\textnormal{and}\;\;\;\|G(w)\|_{S(L^2;I)}+\|\nabla G(w)\|_{S(L^2;I)}\leq K.
\end{equation*}
Therefore, $G$ is well defined on $B(\rho,K)$. The contraction property can be obtained by  similar arguments. Thus, applying the Banach Fixed Point Theorem we obtain a unique solution $w$ on $I\times \mathbb{R}^N$ such that 
$$
\|w\|_{S(\dot{H}^{s_c};I)}\lesssim \varepsilon \;\;\;\textnormal{and}\;\;\;\|w\|_{S(L^2;I)}+\|w\|_{S(L^2;I)} \lesssim M'.
$$ 
Finally, it is easy to see that $u=\widetilde{u}+w$ is a solution to \eqref{INLS} satisfying \eqref{C} and \eqref{C1}. 
 
\ The proof is completed after showing \eqref{ST1}. Indeed, we first show that 
\begin{equation}\label{widetilde{u}}
\|\nabla\widetilde{u}\|_{S(L^2;I)}\lesssim M.
\end{equation}
We get using the same arguments as before
$$
\|\nabla\widetilde{u}\|_{S(L^2;I)}\lesssim  \|\nabla\widetilde{u}_0\|_{L^2}+ \left\|\nabla(|x|^{-b}|\widetilde{u}|^\alpha\widetilde{u})\right\|_{S'(L^2;I)}+\|\nabla e\|_{S'(L^2;I)}.
$$
Furthermore, Lemma \ref{LG3} implies that
\begin{eqnarray*}
\|\nabla\widetilde{u}\|_{S(L^2;I)}&\lesssim& M+ \| \widetilde{u} \|^\theta_{L^\infty_tH^1_x}\| \widetilde{u} \|^{\alpha-\theta}_{S(\dot{H}^{s_c};I)} \left( \|\nabla \widetilde{u} \|_{S(L^2;I)} + \| \widetilde{u} \|_{L^\infty_tH^1_x} \right)+\varepsilon\\
&\lesssim& M+\varepsilon+M^{\theta+1} \varepsilon_0^{\alpha-\theta}+  M^\theta \varepsilon_0^{\alpha-\theta}\|\nabla \widetilde{u} \|_{S(L^2;I)}.
\end{eqnarray*}
Therefore, choosing $\varepsilon_0$ sufficiently small the linear term $M^\theta \varepsilon_0^{\alpha-\theta}\|\nabla \widetilde{u} \|_{S(L^2;I)}$ may be absorbed by the left-hand term and we conclude the proof of \eqref{widetilde{u}}. Similar estimates also imply $\|\widetilde{u}\|_{S(L^2;I)}\lesssim M$.
\end{proof}	
\end{proposition}
\begin{remark}\label{RSP} 
In view of Proposition \ref{STP}, we also obtain the following estimates:
\begin{equation}\label{RSP1}
\|H(\cdot,\widetilde{u},w)\|_{S'(\dot{H}^{-s_c}; I)}\leq C(M,M') \varepsilon
\end{equation}
and
\begin{equation}\label{RSP2}
\|H(\cdot,\widetilde{u},w)\|_{S'(L^2; I)}+\|\nabla H(\cdot,\widetilde{u},w)\|_{S'(L^2; I)}\leq C(M,M')\varepsilon^{\alpha-\theta},
\end{equation}
with $\theta>0$ sufficiently small.

\ Indeed, it follows from \eqref{SP7}, \eqref{SP8} and \eqref{SP9} that
\begin{eqnarray*}
\|H(\cdot,\widetilde{u},w)\|_{S'(\dot{H}^{-s_c}; I)} \lesssim \left(M^\theta\varepsilon^{\alpha-\theta}+K^\theta \rho^{\alpha-\theta}\right)\rho,
\end{eqnarray*}
\begin{eqnarray*}
\|H(\cdot,\widetilde{u},w)\|_{S'(L^2; I)}\lesssim \left(M^\theta\varepsilon^{\alpha-\theta}+K^\theta \rho^{\alpha-\theta}\right)K
 \end{eqnarray*}
and
\begin{align*}
\|\nabla H(\cdot,\widetilde{u},w)\|_{S'(L^2; I)} \lesssim E_1+\left(M^\theta\varepsilon^{\alpha-\theta}+K^\theta \rho^{\alpha-\theta}\right)K,
\end{align*}
where
\begin{eqnarray*}
 E_1\lesssim & \left\{\begin{array}{cl}
\left( M^\theta \varepsilon^{\alpha-1-\theta} + K^\theta \rho^{\alpha-1-\theta} \right) \rho M&\textnormal{if}\;\;\alpha> 1, \vspace{0.2cm} \\
K^\theta \rho^{\alpha-\theta} M\;\;\;\;\;\; \textnormal{if}\;\;\;\;\;\alpha\leq 1.
\end{array}\right.
\end{eqnarray*}
Therefore, the choice $\rho=2c\varepsilon$ and $K=3cM'$ in Proposition \ref{STP} yield \eqref{RSP1} and \eqref{RSP2}. 
\end{remark}

\ Next, using the previous proposition we show the long-time perturbation result. This will be necessary in the
construction of the critical solution below.
\begin{proposition}\label{LTP}{\bf (Long-time perturbation)} 
Let $I\subseteq \mathbb{R}$ be a time interval containing zero and let $\widetilde{u}$ defined on $I\times \mathbb{R}^N$ be a solution to 
\begin{equation*}
i\partial_t \widetilde{u} +\Delta \widetilde{u} + |x|^{-b} |\widetilde{u}|^\alpha \widetilde{u} =e,
\end{equation*}  
with initial data $\widetilde{u}_0\in H^1(\mathbb{R}^N)$, satisfying (for some positive constants $M,L$)
\begin{equation}\label{HLP1} 
\sup_{t\in I}  \|\widetilde{u}\|_{H^1_x}\leq M \;\; \textnormal{and}\;\; \|\widetilde{u}\|_{S(\dot{H}^{s_c}; I)}\leq L.
\end{equation}
 
\indent Let $u_0\in H^1(\mathbb{R}^N)$ such that 
\begin{equation}\label{HLP2}
\|u_0-\widetilde{u}_0\|_{H^1}\leq M'\;\; \textnormal{and}\;\; \|U(t)(u_0-\widetilde{u}_0)\|_{S(\dot{H}^{s_c}; I)}\leq \varepsilon,
\end{equation}
for some positive constant $M'$ and some $0<\varepsilon<\varepsilon_1=\varepsilon_1(M,M',L)$.
In addition, assume also the following conditions
\begin{equation*}
\|e\|_{S'(L^2; I)}+\|\nabla e\|_{S'(L^2; I)}+ \|e\|_{S'(\dot{H}^{-s_c}; I)}\leq \varepsilon.
\end{equation*}
\indent Then, there exists a unique solution $u$ to \eqref{INLS} on $I\times \mathbb{R}^N$ with initial data $u_0$ at the time $t=0$ satisfying 
	
\begin{equation}\label{CLP} 
\|u-\widetilde{u}\|_{S(\dot{H}^{s_c}; I)}\leq C(M,M',L)\varepsilon\;\;\;\;\;\;\;\textnormal{and}
\end{equation}
\begin{equation}\label{CLP1}
\|u\|_{S(\dot{H}^{s_c}; I)} +\|u\|_{S(L^2; I)}+\|\nabla u\|_{S(L^2; I)}\leq C(M,M',L).
\end{equation}
	
\begin{proof} Since $\|\widetilde{u}\|_{S(\dot{H}^{s_c}; I)}\leq L$, given\footnote{$\varepsilon_0$ is given by the previous result and $\varepsilon$ to be determined later.} $\varepsilon<\varepsilon_0(M,2M')$ we can partition $I$ into $n = n(L,\varepsilon)$ intervals $I_j = [t_j ,t_{j+1})$ such that $\|\widetilde{u}\|_{S(\dot{H}^{s_c};I_j)}\leq \varepsilon$, for each $j$. Observe that $M'$ is  being replaced by $2M'$, as the $H^1$-norm of the difference of two different initial data may increase in each iteration.

\ Similarly as before, we can assume $0=\inf I$. Let $w$ be defined by $u = \widetilde{u} + w$, then $w$ solves IVP \eqref{IVPP} with initial time $t_j$. Thus, the integral equation in the interval $I_j = [t_j ,t_{j+1})$ reads as follows
  \begin{equation*}
   w(t)=U(t-t_j)w(t_j)+i\int_{t_j}^{t}U(t-s)(H(x,\widetilde{u},w)+e)(s)ds,
  \end{equation*}
 where $H(x,\widetilde{u},w)=|x|^{-b} \left(|\widetilde{u}+w|^\alpha (\widetilde{u}+w)-|\widetilde{u}|^\alpha \widetilde{u}\right)$.	
	
\ Choosing $\varepsilon_1$ sufficiently small (depending on $n$, $M$, and $M'$), we may apply Proposition \ref{STP} to obtain for each $0\leq j<n$ and all $\varepsilon<\varepsilon_1$, 
\begin{equation}\label{LP1}
\|u-\widetilde{u}\|_{S(\dot{H}^{s_c};I_j)}\leq C(M,M',j)\varepsilon
\end{equation}
and
\begin{equation}\label{LP2}
\|w\|_{S(\dot{H}^{s_c};I_j)}+\|w\|_{S'(L^2;I_j)}+\|\nabla w\|_{S'(L^2;I_j)}\leq C(M,M',j)
\end{equation}
provided we can prove (for each $0\leq j<n$)
\begin{equation}\label{LP3}
 \|U(t-t_j)(u(t_j)-\widetilde{u}(t_j))\|_{S(\dot{H}^{s_c};I_j)}\leq C(M,M',j)\varepsilon\leq \varepsilon_0
 \end{equation}
 and
 \begin{equation}\label{LP4}
 \|u(t_j)-\widetilde{u}(t_j)\|_{H^1_x}\leq 2M'.
 \end{equation}
 
\ Indeed, it follows from Strichartz estimates \eqref{SE2} and \eqref{SE5} that   
\begin{eqnarray*}
 \|U(t-t_j)w(t_j)\|_{S(\dot{H}^{s_c};I_j)}&\lesssim& \|U(t)w_0\|_{S(\dot{H}^{s_c}; I)}+\|H(\cdot,\widetilde{u},w)\|_{S'(\dot{H}^{-s_c};[0,t_j])}\\
  &&+\|e\|_{S'(\dot{H}^{-s_c};I)},
  \end{eqnarray*}
which implies by \eqref{RSP1} that	
$$
\|U(t-t_j)(u(t_j)-\widetilde{u}(t_j))\|_{S(\dot{H}^{s_c}; I_j)}\lesssim \varepsilon+\sum_{k=0}^{j-1}C(k,M,M')\varepsilon.
$$
  
\ In the same way, applying the Strichartz estimates \eqref{SE1}, \eqref{SE3} and \eqref{RSP2} we get
 \begin{eqnarray*}
 \|u(t_j)-\widetilde{u}(t_j)\|_{H^1_x}&\lesssim & \|u_0-\widetilde{u}_0\|_{H^1}+\|e\|_{S'(L^2; I)}+\|\nabla e\|_{S'(L^2;I)}\\
 &&+\| H(\cdot,\widetilde{u},w)\|_{S'(L^2;[0,t_j])}+\|\nabla H(\cdot,\widetilde{u},w)\|_{S'(L^2;[0,t_j])}\\
 &\lesssim& M'+\varepsilon+\sum_{k=0}^{j-1}C(k,M,M')\varepsilon^{\alpha-\theta}.
 \end{eqnarray*}
 Taking $\varepsilon_1=\varepsilon(n,M,M')$ sufficiently small, we see that \eqref{LP3} and \eqref{LP4} hold and so, it implies \eqref{LP1} and \eqref{LP2}.
 
\ We complete the proof summing this over all subintervals $I_j$, that is
$$
\|u-\widetilde{u}\|_{S(\dot{H}^{s_c}; I)}\leq C(M,M',L)\varepsilon
$$
and
$$
\|w\|_{S(\dot{H}^{s_c}; I)}+\|w\|_{S'(L^2; I)}+\|\nabla w\|_{S'(L^2; I)}\leq C(M,M',L).
 $$
 
 \end{proof} 
\end{proposition}

Finally, we show the existence of the Wave Operator. The proof follows the ideas introduced by C\^ote \cite{COTE} for the KdV equation (see also our paper \cite{paper2}).
\begin{proposition}{\bf (Existence of Wave Operator)}\label{PEWO} Assume $\phi \in H^1(\mathbb{R}^N)$ and 
\begin{equation}\label{HEWO}
\|\nabla \phi\|^{2s_c}_{L^2}\|\phi\|^{2(1-s_c)}_{L^2}<\lambda^2\left(\frac{N\alpha+2b}{\alpha s_c}\right)^{s_c}E[Q]^{s_c}M[Q]^{1-s_c},
\end{equation}
for some\footnote{Note that $(\frac{2\alpha s_c}{N\alpha+2b})^{\frac{s_c}{2}}<1$.} $0<\lambda \leq (\frac{2\alpha s_c}{N\alpha+2b})^{\frac{s_c}{2}}$. Then, there exists $u^+_0\in H^1(\mathbb{R}^N)$ such that $u$ solving \eqref{INLS} with initial data $u^+_0$ is global in $H^1(\mathbb{R}^N)$ with
\begin{itemize}
\item [(i)]$M[u]=M[\phi]$,
\item [(ii)] $E[u]=\frac{1}{2}\|\nabla \phi\|^2_{L^2}$,
\item [(iii)] $\lim\limits_{t\rightarrow +\infty} \|u(t)-U(t)\phi\|_{H^1}=0$,
\item [(iv)] $\|\nabla u(t)\|^{s_c}_{L^2}\| u(t)\|^{1-s_c}_{L^2}\leq \lambda\|\nabla Q\|^{s_c}_{L^2}\|Q\|^{1-s_c}_{L^2}$.
\end{itemize}
\begin{proof} 
First, we construct the wave operator for large time. Indeed, let $I_T=[T,+\infty)$ for $T\gg 1$ and define  
\begin{equation*}\label{IEWO1} 
G(w)(t)=- i \int_t^{+\infty} U(t-s)(|x|^{-b}|w+U(t)\phi|^\alpha (w+U(t)\phi)(s)ds,\;\;t\in I_T
\end{equation*}
and
$$
B(T,\rho)=\{w\in C\left(I_T;H^1(\mathbb{R}^N)\right): \|w\|_{T}\leq \rho   \},
$$
where
$$
\|w\|_T=\|w\|_{S(\dot{H}^{s_c};I_T)}+\|w\|_{S(L^2;I_T)}+\|\nabla w\|_{S(L^2;I_T)}. 
$$ 
We show a fixed point for $G$ on $B(T,\rho)$. 
  
\ The Strichartz estimates \eqref{SE3} \eqref{SE5} and Lemmas \ref{LG1}-\ref{LG2}-\ref{LG3} imply that  
\begin{align}\label{EWO1}    
\| G(w) \|_{S(\dot{H}^{s_c};I_T)} \lesssim &  \| w+U(t)\phi\|^\theta_{L^\infty_{T}H^1_x}\| w+U(t)\phi \|^{\alpha-\theta}_{S(\dot{H}^{s_c};I_T)}\| w+U(t)\phi \|_{S(\dot{H}^{s_c};I_T)}
\end{align}
\begin{align}\label{EWO2}
\| G(w) \|_{S(L^2;I_T)} \lesssim &  \| w+U(t)\phi\|^\theta_{L^\infty_{T}H^1_x}\| w+U(t)\phi \|^{\alpha-\theta}_{S(\dot{H}^{s_c};I_T)}\| w+U(t)\phi \|_{S(L^2;I_T)}
\end{align}
and
\begin{align}\label{EWO3}
\|\nabla G(w) \|_{S(L^2;I_T)} \lesssim &  \| w+U(t)\phi\|^\theta_{L^\infty_{T}H^1_x}\| w+U(t)\phi \|^{\alpha-\theta}_{S(\dot{H}^{s_c};I_T)}\| \nabla(w+U(t)\phi) \|_{S(L^2;I_T)}  \nonumber \\
&+\| w+U(t)\phi\|^{1+\theta}_{L^\infty_{T}H^1_x}\|w+U(t)\phi \|^{\alpha-\theta}_{S(\dot{H}^{s_c};I_T)}.
\end{align}	  
Hence,
\begin{eqnarray*}\label{EWO4}
\|G(w)\|_{T} &\lesssim &  \| w+U(t)\phi\|^\theta_{L^\infty_{T}H^1_x}\| w+U(t)\phi \|^{\alpha-\theta}_{S(\dot{H}^{s_c};I_T)}\| w+U(t)\phi \|_{T} \nonumber  \\
 && +\|w+U(t)\phi \|^{\alpha-\theta}_{S(\dot{H}^{s_c};I_T)}\|w+U(t)\phi \|^{\theta+1}_{T}.
\end{eqnarray*}
   
\ Since\footnote{Observe that \eqref{U(t)phi} is possible not true using the norm $L^{\infty}_{I_T}L^{\frac{2N}{N-2s_c}}_x$ and for this reason we remove the pair $\left(\infty,\frac{2N}{N-2s_c}\right)$ in the Definition \ref{Hsdefinition}.} 
\begin{equation}\label{U(t)phi}
\| U(t)\phi \|_{S(\dot{H}^{s_c};I_T)}\rightarrow 0
\end{equation}
as $T\rightarrow +\infty$, we can find $T_0>0$ large enough and $\rho>0$ small enough such that $G$ is well defined on $B(T_0,\rho)$. The same computations show that $G$ is a contraction on $B(T_0,\rho)$. Therefore, $G$ has a unique fixed point, that is $G(w)=w$. 
 
\ Next, using \eqref{EWO1} and since 
$$
\| w+U(t)\phi\|_{L^\infty_{T}H^1_x}\leq \|w\|_{H^1} +\|\phi\|_{H^1}<+\infty,
$$ 
we have
\begin{eqnarray*}    
\| w \|_{S(\dot{H}^{s_c};I_T)} &\lesssim &  \| w+U(t)\phi \|^{\alpha-\theta}_{S(\dot{H}^{s_c};I_T)}\| w+U(t)\phi \|_{S(\dot{H}^{s_c};I_T)}\\
 &\lesssim &  A\| w\|_{S(\dot{H}^{s_c};I_T)} +A\|U(t)\phi \|_{S(\dot{H}^{s_c};I_T)}
\end{eqnarray*}
where $A=\| w+U(t)\phi \|^{\alpha-\theta}_{S(\dot{H}^{s_c};I_T)}$. Moreover, if $\rho$ has been chosen small enough and since $\|U(t)\phi\|_{S(\dot{H}^{s_c};I_T)}$ is also sufficiently small for $T$ large, we deduce
$$
A\leq c\| w\|^{\alpha-\theta}_{S(\dot{H}^{s_c};I_T)}+c\| U(t)\phi \|^{\alpha-\theta}_{S(\dot{H}^{s_c};I_T)}<\frac{1}{2},
$$
and so
$$
\frac{1}{2} \| w \|_{S(\dot{H}^{s_c};I_T)} \lesssim A \|U(t)\phi \|_{S(\dot{H}^{s_c};I_T)},
$$
which implies,
\begin{equation}\label{EWO5}
  \| w \|_{S(\dot{H}^{s_c};I_T)}\rightarrow 0\;\;\;\;\textnormal{as}\;\;\;\; T\rightarrow +\infty.
\end{equation} 
The relations \eqref{EWO2}, \eqref{EWO3} and $\eqref{EWO5}$ also imply that\footnote{Note that $\| w+U(t)\phi \|_{S(\dot{H}^{s_c};I_T)}\leq \| w \|_{S(\dot{H}^{s_c};I_T)}+\| U(t)\phi \|_{S(\dot{H}^{s_c};I_T)} \rightarrow 0$ as $T\rightarrow +\infty$ by \eqref{EWO5} and $\| w+U(t)\phi\|^\theta_{L^\infty_{T}H^1_x}, \| w+U(t)\phi \|_{S(L^2;I_T)}, \| \nabla(w+U(t)\phi) \|_{S(L^2;I_T)} <\infty$ since $w\in B(T,\rho)$ and $\phi \in H^1(\mathbb{R}^N)$.}
$$
 \| w \|_{S(L^{2};I_T)}\;,\,\|\nabla
 w \|_{S(L^{2};I_T)}\rightarrow 0\;\;\;\;\textnormal{as}\;\;\;\; T\rightarrow +\infty,
 $$
and finally
\begin{equation}\label{EWO6}
   \|w\|_{T}\rightarrow 0 \;\; \textnormal{as}\;\; T\rightarrow +\infty.  
 \end{equation}
   
\ On the other hand, we claim that $u(t)=U(t)\phi+w(t)$ satisfies \eqref{INLS} in the time interval $[T_0,\infty)$. To do this, we need to show that
 \begin{equation}\label{CWO} 
 u(t)=U(t-T_0)u(T_0)+i\int_{T_0}^{t}U(t-s)(|x|^{-b}|u|^\alpha u)sds,
 \end{equation}
 for all $t\in [T_0,\infty)$. Indeed, since
  	$$
  	w(t)=- i \int_t^\infty U(t-s)|x|^{-b}|w+U(t)\phi|^\alpha (w+U(t)\phi)(s)ds,
  $$
we deduce
 \begin{eqnarray*}
  U(T_0-t)w(t)&=&- i \int_t^\infty U(T_0-s)|x|^{-b}|w+U(t)\phi|^\alpha (w+U(t)\phi)(s)ds\\
   	&=& i\int_{T_0}^t U(T_0-s)|x|^{-b}|w+U(t)\phi|^\alpha (w+U(t)\phi)(s)ds+w(T_0),
  \end{eqnarray*}
  and so applying $U(t-T_0)$ on both sides, we obtain
  $$
  w(t)=U(t-T_0)w(T_0)+ i \int_{T_0}^t U(t-s)|x|^{-b}|w+U(t)\phi|^\alpha (w+U(t)\phi)(s)ds.
  $$
  Finally, adding $U(t)\phi$ in both sides of the last equation, we deduce \eqref{CWO}.
  
 \ Our goal now is to show relations (i)-(iv). Since $u(t)=U(t)\phi+w$ then  
 \begin{align}\label{EWO7}
 \|u(t)-U(t)\phi\|_{L^\infty_TH^1_x}=\|w\|_{L^\infty_TH^1_x}\leq c\|w\|_{S(L^2;I_T)}+c\|\nabla w\|_{S(L^2;I_T)}\leq c\|w\|_{T},
  \end{align}
which implies (iii) (using $\eqref{EWO2}$). Moreover, it is easy to see, by \eqref{EWO7} 
 \begin{equation}\label{EWO81}
 \lim_{t\rightarrow\infty}\|u(t)\|_{L^2_x}=\| \phi\|_{L^2}.
 \end{equation}
 and
 \begin{equation}\label{EWO8}
 \lim_{t\rightarrow\infty}\|\nabla u(t)\|_{L^2_x}=\|\nabla \phi\|_{L^2}.
 \end{equation}
The mass conservation \eqref{mass} yields $\|u(t)\|_{L^2}=\|u(T_0)\|_{L^2}$ for all $t$, so from \eqref{EWO81} we deduce $\|u(T_0)\|_{L^2}=\|\phi\|_{L^2}$, i.e., item (i) holds. On the other hand, applying Lemma \ref{LEWO} (ii) we deduce 
 \begin{eqnarray*}
\left\| |x|^{-b}|u(t)|^{\alpha+2} \right\|_{L^1_x}&\leq& c\left\| |x|^{-b}|u(t)-U(t)\phi|^{\alpha+2} \right\|_{L^1_x}+c\left\| |x|^{-b}|U(t)\phi|^{\alpha+2} \right\|_{L^1_x}\\
&\leq & c\left\|u(t)-U(t)\phi| \right\|^{\alpha+2}_{H^1_x}+c\left\| |x|^{-b}|U(t)\phi|^{\alpha+2} \right\|_{L^1_x},
\end{eqnarray*}
which goes to zero as $t\rightarrow +\infty$, by item (iii) and Lemma \ref{LEWO} (iii), i.e.
\begin{equation}\label{EWO9}
 \lim_{t\rightarrow\infty}\left\| |x|^{-b}|u(t)|^{\alpha+2} \right\|_{L^1_x}=0.
 \end{equation} 
We have (ii) combining \eqref{EWO8} and \eqref{EWO9}.
 
 \ In view of \eqref{HEWO}, (i) and (ii) it follows that
 $$
 E[u]^{s_c}M[u]^{1-s_c}=\frac{1}{2^{s_c}}\|\nabla \phi\|^{2s_c}_{L^2}\|\phi\|^{2(1-s_c)}_{L^2}<\lambda^2\left(\frac{N\alpha+2b}{2\alpha s_c}\right)^{s_c}E[Q]^{s_c}M[Q]^{1-s_c}
  $$
 and by our choice of $\lambda$ we conclude
 \begin{equation*}\label{EWO10}
 E[u]^{s_c}M[u]^{1-s_c}<E[Q]^{s_c}M[Q]^{1-s_c}.
 \end{equation*}
 Furthermore, from \eqref{EWO81}, \eqref{EWO8} and \eqref{HEWO}
\begin{eqnarray*}
\lim_{t\rightarrow \infty}\|\nabla u(t)\|^{2s_c}_{L^2_x}\|u(t)\|^{2(1-s_c)}_{L^2_x}&=&\|\nabla \phi\|^{2s_c}_{L^2}\|\phi\|^{2(1-s_c)}_{L^2}\\
&<& \lambda^2\left(\frac{N\alpha+2b}{\alpha s_c}\right)^{s_c}E[Q]^{s_c}M[Q]^{1-s_c}\\
&=&\lambda^2\|\nabla Q\|^{2s_c}_{L^2}\|Q\|^{2(1-s_c)}_{L^2}
\end{eqnarray*}
where we have used \eqref{EGS}. Thus, one can take $T_1>0$ sufficiently large such that
\begin{equation*}\label{EWO11}
\|\nabla u(T_1)\|^{s_c}_{L^2_x}\|u(T_1)\|^{1-s_c}_{L^2_x}<\lambda \|\nabla Q\|^{s_c}_{L^2}\|Q\|^{1-s_c}_{L^2}.
\end{equation*}
Therefore, since $\lambda<1$, we deduce that relations \eqref{EMC} and \eqref{GFC} hold with $u_0=u(T_1)$ and so, by Theorem \ref{TG}, we have in fact that $u(t)$ constructed above is a global solution of \eqref{INLS}.  
\end{proof}		
\end{proposition}

\begin{remark}\label{backward}
A similar Wave Operator construction also holds when the time limit is taken as $t\rightarrow -\infty$.
\end{remark}


%



\section{Profile and energy decomposition}

\ We start by recalling some elementary inequalities (see G\'erard \cite{Ge98} inequality (1.10) and Guevara \cite{GUEVARA} page 217). Let $(z_j)\subset\mathbb{C}^M$ with $M\geq 2$. For all $q>1$ there exists $C_{q,M}>0$ such that 
\begin{equation}\label{FI}
 \left|\;\left | \sum_{j=1}^M z_j \right|^q-\sum_{j=1}^M|z_j|^q \right|  \leq C_{q,M}\sum_{j\neq k}^{M} |z_j| |z_k|^{q-1},
\end{equation}  
and for $\beta>0$ there exists a constant $C_{\beta,M}>0$ such that 
\begin{equation}\label{EIerror} 
\left| \left|\sum_{j=1}^{M}z_j\right|^\beta\sum_{j=1}^{M}z_j-\sum_{j=1}^{M}  |z_j|^\beta z_j\right|\leq C_{\beta,M}\sum_{j=1}^{M}\sum_{1\leq j\neq k\leq M}|z_j|^\beta|z_k|.
\end{equation} 

\ Our goal in this section is to establish a profile decomposition result and an Energy Pythagorean expansion for such a decomposition. To this end, we use similar arguments as in our work \cite{paper2}, with $(N,\alpha)=(3,2)$, and for the sake of completeness, we provide the details here. 

\begin{proposition}\label{LPD} {\bf (Profile decomposition)}
	Let $\phi_n(x)$ be a radial uniformly bounded sequence in $H^1(\mathbb{R}^N)$. Then for each $M\in \mathbb{N}$ there exists a  subsequence of $\phi_n$ (also denoted by $\phi_n$), such that, for each $1\leq j\leq M$, there exist a profile $\psi^j$ in $H^1(\mathbb{R}^N)$, a sequence $t_n^j$ of time shifts and a sequence $W_n^M$ of remainders in $H^1(\mathbb{R}^N)$, such that	
\begin{equation}\label{Aproximation}
\phi_n(x)=\sum_{j=1}^{M}U(-t_n^j)\psi^j(x)+W_n^M(x)
\end{equation}
with the following properties: 	  	
\begin{itemize}
\item \textsl{Pairwise divergence} for the time sequences. For $1\leq k\neq j\leq M$,
 \begin{equation}\label{PD} 
 \lim\limits_{n \rightarrow +\infty}|t_n^j-t_n^k|=+\infty.
 \end{equation}  
 \item \textsl{Asymptotic smallness} for the remainder sequence (recalling $s_c=\frac{N}{2}-\frac{2-b}{\alpha}$)
\begin{equation}\label{AS}
   \lim\limits_{M \rightarrow +\infty}\left(\lim\limits_{n \rightarrow +\infty}\|U(t)W_n^M\|_{S(\dot{H}^{s_c})}\right)=0.
  \end{equation}
		
\item \textsl{Asymptotic Pythagoream expansion}. For fixed $M\in \mathbb{N}$ and any $s\in [0,1]$, we have
\begin{equation}\label{PDNHs}
\|\phi_n\|^2_{\dot{H}^s}=\sum_{j=1}^{M}\|\psi^j\|^2_{\dot{H}^s}+\|W_n^M\|^2_{\dot{H}^s}+o_n(1)
\end{equation}
where $o_n(1) \rightarrow 0$ as $n\rightarrow +\infty$.
\end{itemize}
\begin{proof} Consider $\|\phi_n\|_{H^1}\leq C_1$, for some $C_1>0$.
Let $(a,r)$ $\dot{H}^{s_c}$-admissible and define $r_1=2r$, $a_1=\frac{4r}{r(N-2s_c)-N}$. It is easy to see that $(a_1,r_1)$ is also $\dot{H}^{s_c}$-admissible, thus combining the interpolation inequality with $\eta=\frac{N}{r(N-2s_c)-N}\in (0,1)$ and the Strichartz estimate \eqref{SE2}, we deduce
 \begin{eqnarray}\label{LPD1} 
  \|U(t)W_n^M\|_{L_t^aL^r_x}&\leq&\|U(t)W_n^M\|^{1-\eta}_{L_t^{a_1}L^{r_1}_x}\|U(t)W_n^M\|^\eta_{L_t^\infty L^{\frac{2N}{N-2s_c}}_x}\nonumber\\
  & \leq & \|W_n^M\|^{1-\eta}_{\dot{H}^{s_c}}\|U(t)W_n^M\|^\eta_{L_t^\infty L^{\frac{2N}{N-2s_c}}_x}.
\end{eqnarray}
So it will be suffice to conclude (since  $\|W_n^M\|_{\dot{H}^{s_c}}\leq C_1$)
 \begin{equation}\label{LPD2}
\lim\limits_{M \rightarrow +\infty}\left(\limsup\limits_{n \rightarrow +\infty}\|U(t)W_n^M\|_{L_t^\infty L^{\frac{2N}{N-2s_c}}_x}\right)=0.
\end{equation}

\ Indeed, we start by constructing $\psi^1_n$, $t_n^1$ and $W_n^1$. Let
$$
A_1=\limsup \limits_{n\rightarrow +\infty} \|U(t)\phi_n\|_{L_t^\infty L^{\frac{2N}{N-2s_c}}_x}.
$$
If $A_1=0$, we take $\psi^j=0$ for all $j=1,\dots,M$ and the proof is complete. Suppose $A_1>0$. Passing to a subsequence, we may consider $A_1=\lim \limits_{n\rightarrow +\infty} \|U(t)\phi_n\|_{L_t^\infty L^{\frac{2N}{N-2s_c}}_x}$. We claim that there exist a time sequence $t_n^1$ and $\psi^1$ such that $U(t_n^1)\phi_n \rightharpoonup \psi^1$ and
\begin{equation}\label{LPD22}
\beta C_1^ { \frac{N-2s_c}{2s_c(1-s_c)} }\|\psi^1\|_{\dot{H}^{s_c}}\geq A_1^{\frac{N-2s_c^2}{2s_c(1-s_c)}}, 
\end{equation} 
where $\beta>0$ is independent of $C_1$, $A_1$ and $\phi_n$.
Indeed, let $\zeta\in C^\infty_0(\mathbb{R}^N)$ a real-valued and radially symmetric function such that $0\leq \zeta \leq 1$, $\zeta(\xi)=1$ for $|\xi|\leq 1$ and $\zeta(\xi)=0$ for $|\xi|\geq 2$. Given $r>0$, define $\chi_r$ by $\widehat{\chi_r}(\xi)=\zeta(\frac{\xi}{r})$. It follows from Sobolev embedding \eqref{SEI} and since the operator $U(t)$ is an isometry in $H^{s_c}$ that\footnote{Recalling $0<s_c<1$.}    
\begin{align*}\label{LPD23}
\|U(t)\phi_n -\chi_r*U(t)\phi_n\|^2_{L^\infty_tL_x^{\frac{2N}{N-2s_c}}} &\leq c\|U(t)\phi_n -\chi_r*U(t)\phi_n\|^2_{L^\infty_tH_x^{s_c} } \nonumber \\
 &\leq c \int |\xi|^{2s_c}|(1-\widehat{\chi_r})^2|\widehat{\phi}_n(\xi)|^2d\xi  \nonumber \\
 &\leq c r^{-2(1-s_c)}\|\phi\|^2_{\dot{H}^1} \leq c r^{-2(1-s_c)}C_1^2.
\end{align*}
Taking    
\begin{equation}\label{LPD24}
 r=\left(\frac{4\sqrt{c}C_1}{A_1}\right)^{\frac{1}{1-s_c}}
\end{equation}
 and for $n$ large enough we obtain
 \begin{equation}\label{LPD3}
\|\chi_r*U(t)\phi_n\|_{L^\infty_tL_x^{\frac{2N}{N-2s_c}}}\geq \frac{A_1}{2}.
\end{equation} 
Observe that, from the standard interpolation in Lebesgue spaces
\begin{eqnarray}\label{LPD4}
\|\chi_r*U(t)\phi_n\|^N_{L^\infty_tL_x^{\frac{2N}{N-2s_c}}}&\leq&\|\chi_r*U(t)\phi_n\|^{N-2s_c}_{L^\infty_tL_x^2} \|\chi_r*U(t)\phi_n\|^{2s_c}_{L^\infty_tL_x^\infty}\nonumber  \\
&\leq & C_1^{N-2s_c}\|\chi_r*U(t)\phi_n\|^{2s_c}_{L^\infty_tL_x^\infty},
\end{eqnarray}
thus (using \eqref{LPD3} and \eqref{LPD4})
$
 \|\chi_r*U(t)\phi_n\|_{L^\infty_tL_x^\infty}\geq \left(\frac{A_1}{2C_1^{\frac{N-2s_c}{N}}}  \right)^{\frac{N}{2s_c}}
$. Since all $\phi_n$ are radial functions and so are $\chi_r*U(t)\phi_n$, the radial Sobolev Gagliardo-Nirenberg inequality \eqref{RSGN1} leads to
\begin{eqnarray*}
\|\chi_r*U(t)\phi_n\|_{L^\infty_tL_x^\infty(|x|\geq R)}&\leq &\frac{1}{R^{\frac{N-1}{2}}} \|\chi_r*U(t)\phi_n\|^{\frac{1}{2}}_{L^2_x}\|\nabla(\chi_r*U(t)\phi_n)\|^{\frac{1}{2}}_{L^2_x} \\
&\leq& \frac{C_1}{R^{\frac{N-1}{2}}}.
\end{eqnarray*} 
Combining these last inequalities we obtain for $R>0$ large
 $$
 \|\chi_r*U(t)\phi_n\|_{L^\infty_tL_x^\infty(|x|\leq R)}\geq \frac{1}{2} \left(\frac{A_1}{2C_1^{\frac{N-2s_c}{N}}}  \right)^{\frac{N}{2s_c}}.
 $$
Let $t_n^1$ and $x_n^1$, with $|x_n^1|\leq R$, be sequences such that for each $n\in \mathbb{N}$  
 $$
 \left|\chi_r*U(t_n^1)\phi_n(x_n^1)\right|\geq \frac{1}{4}\left(\frac{A_1}{2C_1^{\frac{N-2s_c}{N}}}  \right)^{\frac{N}{2s_c}}
 $$
 or 
\begin{equation}\label{A_1}
 \frac{1}{4}\left(\frac{A_1}{2C_1^{\frac{N-2s_c}{N}}}  \right)^{\frac{N}{2s_c}}\leq \left|\int \chi_r(x_n^1-y)U(t_n^1)\phi_n(y)dy\right|. 
 \end{equation}
Since $\|U(t_n^1)\phi_n\|_{H^1}=\|\phi_n\|_{H^1}\leq C_1$  then $U(t^1_n)\phi_n$ converges weakly in $H^1$ ($U(t_n^1)\phi_n$, that is there exists $\psi^1$ a radial function such that $U(t_n^1)\phi_n \rightharpoonup \psi^1$ in $H^1$ and $\|\psi^1\|_{H^1}\leq \limsup \limits_{n\rightarrow +\infty}\|\phi_n\|_{H^1}\leq C_1$. Moreover, $x_n^1\rightarrow x^1$ since $x_n^1$ is bounded. Thus, the inequality \eqref{A_1}, the Plancherel formula and the Cauchy-Schwarz inequality imply
 $$ 
 \frac{1}{8}\left(\frac{A_1}{2C_1^{\frac{N-2s_c}{N}}}  \right)^{\frac{N}{2s_c}}\leq \left|\int \chi_r(x^1-y) \psi^1(y)dy \right|\leq \|\chi_r\|_{\dot{H}^{-s_c}}\|\psi^1\|_{\dot{H}^{s_c}},
 $$
 which implies (using $
 \|\chi_r\|_{\dot{H}^{-s_c}}\leq c\left(\int_{0}^{2r}\rho^{-2s_c}\rho^{N-1}d\rho\right)^\frac{1}{2}\leq cr^{\frac{N-2s_c}{2}}$)
$$ 
 \frac{1}{8}\left(\frac{A_1}{2C_1^{\frac{N-2s_c}{N}}}  \right)^{\frac{N}{2s_c}}\leq cr^{\frac{N-2s_c}{2}}\|\psi^1\|_{\dot{H}^{s_c}}.
$$
Therefore in view of our choice of $r$ (see \eqref{LPD24}) we deduce \eqref{LPD22}, concluding the claim. 
 
 \ Define $W^1_n =\phi_n-U(-t_n^1)\psi^1$. Given any $0\leq s\leq 1$, it follows that  
 \begin{itemize}
\item $U(t_n^1)W^1_n \rightharpoonup 0$ in $H^1$ (since $U(t_n^1)\phi_n \rightharpoonup \psi^1$),
\item  $\langle \phi_n,U(-t^1_n)\psi^1 \rangle_{\dot{H}^s}=\langle U(t^1_n)\phi_n,\psi^1 \rangle_{\dot{H}^s}\rightarrow \|\psi^1\|^2_{\dot{H}^s}$,
\item $\|W_n^1\|^2_{\dot{H}^s}=\|\phi_n\|^2_{\dot{H}^s}-\|\psi^1\|^2_{\dot{H}^s}+o_n(1)$.
\end{itemize}
The last item, with $s=0$ and $s=1$, implies $\|W_n^1\|_{H^1}\leq C_1$.

\ Next, let $A_2=\limsup \limits_{n\rightarrow +\infty}\|U(t)W_n^1\|_{L_t^\infty L^{\frac{2N}{N-2s}}_x}$. If $A_2=0$, there is nothing to prove. Again the only case we need to consider is $A_2>0$. Repeating the above procedure, with $\phi_n$ replaced by $W_n^1$ we obtain a sequence $t_n^2$ and a function $\psi^2$ such that $U(t_n^2)W_n^1\rightharpoonup \psi^2$ in $H^1$ and 
\begin{equation*}
\beta C_1^ { \frac{N-2s_c}{2s_c(1-s_c)} }\|\psi^2\|_{\dot{H}^{s_c}}\geq A_2^{\frac{N-2s_c^2}{2s_c(1-s_c)}}.
\end{equation*}

\ We now show that $|t_n^2-t_n^1|\rightarrow +\infty$. We suppose that $t_n^2-t_n^1\rightarrow t^*$ finite, then 
$$
 U(t_n^2-t_n^1)\left(U(t_n^1)\phi_n-\psi^1 \right)=U(t_n^2)\left(\phi_n-U(-t_n^1)\psi^1 \right)=U(t_n^2)W_n^1\rightharpoonup \psi^2.
$$
On the other hand, since $U(t_n^1)\phi_n\rightharpoonup \psi^1$, the left side of the above expression converges weakly to $0$, and thus $\psi^2=0$, a contradiction. Let $W_n^2=W_n^1-U(-t_n^2)\psi^2$. We get for any $0\leq s\leq 1$ (using the fact that $|t_n^1-t_n^2|\rightarrow +\infty$)
\begin{eqnarray*}
\langle \phi_n,U(-t_n^2)\psi^2 \rangle_{\dot{H}^{s}}&=&\langle U(t_n^2)\phi_n,\psi^2 \rangle_{\dot{H}^{s}}\\
&=&\langle U(t_n^2)W_n^1,\psi^2 \rangle_{\dot{H}^{s}}+\langle U(t_n^2-t_n^1)\psi^1,\psi^2 \rangle_{\dot{H}^{s}} \\
&\rightarrow& \|\psi^2\|^2_{\dot{H}^{s}}.
\end{eqnarray*}
The definition of $W_n^2$ also yields that 
$$
\|W_n^2\|^2_{\dot{H}^s}=\|W_n^1\|^2_{\dot{H}^{s_c}}-\|\psi^2\|^2_{\dot{H}^s}+o_n(1)\;\;\;\textnormal{and}\;\;\;\|W_n^2\|_{H^1}\leq C_1
$$

\ By induction we can construct $\psi^M$, $t_n^M$ and $W_n^M$ such that $U(t_n^M)W_n^{M-1}\rightharpoonup \psi^M$ in $H^1$ and 
\begin{equation}\label{LPD5}
\beta C_1^ { \frac{N-2s_c}{2s_c(1-s_c)} }\|\psi^M\|_{\dot{H}^{s_c}}\geq A_M^{\frac{N-2s_c^2}{2s_c(1-s_c)}},
\end{equation} 
where $A_M=\lim \limits_{n\rightarrow +\infty}\|U(t)W_n^{M-1}\|_{L_t^\infty L^{\frac{2N}{N-2s_c}}_x}$. 

\ Next, we prove \eqref{PD}. Assume $1\leq j<M$, we show that $|t^M_n-t_n^j|\rightarrow +\infty$ by induction assuming $|t^M_n-t_n^k|\rightarrow +\infty$ for $k=j+1, \dots, M-1$. Indeed, let $t^M_n-t_n^j\rightarrow t_0$ finite then 
$$
U(t_n^M-t_n^j)\left(U(t_n^j)W_n^{j-1}-\psi^j\right)-U(t_n^M-t_n^{j+1})\psi^{j+1}-...-U(t_n^M-t_n^{M-1})\psi^{M-1}
$$
$$
=U(t_n^M)W_n^{M-1}\rightharpoonup \psi^M.
$$
 Since the left side converges weakly to $0$, we have $\psi^M=0$, a contradiction. 
 
\ We now consider $
W_n^M=\phi_n-U(-t_n^1)\psi^1-U(-t_n^2)\psi^2-...-U(-t_n^M)\psi^M$.
As before, by \eqref{PD} we get 
$$
\langle \phi_n,U(-t_n^M)\psi^M \rangle_{\dot{H}^s}=\langle U(t_n^M)W_n^{M-1},\psi^M \rangle_{\dot{H}^s}+o_n(1),
$$
and $\langle \phi_n,U(-t_n^M)\psi^M \rangle_{\dot{H}^s}\rightarrow \|\psi^M\|^2_{\dot{H}^s}$. Hence expanding $\|W_n^M\|^2_{\dot{H}^s}$ we conclude that \eqref{PDNHs} also holds. 

 \ The relations \eqref{LPD5} and \eqref{PDNHs} yield 
$
\sum_{M\geq 1} \left(\frac{A_M^{\frac{N-2s_c^2}{s_c(1-s_c)}}}{\beta^2C_1^{ \frac{N-2s_c}{s_c(1-s_c)}  }}\right)\leq \lim\limits_{n\rightarrow+\infty}\|\phi_n\|^2_{\dot{H}^{s_c}}<+\infty,
$
 which implies that $A_M\rightarrow 0$ as $M\rightarrow +\infty$ i.e., \eqref{LPD2} holds. Therefore, we get \eqref{AS} by \eqref{LPD1}. This completes the proof.
\end{proof} 
\end{proposition}

\ The next proposition contains an energy Pythagoream expansion. To this end, we use the following remark.
\begin{remark}\label{RLPD} 
It follows from the proof of Proposition \ref{LPD} that 
\begin{equation}\label{RLPD1}
\lim\limits_{M,n\rightarrow \infty} \|W_n^{M}\|_{L^{p}} =0,
\end{equation}
where\footnote{Recalling $\widetilde{p}$ is defined in \eqref{defr*}.} $2<p<\widetilde{p}$. Indeed, it is easy to see that
\begin{equation}\label{RLPD2}
\lim\limits_{M \rightarrow +\infty}\left(\lim\limits_{n \rightarrow +\infty}\|U(t)W_n^M\|_{L_t^\infty L^p_x}\right)=0.
\end{equation}
We have $\dot{H}^{s}\hookrightarrow L^{p}$, where $s=\frac{N}{2}-\frac{N}{p}$ (see inequality \eqref{SEI}). Since $2<p<\widetilde{p}$ then $0<s<1$, thus repeating the argument used for showing \eqref{LPD2} with $\frac{2N}{N-2s_c}$ replaced by $p$ and $s_c$ by $s$, we get \eqref{RLPD2}. In addition, \eqref{RLPD1} follows directly from \eqref{RLPD2} and the inequality
$$
\|W_n^{M}\|_{L_x^{p}}\leq \|U(t)W_n^{M}\|_{L^\infty_tL_x^{p}},
$$
since $W_n^{M}=U(0)W_n^{M}$.
\end{remark}
\begin{proposition}\label{EPE}{\bf (Energy Pythagoream Expansion)} Under the hypothesis of Proposition \ref{LPD} we obtain
\begin{equation}\label{PDE} 
E[\phi_n]=\sum_{j=1}^{M}E[U(-t_n^j)\psi^j]+E[W_n^M]+o_n(1).
\end{equation}
\begin{proof} We get (using \eqref{PDNHs} with $s=1$)
$$
E[\phi_n]-\sum_{j=1}^{M}E[U(-t_n^j)\psi^j]-E[W_n^M]=-\frac{A_n}{\alpha+2}+o_n(1),
$$
where 
$$
A_n=\left\|  |x|^{-b}|\phi_n|^{\alpha+2}  \right\|_{L^1}-\sum_{j=1}^{M}\left\|  |x|^{-b} |U(-t_n^j)\psi^j|^{\alpha+2} \right\|_{L_x^1}-\left\|  |x|^{-b}|W_n^M|^{\alpha+2}  \right\|_{L^1}.
$$

\ For a fixed $M\in \mathbb{N}$, if $A_n \rightarrow 0$ as $n\rightarrow +\infty$ then \eqref{PDE} holds. Indeed, pick $M_1\geq M$ and rewrite the last expression as
\begin{eqnarray*}
A_n&=& \int \left(  |x|^{-b}|\phi_n|^{\alpha+2}  -  \sum_{j=1}^{M}   |x|^{-b}  |U(-t_n^j)\psi^j|^{\alpha+2} -  |x|^{-b} |W_n^M|^{\alpha+2}   \right)dx     \\
&=&  I^1_n+I^2_n+I^3_n,     
\end{eqnarray*}   
where 
\begin{eqnarray*}
I^1_n&=&  \int   |x|^{-b} \left[ |\phi_n|^{\alpha+2}-|\phi_n-W_n^{M_1}|^{\alpha+2}\right]dx \\
I^2_n&=& \int |x|^{-b} \left[ |W_n^{M_1}-W_n^M|^{\alpha+2} -|W_n^M|^{\alpha+2}\right]dx   
\end{eqnarray*}
\begin{eqnarray*}
I^3_n&=& \int  |x|^{-b}\left[ |\phi_n-W_n^{M_1}|^{\alpha+2}- \sum_{j=1}^{M} |U(-t_n^j)\psi^j|^{\alpha+2}- |W_n^{M_1}-W_n^M|^{\alpha+2}   \right] dx .
\end{eqnarray*}

\ We start by estimating $I^1_n$. Lemma \ref{LEWO} (i)-(ii) and \eqref{FI} imply that 
\begin{eqnarray*}
|I^1_n|&\lesssim& \int |x|^{-b}\left( |\phi_n|^{\alpha+1}|W_n^{M_1}|+|\phi_n||W_n^{M_1}|^{\alpha+1}  + |W_n^{M_1}|^{\alpha+2}  \right)dx\\
&\lesssim&  \left(  \|\phi_n\|^{\alpha+1}_{L^r}\|W_n^{M_1}\|_{L^r}+ \|\phi_n\|_{L^r}\|W_n^{M_1}\|^{\alpha+1}_{L^{r}} +  \|W_n^{M_1}\|^{\alpha+2}_{L^{r}} \right)+\\
& & \left(  \|\phi_n\|^{\alpha+1}_{L^{\alpha+2}}\|W_n^{M_1}\|_{L^{\alpha+2}}+ \|\phi_n\|_{L^{\alpha+2}}\|W_n^{M_1}\|^{\alpha+1}_{L^{\alpha+2}} + \|W_n^{M_1}\|^{\alpha+2}_{{L^{\alpha+2}}} \right)\\
&\lesssim&    \|\phi_n\|^{\alpha+1}_{H^1}\|W_n^{M_1}\|_{L^r}+ \|\phi_n\|_{H^1}\|W_n^{M_1}\|^{\alpha+1}_{L^r} +  \|W_n^{M_1}\|^{\alpha+2}_{L^r} +\\
 & &  \|\phi_n\|^{\alpha+1}_{H^1}\|W_n^{M_1}\|_{L^{\alpha +2}}+ \|\phi_n\|_{H^1}\|W_n^{M_1}\|^{\alpha+1}_{L^{\alpha +2}} +  \|W_n^{M_1}\|^{\alpha+2}_{L^{\alpha +2}} ,
\end{eqnarray*}
where $\frac{N(\alpha+2)}{N-b}<r<\widetilde{p}$ (recall that $\widetilde{p}$ is defined in \eqref{defr*}). Using \eqref{RLPD1}(we can apply Remark \ref{RLPD} since $r$ and $\alpha+2\in (2,\widetilde{p})$) and since $\{\phi_n\}$ is uniformly bounded in $H^1$, we obtain
$$
I^1_n\rightarrow +\infty\;\;\textnormal{as}\;\;n, M_1\rightarrow +\infty.
$$  

\ In the same way (replacing $\phi_n$ by $W_n^{M}$) we also get 
$$
I^2_n\rightarrow +\infty\;\;\textnormal{as}\;\;n, M_1\rightarrow +\infty.
$$

\ Finally we consider the term $I^3_n$. Since,
$$
\phi_n-W_n^{M_1}=\sum\limits_{j=1}^{M_1}U(-t_n^j)\psi^j\;\;\textnormal{and}\;\;W_n^{M}-W_n^{M_1}=\sum\limits_{j=M+1}^{M_1}U(-t_n^j)\psi^j,
$$
we can rewrite $I^3_n$ as
\begin{eqnarray*}
I^3_n=\int |x|^{-b} \left( \left| \sum\limits_{j=1}^{M_1} U(-t_n^j)\psi^j \right|^{\alpha+2}- \sum\limits_{j=1}^{M_1}| U(-t_n^j)\psi^j |^{\alpha+2}\right)dx
\end{eqnarray*}
$$
-\int |x|^{-b} \left( \left| \sum\limits_{j=M+1}^{M_1} U(-t_n^j)\psi^j \right|^{\alpha+2}- \sum\limits_{j=M+1}^{M_1}| U(-t_n^j)\psi^j |^{\alpha+2}\right)dx.
$$

To complete the proof we make use of the following claim.

 \textit{Claim.} For a fixed $M_1\in \mathbb{N}$ and for some $j_0\in \mathbb{N}$ ($j_0< M_1$), we get
$$
D_n= \left\||x|^{-b} \left|\sum_{j=j_0}^{M_1}  U(-t_n^j)\psi \right|^{\alpha+2}  \right\|_{L^1_x}  -    \sum_{j=j_0}^{M_1} \left\|   |x|^{-b} |U(-t_n^j)\psi^j|^{\alpha+2}   \right\|_{L^1_x}\rightarrow 0,
$$
as $n\rightarrow +\infty$. 

Indeed, it is clear that the last limit implies that $I^3_n\rightarrow 0\;\textnormal{as}\;n\rightarrow +\infty$ completing the proof of relation \eqref{PDE}.

We now show the claim. Observe that \eqref{FI} implies
$$
D_n\leq \sum_{j\neq k}^{M_1}\int |x|^{-b}|U(-t_n^j)\psi^j||U(-t_n^k)\psi^k|^{\alpha+1}dx.
$$
Setting $E^{j,k}_n=\int |x|^{-b}|U(-t_n^j)\psi^j||U(-t_n^k)\psi^k|^{\alpha+1}dx$ and using Lemma \ref{LEWO} (i), we deduce
\begin{equation*}\label{lematecnico1}
E^{j,k}_n\leq  c  \| U(-t_n^k)\psi^k\|^{\alpha+1}_{L^{\alpha+2}_x}\| U(-t_n^j)\psi^j\|_{L^{\alpha+2}_x}+c\| U(-t_n^k)\psi^k\|^{\alpha+1}_{L^r_x}\| U(-t_n^j)\psi^j\|_{L^r_x},
\end{equation*}
where $2<\frac{N(\alpha+2)}{N-b}< r<\widetilde{p}$. Since \eqref{PD} we can consider that $t_n^k$, $t_n^j$ or both go to infinity as $n$ goes to infinity. If $t_n^j\rightarrow +\infty$ as $n\rightarrow +\infty$ then 
\begin{eqnarray*}
E^{j,k}_n & \leq&  c \|\psi^k\|^{\alpha+1}_{H^1}\| U(-t_n^j)\psi^j\|_{L^{\alpha+2}_x} +c\| \psi^k\|^{\alpha+1}_{H^1} \| U(-t_n^j)\psi^j\|_{L^r_x}\\
& \leq&  c \| U(-t_n^j)\psi^j\|_{L^{\alpha+2}_x} +c\| U(-t_n^j)\psi^j\|_{L^r_x},
\end{eqnarray*}
where in the last inequality we have used that $(\psi^k)_{k\in\mathbb{N}}$ is a uniformly bounded sequence in $H^1$. Thus, if $n\rightarrow +\infty$ we have $t_n^j\rightarrow +\infty$ and by \eqref{LEWO3} with $t=t_n^j$ and $f=\psi^j$ we conclude that $E^{j,k}_n\rightarrow 0$ as $n\rightarrow +\infty$. Similarly, for the case $t^k_n\rightarrow +\infty$ as $n\rightarrow +\infty$, we have $E^{j,k}_n\rightarrow 0$. Finally, in view of $D_n$ is a finite sum of terms in the form of $E^{j,k}$, we conclude that $D_n\rightarrow 0$ as $n\rightarrow +\infty$. 
\end{proof}
\end{proposition}

\section{Critical solution}\label{CCS}

\ In this section, we study a critical solution of \eqref{INLS} (denoted by $u_c$). First, assuming that $\delta_c<E[u]^{s_c}M[u]^{1-s_c}$ (see \eqref{deltac}), we construct $u_c$ of \eqref{INLS} with infinite Strichartz norm $\|\cdot\|_{S(\dot{H}^{s_c})}$ satisfying $$
E[u_c]^{s_c}M[u_c]^{1-s_c}=\delta_c.
$$ 
After that, we show that the flow associated to this critical solution is precompact in $H^1(\mathbb{R}^N)$. The key ingredients here are the results of the previous section and the long time perturbation theory (Proposition \ref{LTP}). 
\begin{proposition}\label{ECS}{\bf (Existence of $u_c$)} Let $N\geq 2$, $ \frac{4-2b}{N}<\alpha<2_*$ and $0<b<\min\{\frac{N}{3},1\}$. 
If
\begin{equation*}
\delta_c<E[Q]^{s_c}M[Q]^{1-s_c},
\end{equation*}
then there exists a radial function $u_{c,0}\in H^1(\mathbb{R}^N)$ such that the corresponding solution $u_c$ of the IVP \eqref{INLS} is global in $H^1(\mathbb{R}^N)$. Moreover the following properties hold
\begin{itemize}
\item [(i)] $M[u_c]=1$,
\item [(ii)] $E[u_c]^{s_c}=\delta_c$,
\item [(iii)] $\|  \nabla u_{c,0} \|_{L^2}^{s_c} \|u_{c,0}\|_{L^2}^{1-s_c}<\|\nabla Q \|_{L^2}^{s_c} \|Q\|_{L^2}^{1-s_c}$,
\item [(iv)] $\|u_{c}\|_{S(\dot{H}^{s_c})}=+\infty$.
\end{itemize}
\begin{proof}
There exists a sequence of solutions $u_n$ to \eqref{INLS} with $H^1$
 initial data $u_{n,0}$, with $\|u_n\|_{L^2} = 1$ for all $n\in \mathbb{N}$, such that (see section \ref{SPMR}) 
 \begin{equation}\label{PCS1}   
 \|\nabla u_{n,0}\|^{s_c}_{L^2} <  \|\nabla Q\|^{s_c}_{L^2}\|Q\|^{1-s_c}_{L^2}
 \end{equation}
 and
 \begin{equation*}\label{PCS2}
 E[u_n] \searrow \delta_c^{\frac{1}{s_c}}\;\; \textnormal{as}\;\; n \rightarrow +\infty.
 \end{equation*} 	
Also 
\begin{equation}\label{un}
\|u_n\|_{S(\dot{H}^{s_c})}=+\infty
\end{equation}
for every $n\in \mathbb{N}$. Since $\delta_c<E[Q]^{s_c}M[Q]^{1-s_c}$, there exists $a \in (0,1)$ such that
\begin{equation}\label{PCS21}
 E[u_n]\leq  a E[Q]M[Q]^{\sigma},
\end{equation}
where $\sigma=\frac{1-s_c}{s_c}$. Moreover, Lemma \ref{LGS} (ii) and \eqref{PCS1} yield
\begin{equation*}
 \|\nabla u_{n,0}\|^{2}_{L^2} \leq w^{\frac{1}{s_c}} \|\nabla Q\|^{2}_{L^2}\|Q\|^{2\sigma}_{L^2},
\end{equation*}
where $w=\frac{E[u_n]^{s_c}M[u_n]^{1-s_c}}{E[Q]^{s_c}M[Q]^{1-s_c}}$, thus we deduce from \eqref{PCS21} and $\|u_n\|_{L^2} = 1$ that $w^{\frac{1}{s_c}}\leq a$ which implies
\begin{equation}\label{PCS22}
\|\nabla u_{n,0}\|^{2}_{L^2} \leq a \|\nabla Q\|^{2}_{L^2}\|Q\|^{2\sigma}_{L^2}.
\end{equation}

\ On the other hand, we have using the linear profile decomposition (Proposition \ref{LPD}) applied to $u_{n,0}$, which is uniformly bounded in $H^1(\mathbb{R}^N)$ by \eqref{PCS22} that
 \begin{equation}\label{PCS3}
 u_{n,0}(x)=\sum_{j=1}^{M}U(-t_n^j)\psi^j(x)+W_n^M(x),
 \end{equation} 
where $M$ will be taken large later. By the Pythagorean expansion \eqref{PDNHs}, with $s=0$, that for all $M\in \mathbb{N}$ we deduce
\begin{equation}\label{PCS4}
\sum_{j=1}^{M}\|\psi^j\|^2_{L^2}+\lim_{n\rightarrow +\infty}\|W_n^M\|^2_{L^2}\leq \lim_{n\rightarrow +\infty}\|u_{n,0}\|^2_{L^2}= 1,
\end{equation}
which implies 
\begin{equation}\label{PCS41}
\sum_{j=1}^{M}\|\psi^j\|^2_{L^2}\leq 1.
\end{equation}
Another application of \eqref{PDNHs}, with $s=1$, and \eqref{PCS22} lead to
\begin{equation}\label{Sumpsij}
\sum_{j=1}^{M}\|\nabla \psi^j\|^2_{L^2}+\lim_{n\rightarrow +\infty}\|\nabla W_n^M\|^2_{L^2}\leq \lim_{n\rightarrow +\infty}\|\nabla u_{n,0}\|^2_{L^2}\leq a\|\nabla Q\|^2_{L^2}\|Q\|^{2\sigma}_{L^2},
\end{equation}
and so
\begin{equation}\label{PCS5}
\|\nabla \psi^j\|^{s_c}_{L^2}\leq a^{\frac{s_c}{2}}\|\nabla Q\|^{s_c}_{L^2}\| Q\|^{1-s_c}_{L^2},\;\;j=1,\dots,M.
\end{equation}

\ Let $\{t^j_n\}_{n\in \mathbb{N}}$ be the sequence given by Proposition \ref{LPD}. Combining \eqref{PCS41} and \eqref{PCS5} we obtain\footnote{Recalling that $U(t)$ is an isometry in $L^2(\mathbb{R}^N)$ and $\dot{H}^1(\mathbb{R}^N)$.}
 \begin{equation*}\label{U(-t_n^j)}
 \|U(-t_n^j)\psi^j\|^{1-s_c}_{L^2_x}\|\nabla U(-t_n^j)\psi^j\|^{s_c}_{L^2_x}\leq a^{\frac{s_c}{2}} \|\nabla Q\|^{s_c}_{L^2}\| Q\|^{1-s_c}_{L^2}.
 \end{equation*}
Also, we have by Lemma \ref{LGS} (i)
\begin{equation}\label{PCS6}
E[U(-t_n^j)\psi^j]\geq c(N,b,\alpha)\|\nabla \psi^j\|_{L^2}\geq 0.
\end{equation}
 
\ Similarly as before, for all $M\in \mathbb{N}$ we also get
 $$
 \lim_{n\rightarrow +\infty}\|W_n^M\|^2_{L^2}\leq 1,
 $$	
 $$
 \lim_{n\rightarrow +\infty}\|\nabla W_n^M\|^{s_c}_{L^2}\leq a^{\frac{s_c}{2}}\|\nabla Q\|^{s_c}_{L^2}\|Q\|^{1-s_c}_{L^2},
 $$	
and for $n$ large 
\begin{equation}\label{PCS7}
E[W_n^M]\geq 0.
\end{equation}

 \ The energy Pythagorean expansion (Proposition \ref{EPE}) allows us to deduce that
$$
\sum_{j=1}^{M}\lim_{n\rightarrow+\infty}E[U(-t_n^j)\psi^j]+\lim_{n\rightarrow+\infty}E[W_n^M]=\lim_{n\rightarrow+\infty}E[u_{n,0}]=\delta_c^{\frac{1}{s_c}},
$$
which implies (using \eqref{PCS6} and \eqref{PCS7}) that  
\begin{equation}\label{PCS8}
\lim_{n\rightarrow\infty}E[U(-t_n^j)\psi^j]\leq \delta_c^{\frac{1}{s_c}},\;\textnormal{for all}\;\;j=1,...,M.
\end{equation}

\ We now analyze two cases: if more than one $\psi^j\neq 0$ and only one profile is nonzero.   
 
\ If more than one $\psi^j\neq 0$, we prove a contradiction. Indeed, by \eqref{PCS4} we must have $M[\psi^j]<1$ for each $j$. Passing to a subsequence, if necessary, again we have two cases to consider.

 \ {\bf Case $1$}. $t^j_n\rightarrow t^*$ finite.\footnote{Note that, at most only one such $j$ exists by \eqref{PD}).} By the continuity of the linear flow in $H^1(\mathbb{R}^N)$    
 \begin{equation}\label{widetildepsi}
  U(-t_n^j)\psi^j\rightarrow U(-t^*)\psi^j\;\;\;\;\textnormal{strongly in}\;H^1.
 \end{equation}
We denote the solution of \eqref{INLS} with initial data $\psi$ by INLS$(t)\psi$. Set $\widetilde{\psi}^j=\textnormal{INLS}(t^*)(U(-t^*)\psi^j)$ so that $\mbox{INLS}(-t^*)\widetilde{\psi}^j=U(-t^*)\psi^j$. In view of the set
\begin{equation*}
\mathcal{K}:=\left\{u_0\in H^1(\mathbb{R}^N):\;  \textrm{relations} \; \eqref{EMC} \; \textrm{and} \; \eqref{GFC} \; \textrm{hold}\;\right\}
\end{equation*}
is closed in $H^1(\mathbb{R}^N)$ then $\widetilde{\psi}^j\in \mathcal{K}$, that is, INLS$(t)\widetilde{\psi}^j$ is a global solution by Theorem \ref{TG}. In addition, the relations \eqref{U(-t_n^j)}, \eqref{PCS8} and $M[\psi^j]<1$ implies that
$$
\|\widetilde{\psi}^j\|^{1-s_c}_{L^2_x}\|\nabla \widetilde{\psi}^j\|^{s_c}_{L^2_x}\leq \|\nabla Q\|^{s_c}_{L^2}\| Q\|^{1-s_c}_{L^2}\;\;\textnormal{and}\;\;E[\widetilde{\psi}^j]^{s_c}M[\widetilde{\psi}^j]^{1-s_c}< \delta_c.
$$
So, using the definition of $\delta_c$ (see \eqref{deltac}) we have
\begin{equation}\label{CSCP} 
\|\textnormal{INLS}(t)\widetilde{\psi}^j\|_{S(\dot{H}^{s_c})}<+\infty.
\end{equation}
Finally, it is easy to see by \eqref{widetildepsi} 
\begin{equation}\label{CSWO3}
\lim_{n\rightarrow+\infty}\|\textnormal{INLS}(-t_n^j)\widetilde{\psi}^j-U(-t_n^j)\psi^j\|_{H^1_x}=0.
\end{equation}

\ {\bf Case $2$.} If $|t^j_n|\rightarrow+\infty$ then by Lemma \ref{LEWO} (iii)
$$
\left\||x|^{-b}|U(-t_n^j)\psi^j|^{\alpha+2}\right\|_{L^1_x}\rightarrow 0,
$$
and thus, using the fact that $U(t)$ is an isometry in $\dot{H}^1(\mathbb{R}^N)$ and \eqref{PCS8}
\begin{equation}\label{CS0}
\left( \frac{1}{2}\|\nabla \psi^j\|^2_{L^2}\right)^{s_c}= \lim_{n \rightarrow \infty}  E [ U(-t_n^j) \psi^j]^{s_c} \leq  \delta_c < E[Q]^{s_c}M[Q]^{1-s_c},
\end{equation}
Therefore, by the existence of wave operator, Proposition \ref{PEWO} with $\lambda=(\frac{2\alpha s_c}{N\alpha+2b})^{\frac{s_c}{2}}<1$, there exists $\widetilde{\psi}^j\in H^1(\mathbb{R}^N)$ such that
\begin{equation}\label{CSWO1}
M[\widetilde{\psi}^j]=M[\psi^j]\;\;\;\textrm{ and }\;\;\;\;E[\widetilde{\psi}^j]=\frac{1}{2}\|\nabla \psi^j\|^2_{L^2},
\end{equation}
\begin{equation}\label{CSWO2}
 \|\nabla \textnormal{INLS}(t)\widetilde{\psi}^j\|^{s_c}_{L^2_x}\|\widetilde{\psi}^j\|^{1-s_c}_{L^2}<\|\nabla Q\|^{s_c}_{L^2}\| Q\|^{1-s_c}_{L^2} 
  \end{equation}
  and \eqref{CSWO3} also holds in this case.
  
Since $M[{\psi}^j]<1$ and using \eqref{CS0}-\eqref{CSWO1}, we get $E[\widetilde{\psi}^j]^{s_c}M[\widetilde{\psi}^j]^{1-s_c}<\delta_c$. Thus, the definition of $\delta_c$  together with \eqref{CSWO2} also imply \eqref{CSCP}.
 
In either case, we have a new profile $\widetilde{\psi}^j$ for the given $\psi^j$ such that \eqref{CSWO3} \eqref{CSCP} hold. 
	
 We now define
$$
 u_n(t)=\textnormal{INLS}(t)u_{n,0},
 $$	
 $$
 v^j(t)=\textnormal{INLS}(t)\widetilde{\psi}^j,
$$
 $$
 \widetilde{u}_n(t)=\sum_{j=1}^{M}v^j(t-t_n^j),
 $$
\begin{equation}\label{CSR}
\widetilde{W}_n^M=\sum_{j=1}^{M}\left[ U(-t_n^j)\psi^j-\textnormal{INLS}(-t_n^j)\widetilde{\psi}^j \right]+W_n^M.
\end{equation}
Then $\widetilde{u}_n(t)$ solves the following equation
\begin{equation}\label{widetildeun}
 i\partial_t\widetilde{u}_n+\Delta \widetilde{u}_n+|x|^{-b}|\widetilde{u}_n|^{\alpha}\widetilde{u}_n=e_n^M,
\end{equation}
where
\begin{equation}\label{CSR1}
 e_n^M=|x|^{-b}\left( |\widetilde{u}_n|^{\alpha}\widetilde{u}_n-\sum_{j=1}^{M}|v^j(t-t_n^j)|^{\alpha}v^j(t-t_n^j) \right).
\end{equation}
By definition of $\widetilde{W}_n^M$ in \eqref{CSR} and \eqref{PCS3}we can write
  \begin{equation*}\label{aproximation1}
    u_{n,0}=\sum_{j=1}^{M}\textnormal{INLS}(-t_n^j)\widetilde{\psi}^j+\widetilde{W}_n^M,
  \end{equation*}
so $u_{n,0}-\widetilde{u}_n(0)=\widetilde{W}_n^M$. Combining \eqref{CSR} together with the Strichartz inequality \eqref{SE2}, we estimate
\begin{equation*}
\|U(t)\widetilde{W}_n^M\|_{S(\dot{H}^{s_c})}\leq c\sum_{j=1}^{M}\|\textnormal{INLS}(-t_n^j)\widetilde{\psi}^j-U(-t_n^j)\psi^j\|_{H^1}+\|U(t)W_n^M\|_{S(\dot{H}^{s_c})},
  \end{equation*}
  which implies (using \eqref{AS} and \eqref{CSWO3})
 \begin{equation}\label{CSR2}
 \lim_{M\rightarrow +\infty} \left[\lim_{n\rightarrow +\infty}  \|U(t)(u_{n,0}-\widetilde{u}_{n,0})\|_{S(\dot{H}^{s_c})}\right]=0, 
\end{equation}
  
\ Next we approximate $u_n$ by $\widetilde{u}_n$. Then, it follow from the long time perturbation theory (Proposition \ref{LTP}) and \eqref{CSCP} that 
$$
\|u_n\|_{S(\dot{H}^{s_c})}<+\infty,
$$    
for $n$ large enough, which is a contradiction with \eqref{un}. Indeed, assume the following two claims\footnote{These claims will be proved in the next subsection.} to conclude the proof.\\
{\bf Claim $1$.} For each $M$ and $\varepsilon>0$, there exists $n_0=n_0(M,\varepsilon)$ such that
 \begin{equation}\label{claim2}
n>n_0\;\; \Rightarrow\;\;   \|e_n^M\|_{S'(\dot{H}^{-s_c})}+\|e_n^M\|_{S'(L^2)}+\|\nabla e_n^M\|_{S'(L^2)}\leq\varepsilon.
\end{equation}
\ {\bf Claim $2$.} There exist $L>0$ and $S>0$ independent of $M$ such that for any $M$, there exists $n_1=n_1(M)$ such that
 \begin{equation}\label{claim1}
 n>n_1\;\; \Rightarrow\;\;   \|\widetilde{u}_n\|_{S(\dot{H}^{s_c})}\leq L\;\;\textnormal{and}\;\;\|\widetilde{u}_n\|_{L^\infty_tH^1_x}\leq S.
 \end{equation}

  \ By \eqref{CSR2}, there exists $M_1=M_1(\varepsilon)$ such that for each $M>M_1$ there exists $n_2=n_2(M)$ such that
  $$ 
n>n_2\;\; \Rightarrow\;\;  \|U(t)(u_{n,0}-\widetilde{u}_{n,0})\|_{S(\dot{H}^{s_c})}\leq \varepsilon,
  $$
  with $\varepsilon<\varepsilon_1$ as in Proposition \ref{LTP}. Hence, if the two claims hold true, using Proposition \ref{LTP}, for $M$ large enough and $n>\max\{n_0,n_1,n_2\}$, we conclude  $\|u_n\|_{S(\dot{H}^{s_c})}<+\infty$, reaching the desired contradiction. 
  
 \ Up to now, we have reduced the profile expansion to the case where $\psi^1\neq 0$ and $\psi^j= 0$ for all $j\geq 2$. We begin to show the existence of a critical solution. Using the same arguments as before, we can find $\widetilde{\psi}^1$ such that 
\begin{equation*}
u_{n,0}=\textnormal{INLS}(-t_n^1)\widetilde{\psi}^1+\widetilde{W}_n^M,
\end{equation*}   
with
\begin{equation}\label{CSWO11}
M[\widetilde{\psi}^1]=M[\psi^1]\leq 1
\end{equation}
\begin{equation}\label{CSWO221}
 E[\widetilde{\psi}^1]^{s_c}=\left(\frac{1}{2} \|\nabla \psi^1 \|^2_{L^2}\right)^{s_c}\leq \delta_c
  \end{equation}
  \begin{equation}\label{CSWO22}
 \|\nabla \textnormal{INLS}(t)\widetilde{\psi}^1\|^{s_c}_{L^2_x}\|\widetilde{\psi}^1\|^{1-s_c}_{L^2}< \|\nabla Q\|^{s_c}_{L^2}\| Q\|^{1-s_c}_{L^2}
\end{equation}
\begin{equation}\label{CSCP1} 
\lim_{n\rightarrow +\infty} \|U(t)(u_{n,0}-\widetilde{u}_{n,0})\|_{S(\dot{H}^{s_c})}=\lim_{n\rightarrow +\infty}\|U(t)\widetilde{W}_{n}^M\|_{S(\dot{H}^{s_c})}=0.
\end{equation}
   
 \ Set $\widetilde{\psi}^1=u_{c,0}$ and $u_c$ be the global solution\footnote{The global solution is guaranteed by Theorem \ref{TG} and inequalities \eqref{CSWO11}-\eqref{CSWO22}.} to \eqref{INLS} with initial data $\widetilde{\psi}^1$, that is, $u_c(t)=\textnormal{INLS}(t)\widetilde{\psi}^1$. We claim that
\begin{equation}\label{claimfinal}
\|u_c\|_{S(\dot{H}^{s_c})}=+\infty.
\end{equation}    
	
\ Assume, by contradiction, that $\|u_c\|_{S(\dot{H}^{s_c})}<+\infty$. Let $
 \widetilde{u}_n(t)=\textnormal{INLS}(t-t_n^j)\widetilde{\psi}^1,
$ 
then $
\|\widetilde{u}_n(t)\|_{S(\dot{H}^{s_c})}=\|\textnormal{INLS}(t-t_n^j)\widetilde{\psi}^1\|_{S(\dot{H}^{s_c})}=\|u_c\|_{S(\dot{H}^{s_c})}<+\infty$. Also, it follows from \eqref{CSWO11}-\eqref{CSCP1} that 
$$
\sup_{t\in \mathbb{R}}\|\widetilde{u}_n\|_{H^1_x}=\sup_{t\in \mathbb{R}}\|u_c\|_{H^1_x}<+\infty\;\;\;\textnormal{and}\;\;\;\|U(t)(u_{n,0}-\widetilde{u}_{n,0})\|_{S(\dot{H}^{s_c})}\leq \varepsilon,
$$
for $n$ large enough. Therefore, using the long time perturbation theory (Proposition \ref{LTP}) with $e=0$, we deduce $\|u_n\|_{S(\dot{H}^{s_c})}<+\infty$, which is a contradiction with \eqref{un}. 
	
\ On the other hand, the relation \eqref{claimfinal} yields $E[u_c]^{s_c}M[u_c]^{1-s_c}=\delta_c$ (see \eqref{deltac}). Thus, we conclude from \eqref{CSWO11} and \eqref{CSWO221} that
$$
M[u_c]=1\;\;\;\;\textnormal{and}\;\;\;\;E[u_c]^{s_c}=\delta_c.
$$ 
Also note that \eqref{CSWO22} implies (iii) in the statement of the Proposition \ref{ECS}.

\subsection{Proof of Claim 1 and 2} In this subsection we complete the proof of Proposition \ref{ECS}. We show Claims $1$ and $2$ (see \eqref{claim1} and \eqref{claim2}). To this end, we use the same admissible pairs used by the second author in \cite{CARLOS} to prove global well-posedness.
\begin{equation}\label{PHsA1}  
\widehat{q}=\frac{4\alpha(\alpha+2-\theta)}{\alpha(N\alpha+2b)-\theta(N\alpha-        4+2b)},\;\;\;\widehat{r}\;=\;\frac{N\alpha(\alpha+2-\theta)}{\alpha(N-b)-\theta(2-b)},
\end{equation}
and
\begin{equation}\label{PHsA2}
\widetilde{a}\;=\;\frac{2\alpha(\alpha+2-\theta)}{\alpha[N(\alpha+1-\theta)-2+2b]-(4-2b)(1-\theta)},\;\;\;  \widehat{a}=\frac{2\alpha(\alpha+2-\theta)}{4-2b-(N-2)\alpha}.
\end{equation}
We have that $(\widehat{q},\widehat{r})$ is $L^2$-admissible, $(\widehat{a},\widehat{r})$ is $\dot{H}^{s_c}$-admissible and $(\widetilde{a},\widehat{r})$ is $\dot{H}^{-s_c}$-admissible (for more details see \cite[Subsection $4.2$]{CARLOS}).\\

\ {\bf Proof of Claim $1$.} First, we prove that for each $M$ and $\varepsilon>0$, there exists $n_0=n_0(M,\varepsilon)$ such that $\|e_n^M\|_{S'(\dot{H}^{-s_c})}< \frac{\varepsilon}{3}$. It follows from \eqref{CSR1} and \eqref{EIerror} that
\begin{equation}\label{ec21} 
\|e_n^M\|_{S'(\dot{H}^{-s_c})}\leq C_{\alpha,M}\sum_{j=1}^{M}\sum_{1\leq j\neq k\leq M}  \left\||x|^{-b}|v^k|^\alpha|v^j|\right\|_{L^{\widetilde{a}'}_tL^{\widehat{r}'}_x}.
\end{equation}
We claim that the norm in the right hand side of \eqref{ec21} goes to $0$ as $n\rightarrow +\infty$. Indeed, using the relation $(4.13)$ of \cite{CARLOS}, with $s=1$, we get
 \begin{align}\label{ec22}
  \left\||x|^{-b}|v^k|^\alpha|v^j|\right\|_{L^{\widetilde{a}'}_tL^{\widehat{r}'}_x} \leq& c \|v^k\|^{\theta}_{L^\infty_tH^1_x} \left\|\|v^k(t-t_n^k)\|^{\alpha-\theta}_{L_x^{\widehat{r}}} \|v^j(t-t_n^j)\|_{L^{\widehat{r}}_x}\right\|_{L^{\widetilde{a}'}_t}.   
\end{align}
Fix $1\leq j\neq k\leq M$. Combining \eqref{CSWO1} and  \eqref{CSWO2} we deduce $\|v^k\|_{H^1_x}<+\infty$ and by \eqref{CSCP} $v^j$, $v^k\in S(\dot{H^{s_c}})$, thus we can approximate $v^j$ by functions of $C_0^\infty(\mathbb{R}^{N+1})$. Define $g_n(t)= \|v^k(t)\|^{\alpha-\theta}_{L_x^{\widehat{r}}} \|v^j(t-(t_n^j-t_n^k))\|_{L^{\widehat{r}}_x}$, we have $g_n\in L^{\widetilde{a}'}_t$. Indeed, it follows from H\"older inequality (since $\frac{1}{\widetilde{a}'}=\frac{\alpha-\theta}{\widehat{a}}+\frac{1}{\widehat{a}}$) that 
$$
\|g_n\|_{L^{\widetilde{a}'}_t}\leq \|v^k\|^{\alpha-\theta}_{L^{\widehat{a}}_t L_x^{\widehat{r}}} \|v^j\|_{L^{\widehat{a}}_tL^{\widehat{r}}_x}\leq \|v^k\|^{\alpha-\theta}_{S(\dot{H}^{s_c})} \|v^j\|_{S(\dot{H}^{s_c})}<+\infty.
$$
Moreover, by \eqref{PD} we obtain $g_n(t)\rightarrow 0$ as $n\rightarrow +\infty$. On the other hand, $|g_n(t)|\leq  KI_{supp(v^j)}\|v^k(t)\|^{\alpha-\theta}_{L_x^{\widehat{r}}}\equiv g(t)$ for all $n$, where $K>0$ and $I_{supp(v^j)}$ is the characteristic function of $supp(v^j)$. Similarly as (i), we get
$$
\|g\|_{L^{\widetilde{a}'}_t}\leq \|v^k\|^{\alpha-\theta}_{L^{\widehat{a}}_t L_x^{\widehat{r}}}\|I_{supp(v^j)}\|_{L^{\widehat{a}}_t L_x^{\widehat{r}}}<+\infty.
$$
That is, $g\in L^{\widetilde{a}'}_t$. Then, the Dominated Convergence Theorem yields $\|g_n\|_{L^{\widetilde{a}'}_t}\rightarrow 0$ as $n\rightarrow +\infty$, which implies by \eqref{ec22} the first estimate. 
 
\ Next, using the same argument as before, we show $\|e_n^M\|_{S'(L^2)}<\frac{\varepsilon}{3}$. Indeed, again the elementary inequality \eqref{EIerror} yields
\begin{equation*}
\|e_n^M\|_{S'(L^2)}\leq C_{\alpha,M}\sum_{j=1}^{M}\sum_{1\leq j\neq k\leq M}  \left\||x|^{-b}|v^k|^\alpha|v^j|\right\|_{L^{\widehat{q}'}_tL^{\widehat{r}'}_x}.
\end{equation*}   
We also obtain (see proof of \cite[Lemma $4.2$ with $s=1$]{CARLOS}) 
\begin{eqnarray*}
\left\|  |x|^{-b}|v^k|^\alpha|v^j \right \|_{L_t^{\widehat{q}'}L^{\widehat{r}'}_x}&\leq&  c \|v^k\|^{\theta}_{L^\infty_tH^1_x}\left \| \|v^k(t-t_n^k)\|^{\alpha-\theta}_{ L_x^{\widehat{r}}} \|v^j(t-t_n^j)\|_{L^{\widehat{r}}_x}\right\|_{L_t^{\widehat{q}'}}\\
  &\leq & c\|v^k\|^{\theta}_{L^\infty_tH^1_x} \|v^k\|^{\alpha-\theta}_{L^{\widehat{a}}_t L_x^{\widehat{r}}} \|v^j\|_{L_t^{\widehat{q}}L^{\widehat{r}}_x}\\
 &\leq & c\|v^k\|^{\theta}_{L^\infty_tH^1_x}  \|v^k\|^{\alpha-\theta}_{S(\dot{H}^{s_c})}  \|v^j\|_{S(L^2)}.
\end{eqnarray*}
This implies that the right hand side of the last inequality is finite (since\footnote{Note that, $v^j\in S(\dot{H}^{s_c})$ implies that $\|v^j\|_{S(L^2)}$ and $\|\nabla v^j\|_{S(L^2)}<+\infty$, by \eqref{SCATTER1}.} $\|v^j\|_{S(L^2)}$ and $\|\nabla v^j\|_{S(L^2)}<+\infty$ ) and so
$$
 \left \| \|v^k(t-t_n^k) \|^{\alpha-\theta}_{ L_x^{\widehat{r}}} \|v^j(t-t_n^j) \|_{L^{\widehat{r}}_x}\right\|_{L_t^{\widehat{q}'}}\rightarrow 0, 
  $$
  as $n\rightarrow +\infty$, which lead to 
$
\left\|  |x|^{-b}|v^k|^\alpha|v^j \right \|_{L_t^{\widehat{q}'}L^{\widehat{r}'}_x}\rightarrow 0.
$
  
\ Finally, we prove $\|\nabla e_n^M\|_{S'(L^2)}<\frac{\varepsilon}{3}$. Observe that
\begin{eqnarray}\label{ec23}
\nabla e_n^M&=&\nabla (|x|^{-b})\left( f(\widetilde{u}_n)-\sum_{j=1}^M f(v^j)  \right)+|x|^{-b}\nabla \left(  f(\widetilde{u}_n)-\sum_{j=1}^M f(v^j)\right)\nonumber \\
&\equiv& R^1_n+R^2_n,
\end{eqnarray}
where $f(v)=|v|^\alpha v$. We start by considering $R^1_n$. Applying \eqref{EIerror} we estimate 
$$
\| R^1_n\|_{S'(L^2)}\leq c\; C_{\alpha,M}\sum_{j=1}^{M}\sum_{1\leq j\neq k\leq M}  \left\||x|^{-b-1}|v^k|^\alpha|v^j|\right\|_{S'(L^2)}
$$
and by Remark \ref{RSglobal} we deduce that $\left\||x|^{-b-1}|v^k|^\alpha|v^j|\right\|_{S'(L^2)}$ is finite, then by the same argument as before we obtain 
$$
 \left\||x|^{-b-1}|v^k(t-t_n^k)|^\alpha|v^j(t-t_n^j)|\right\|_{S'(L^2)}\rightarrow 0\;\;\textnormal{as}\;\;n\rightarrow+\infty.
$$
Therefore, the last two relations yield $\|R^1_n\|_{S'(L^2)}\rightarrow 0$ as $n\rightarrow+\infty$.

\ On the other hand, note that
\begin{eqnarray}\label{ec231}
\nabla (  f(\widetilde{u}_n)-\sum_{j=1}^M f(v^j))&=&f'(\widetilde{u}_n)\nabla \widetilde{u}_n-\sum_{j=1}^M f'(v^j)\nabla v^j\nonumber\\
&=&\sum_{j=1}^M (f'(\widetilde{u}_n)- f'(v^j))\nabla v^j.
\end{eqnarray}
In view of (by Remark \ref{nonlinerity})
$$
|f'(\widetilde{u}_n)- f'(v^j)|\leq C_{\alpha,M}\sum_{1\leq k\neq j\leq M}|v^k|(|v^j|^{\alpha-1}+|v^k|^{\alpha-1})\;\;\;\textnormal{if}\;\;\;\alpha>1
$$
and
$$
|f'(\widetilde{u}_n)- f'(v^j)|\leq C_{\alpha,M}\sum_{1\leq k\neq j\leq M}|v^k|^\alpha\;\;\;\textnormal{if}\;\;\;\alpha\leq 1,
$$
we have using the last two relations together with \eqref{ec23} and \eqref{ec231} 
$$
\|R_n^2\|_{S'(L^2)}\lesssim \sum_{j=1}^M\sum_{1\leq k\neq j\leq M} \left\||x|^{-b}  |v^k|(|v^j|^{\alpha-1}+|v^k|^{\alpha-1})|\nabla v^j|\right\|_{S'(L^2)}\;\;\;\textnormal{if}\;\;\;\alpha>1,
$$
and
$$
\|R_n^2\|_{S'(L^2)}\lesssim \sum_{j=1}^M\sum_{1\leq k\neq j\leq M} \left\||x|^{-b}  |v^k|^\alpha|\nabla v^j|\right\|_{S'(L^2)}\;\;\;\textnormal{if}\;\;\;\alpha\leq 1.
$$
Therefore, from Lemma \ref{LG2} (see also Remark \ref{RGP}) we have that the right hand side of the last two inequalities are finite quantities and, by an analogous argument as before, we conclude that 
$$
\|R_n^2\|_{S'(L^2)}\rightarrow 0\;\;\;\textnormal{as}\;\;\;n\rightarrow+\infty.
$$
This completes the proof of Claim $1$.\\
   
\ {\bf Proof of Claim $2.$} To this end, we first prove that $\|\widetilde{u}_n\|_{L^\infty_tH^1_x}$ and $\|\widetilde{u}_n\|_{L^\gamma_tL^\gamma_x}$ are bounded quantities where $\gamma=\frac{2(N+2)}{N}$. Indeed, we already know (see \eqref{PCS41} and \eqref{Sumpsij}) that there exists $C_0$ such that
$$
\sum_{j=1}^{\infty}\|\psi^j\|^2_{H^1_x}\leq C_0,
$$
then choosing $M_0\in \mathbb{N}$ large enough such that
\begin{equation}\label{SP} 
\sum_{j=M_0}^{\infty}\|\psi^j\|^2_{H^1_x}\leq \frac{\delta}{2},
\end{equation}
where $\delta>0$ is a sufficiently small.\\
Fix $M\geq M_0$. From \eqref{CSWO3}, there exists $n_1(M)\in \mathbb{N}$ where for all $n> n_1(M)$, it follows that (using \eqref{SP})
$$
\sum_{j=M_0}^{M}\|\textnormal{INLS}(-t_n^j)\widetilde{\psi}^j\|^2_{H^1_x}\leq \delta,
$$    
This is equivalent to 
\begin{equation}\label{claim11}
\sum_{j=M_0}^{M}\|v^j(-t_n^j)\|^2_{H^1_x}\leq \delta.
\end{equation}
Then, by the Small Data Theory (Proposition \ref{GWPH1})
$$
\sum_{j=M_0}^{M}\|v^j(t-t_n^j)\|^2_{L_t^{\infty}H^1_x}\leq c\delta\;\;\textnormal{for}\;n\geq n_1(M).
$$
Observe that,  
$$
\left\|\sum_{j=M_0}^{M} v^j(t-t_n^j)\right\|^2_{H^1_x}=\sum_{j=M_0}^{M}\|v^j(t-t_n^j)\|^2_{H_x^1}+2\sum_{M_0\leq l \neq k\leq M}\langle v^l(t-t_n^l),v^k(t-t_n^k)\rangle_{H^1_x},
$$ 
so, for $l\neq k$ we deduce from $\eqref{PD}$ that (see \cite[Corollary $4.4$]{JIANCAZENAVE} for more details)
$$ \sup_{t\in\mathbb{R}} |\langle v^l(t-t_n^l),v^k(t-t_n^k)\rangle_{H^1_x}|\rightarrow 0\;\;\textnormal{as}\;\;n\rightarrow +\infty.
$$
In view of $\|v^j\|_{L^{\infty}_tH_x^1}$ is bounded (see \eqref{CSWO1} - \eqref{CSWO2}), by definition of $\widetilde{u}_n$ there exists $S>0$ (independent of $M$) such that
 \begin{equation}\label{claim12}
 \sup_{t\in\mathbb{R}}\|\widetilde{u}_n\|^2_{H^1_x}\leq S \;\,\textnormal{for}\;\;n>n_1(M).
\end{equation}
   
   \ We now show $\|\widetilde{u}_n\|_{L^\gamma_tL^\gamma_x}\leq L_1$. Using again \eqref{claim11} with $\delta$ small enough and the Small Data Theory (noting that $(\gamma,\gamma)$ is $L^2$-admissible and $\gamma >2$), we deduce 
\begin{equation}\label{claim111}
 \sum_{j=M_0}^{M}\|v^j(t-t_n^j)\|^{\gamma}_{L^\gamma_tL^\gamma_x}\leq c  \sum_{j=M_0}^{M}\|v^j(-t_n^j)\|^{\gamma}_{H^1_x}\leq c  \sum_{j=M_0}^{M}\|v^j(-t_n^j)\|^2_{H^1_x}\leq c\delta,
\end{equation}
for $n\geq n_1(M)$.

 On the other hand, since \eqref{FI} we have that 
  $$
   \left\|\sum_{j=M_0}^{M} v^j(t-t_n^j)\right\|^\gamma_{L^\gamma_tL^\gamma_x}\leq \sum_{j=M_0}^{M}\|v^j\|^\gamma_{L^\gamma_tL^\gamma_x}+C_M\sum_{M_0\leq j \neq k\leq M}\int_{\mathbb{R}^{N+1}} |v^j||v^k||v^k|^{\gamma-2}
 $$ 
for all $M>M_0$. If for a given $j$ such that $M_0\leq j\neq k\leq M$, it follows from H\"older inequality that 
\begin{eqnarray}\label{claim112}
 \int_{\mathbb{R}^{N+1}} |v^j||v^k||v^k|^{\gamma-2}& \leq &\|v^k(t-t_n^k)\|_{L^\gamma_{t}L^\gamma_x}\left( \int_{\mathbb{R}^{N+1}} |v^j|^{\frac{\gamma}{2}}|v^k|^{\frac{\gamma}{2}} \right)^{\frac{2}{\gamma}} \nonumber \\
 &\leq& c \|v^j(-t_n^j)\|_{H^1_x}\left( \int_{\mathbb{R}^{N+1}} |v^j|^{\frac{\gamma}{2}}|v^k|^{\frac{\gamma}{2}} \right)^{\frac{2}{\gamma}}.
\end{eqnarray}
In view of $v^j$ and $v^k\in L^\gamma_tL^\gamma_x$ we get that the right hand side of \eqref{claim112} is bounded and so by similar arguments as in the previous claim, we conclude from \eqref{PD} that the integral in the right hand side of the previous inequality goes to $0$ as $n\rightarrow +\infty$. This implies that there exists $L_1$ (independent of $M$) such that (\eqref{claim111})
\begin{equation}\label{claim113}
 \|\widetilde{u}_n\|_{L^\gamma_tL^\gamma_x} \leq \sum_{j=1}^{M_0}\|v^j\|_{L^\gamma_tL^\gamma_x}+\left\|\sum_{j=M_0}^{M} v^j\right\|_{L^\gamma_tL^\gamma_x}\leq L_1\;\;\;\textnormal{for}\;n\geq n_1(M).
\end{equation}

\ To complete the proof of the Claim $2$ we will prove the following inequalities 
\begin{equation}\label{ineq1}
\left\|  |x|^{-b}  |\widetilde{u}_n|^\alpha \widetilde{u}_n \right\|_{L^{\bar{a}'}_tL_x^{\bar{r}'}} \leq c\|\widetilde{u}_n\|^\theta_{L^\infty_tH^1_x}\|\widetilde{u}_n\|^{\alpha-\theta+1}_{L^a_tL^r_x}
\end{equation}
and 
\begin{equation}\label{ineq2}
\|\widetilde{u}_n\|_{L^a_tL^r_x}\leq \|\widetilde{u}_n\|^{1-\frac{\gamma}{a}}_{L^\infty_tH^1_x}\|\widetilde{u}_n\|^{\frac{\gamma}{a}}_{L^\gamma_tL^\gamma_x},
\end{equation}
where $\theta\in (0,\alpha)$ is a small enough and the pairs $(\bar{a},\bar{r})$ and $(a,r)$ are $\dot{H}^{-s_c}$-admissible and $\dot{H}^{s_c}$-admissible, respectively.  

\ Observe that, combining \eqref{claim12} and \eqref{claim113} we deduce from \eqref{ineq2} that
$$
\|\widetilde{u}_n\|_{L^a_tL^r_x}\leq S^{1-\frac{\gamma}{a}}L_1^{\frac{\gamma}{a}}=L_2,\;\;\;\textnormal{for}\;n\geq n_1(M).
$$
Then, since $\widetilde{u}_n$ satisfies the perturbed equation \eqref{widetildeun} we can apply the Strichartz estimates (Lemma \ref{Lemma-Str}) and \eqref{ineq1} to the integral formulation and conclude (using also Claim $1$)
\begin{eqnarray*}
\|\widetilde{u}_n\|_{S(\dot{H}^{s_c})}&\leq &c\|\widetilde{u}_{n,0}\|_{H^1_x}+c\left\|  |x|^{-b}  |\widetilde{u}_n|^\alpha \widetilde{u}_n \right\|_{L^{\bar{a}'}_tL_x^{\bar{r}'}}+\|e^M_n\|_{S'(\dot{H}^{-s_c})}\\
&\leq & cS+cL_2+\varepsilon = L,
\end{eqnarray*}
for $n\geq n_1(M)$, which completes the proof of the Claim $2$.  

\ It remains to prove the inequalities \eqref{ineq1} and \eqref{ineq2}. Indeed, we divide in two cases: $N\geq 3$ and $N=2$, since we will make use of the Sobolev embeddings in Lemma \ref{SI}.\\
{\bf Case $N\geq 3$}: We use the following numbers:
\begin{equation}\label{claim21}
a=\frac{4\alpha(N+2)}{ND}\hspace{1.5cm}\,r=\frac{2\alpha N(N+2)}{(4-2b)(N+2)-ND}
\end{equation}
\begin{equation}\label{claim22}
\bar{a}=\frac{4\alpha(N+2)}{4\alpha(N+2)-(\alpha+1-\theta)ND}\hspace{0.2cm}\
\end{equation}
and
\begin{equation}\label{claim23}
\bar{r}=\frac{2\alpha N(N+2)}{2(N+2)\left(\alpha(N-2)-(2-b)\right)+ND(\alpha+1-\theta)},
\end{equation}
where $D=4-2b-\alpha(N-2)$ and $\theta\in (0,\alpha)$ to be chosen below. 


\ It is easy to see that $(a,r)$ is $\dot{H}^{s_c}$-admissible and $(\bar{a},\bar{r})$ is $\dot{H}^{-s_c}$-admissible. In Appendix, we will verify the conditions of admissible pair.   
 
\ We first show the inequality \eqref{ineq2}. Indeed, by interpolation we have
$$
\|\widetilde{u}_n\|_{L^a_tL^r_x}\leq \|\widetilde{u}_n\|^{1-\frac{\gamma}{a}}_{L^\infty_tL^p_x}\|\widetilde{u}_n\|^{\frac{\gamma}{a}}_{L^\gamma_tL^\gamma_x},
$$
where
$$
\frac{1}{r}=\left(1-\frac{\gamma}{a}\right)\left(\frac{1}{p}\right)+\frac{1}{a},
$$
which is equivalent to (recall that $\gamma=\frac{2(N+2)}{N}$)
\begin{eqnarray*}
\left(1-\frac{\gamma}{a}\right)\left(\frac{1}{p}\right)&=&\frac{1}{r}-\frac{1}{a}\\
\frac{2\alpha-D}{p}&=&\frac{2(4-2b)-ND}{2N}\\
p&=& \frac{2N}{N-2}.
\end{eqnarray*}
Hence, since $H^1 \hookrightarrow L^{\frac{2N}{N-2}}$ (see inequality \eqref{SEI22} with $s=1$) we obtain the desired result. 

\ On the other hand, the proof of inequality \eqref{ineq1} follows from similar ideas as in Lemma \ref{LG3}. We divide the estimate in $B$ and $B^C$. Let $A\subset \mathbb{R}^N$. From the H\"older inequality we deduce
\begin{eqnarray*}
\left\|  |x|^{-b}  |\widetilde{u}_n|^\alpha \widetilde{u}_n \right\|_{L^{\bar{a}'}_tL_x^{\bar{r}'}(A)}&\leq & \left\| \| |x|^{-b}\|_{L^d(A)}  \|\widetilde{u}_n\|^{\theta}_{L_x^{\theta r_1}}\|\widetilde{u}_n\|^{\alpha+1-\theta}_{L_x^{(\alpha+1-\theta)r_2}}\right\|_{L^{\bar{a}'}_t}     \\
& \leq& \| |x|^{-b}\|_{L^d(A)} \|\widetilde{u}_n\|^{\theta}_{L_x^{\theta r_1}}   \|\widetilde{u}_n\|^{\alpha+1-\theta}_{L_t^{(\alpha+1-\theta)\bar{a}'}L_x^{(\alpha+1-\theta)r_2}}\\\
&=& \| |x|^{-b}\|_{L^d(A)} \|\widetilde{u}_n\|^{\theta}_{L_x^{\theta r_1}}\|\widetilde{u}_n\|^{\alpha-\theta+1}_{L^a_tL^r_x},
\end{eqnarray*}
where
$$
\frac{1}{\bar{r}'}=\frac{1}{d}+\frac{1}{r_1}+\frac{1}{r_2}\,\hspace{0.5cm}\,r=(\alpha+1-\theta)r_2\,\hspace{0.5cm}\,a=(\alpha+1-\theta)\bar{a}'.
$$
Using the values of $a$ and $\bar{a}$ above defined, it is easy to check $a=(\alpha+1-\theta)\bar{a}'$. Moreover, to show that $\| |x|^{-b}\|_{L^d(A)}$ is a bounded quantity we need $\frac{N}{d}-b>0$ if $A=B$ and $\frac{N}{d}-b<0$ if $A=B^C$, see Remark \ref{RIxb}. Indeed, the last relation implies
\begin{eqnarray*}
\frac{N}{d}-b&=&N-b-\frac{N}{r_1}-\frac{N}{\bar{r}}-\frac{N(\alpha+1-\theta)}{r}\\
&=&\frac{\theta(2-b)}{\alpha}-\frac{N}{r_1}.
\end{eqnarray*}
Choosing $\theta r_1=2$ we have $\frac{N}{d}-b=-\theta s_c<0$, so $|x|^{-b}\in L^d(B^C)$ and if $\theta r_1=\frac{2N}{N-2}$ then $\frac{N}{d}-b=\theta(1-s_c)>0$, i.e., $|x|^{-b}\in L^d(B)$. Therefore, since in both cases $\theta r_1\in \left[2, \frac{2N}{N-2}\right]$, by the Sobolev embedding \eqref{SEI22} we complete the proof of \eqref{ineq1}.\\
{\bf Case $N=2$.} We start by defining the following numbers.
\begin{equation}\label{C2N21}
a=\frac{2\alpha(\alpha+1-\theta)}{2-b+\varepsilon}\hspace{1.5cm}\,r=\frac{2\alpha(\alpha+1-\theta)}{(2-b)(\alpha-\theta)-\varepsilon} 
\end{equation}
and
\begin{equation}\label{C2N22}
\bar{a}=\frac{2\alpha}{2\alpha-(2-b)-\varepsilon}\hspace{1.5cm}\,\bar{r}=\frac{2\alpha}{\varepsilon},
\end{equation}
where $\theta\in (0,\alpha)$ and $\varepsilon>0$ are sufficiently small numbers. A simple computation shows that $(a,r)$ is $\dot{H}^{s_c}$-admissible and $(\bar{a},\bar{r})$ is $\dot{H}^{-s_c}$ admissible.

\ The interpolation inequality implies that (in this case $\gamma=4$)
$$
\|\widetilde{u}_n\|_{L^a_tL^r_x}\leq \|\widetilde{u}_n\|^{1-\frac{\gamma}{a}}_{L^\infty_tL^p_x}\|\widetilde{u}_n\|^{\frac{\gamma}{a}}_{L^\gamma_tL^\gamma_x},
$$
where
$$
\frac{1}{r}=\left(1-\frac{\gamma}{a}\right)\left(\frac{1}{p}\right)+\frac{1}{a}.
$$
This is equivalent to 
\begin{eqnarray*}
\left(1-\frac{4}{a}\right)\left(\frac{2}{p}\right)&=&\frac{2}{r}-\frac{2}{a}\\
&=&\frac{2-b}{\alpha}-\frac{4}{a}\\
&=&\frac{(2-b)(\alpha-\theta+1)-2(2-b-\varepsilon)}{\alpha(\alpha-\theta+1)}.
\end{eqnarray*}
Thus
$$
p=2\frac{\alpha(\alpha-\theta+1)-2[(2-b)-\varepsilon]}{(2-b)(\alpha+1-\theta)-2[(2-b)-\varepsilon]}.
$$
Since we are assuming $\alpha>2-b$ we have $p>2$, thus by the Sobolev embedding $H^1 \hookrightarrow L^p$ (see \eqref{SEI1} with $N=2$) the inequality \eqref{ineq2} holds. To  show the inequality \eqref{ineq1} we use the same argument as the previous case, that is
\begin{eqnarray*}
\left\|  |x|^{-b}  |\widetilde{u}_n|^\alpha \widetilde{u}_n \right\|_{L^{\bar{a}'}_tL_x^{\bar{r}'}(A)}
& \leq& \| |x|^{-b}\|_{L^d(A)} \|\widetilde{u}_n\|^{\theta}_{L_x^{\theta r_1}}   \|\widetilde{u}_n\|^{\alpha+1-\theta}_{L_t^{(\alpha+1-\theta)\bar{a}'}L_x^{(\alpha+1-\theta)r_2}}\\\
&=& \| |x|^{-b}\|_{L^d(A)} \|\widetilde{u}_n\|^{\theta}_{L_x^{\theta r_1}}\|\widetilde{u}_n\|^{\alpha-\theta+1}_{L^a_tL^r_x},
\end{eqnarray*}
where $A=B$ or $B^C$ and
$$
\frac{1}{\bar{r}'}=\frac{1}{d}+\frac{1}{r_1}+\frac{1}{r_2}\,\hspace{0.5cm}\,r=(\alpha+1-\theta)r_2\,\hspace{0.5cm}\,a=(\alpha+1-\theta)\bar{a}'.
$$
Moreover, we obtain
\begin{eqnarray*}
\frac{2}{d}-b&=&2-b-\frac{2}{r_1}-\frac{2}{\bar{r}}-\frac{2(\alpha+1-\theta)}{r}\\
&=&\frac{\theta(2-b)}{\alpha}-\frac{2}{r_1}.
\end{eqnarray*}
If we choose $\theta r_1\in \left(2,\frac{2\alpha}{2-b} \right)$ then $\frac{2}{d}-b<0$ (so $|x|^{-b}\in L^d(B^C)$) and if $\theta r_1\in \left(\frac{2\alpha}{2-b},+\infty \right)$ we have $\frac{2}{d}-b<0$ (so $|x|^{-b}\in L^d(B)$). Therefore $|x|^{-b}\in L^d(A)$ and so by the Sobolev inequality \eqref{SEI1} with $s=1$, we complete the proof of the inequality \eqref{ineq1}.
\end{proof}
\end{proposition}
\begin{remark}
To show that $\bar{r}$ defined in \eqref{claim22} satisfies the condition \eqref{H-s}, that is $\frac{2N}{N-2s_c}<\bar{r}< \frac{2N}{N-2}$, we need the assumptions $b<\min\{\frac{N}{3},1\}$ and $\alpha<2_*$. Indeed $\bar{r}<\frac{2N}{N-2}$ is equivalent to  
$$
\alpha(N+2)(N-2)<2(N+2)\left( \alpha(N-2)-(2-b)  \right)+ND(\alpha+1-\theta)\Leftrightarrow 
$$
$$
(N+2)D<ND(\alpha+1-\theta)\Leftrightarrow N(\alpha+1-\theta)>N+2\;\Leftrightarrow \;\alpha N-2-\theta N>0.
$$
Since $\alpha>(4-2b)/N$ we have $\alpha N-2-\theta N>4-2b-2-\theta N=2(1-b)-\theta N$ and this is positive choosing $\theta<\frac{2(1-b)}{N}$ (here we use the condition $0<b<\min\{\frac{N}{3},1\}$ to guarantee that $\theta$ can be chosen to be a positive number). Therefore, since $\alpha N-2-\theta N>0$ one gets $\bar{r}<\frac{2N}{N-2}$. On the other hand, $\bar{r}>\frac{2N}{N-2s_c}=\frac{N\alpha}{2-b}$ is equivalent to 
$$
(N+2)(4-2b)>2(N+2)\left( \alpha(N-2)-(2-b) \right) +ND(\alpha+1-\theta)\Leftrightarrow
$$
$$
2(N+2)D>ND(\alpha +1-\theta)\;\Leftrightarrow\;\alpha <\frac{N+4+\theta N}{N}. 
$$
Since $\alpha <2_*$ (defined in \eqref{def2_*}) we need to verify that $\frac{4-2b}{N-2}\leq \frac{N+4+\theta N}{N}$ for $N\geq 4$ and $3-2b\leq \frac{7+3\theta }{3}$ for $N=3$. The first inequality is equivalent to $N(4-2b)\leq (N+4+\theta N)(N-2)$ and this is always true since $N\geq 4$. The second case is also true choosing\footnote{In the particular case when $N=3$, we need to choose $\theta>0$ such that $\max \left\{ 0,\frac{2(1-3b)}{3}\right\}<\theta<\frac{2(1-b)}{3}$, since also need $\theta<\frac{2(1-b)}{N}$ to obtain $\bar{r}<\frac{2N}{N-2}$.} $\theta > \max \left\{ 0,\frac{2(1-3b)}{3}\right\}$.\\
\end{remark}

\ In the next proposition, we prove the precompactness of the flow associated to the critical solution $u_c$. 
\begin{proposition}\label{PSC}{\bf (Precompactness of the flow of the critical solution)} Let $u_c$ be as in Proposition \ref{ECS} and define
 $$
 K=\{u_c(t)\;:\;t\in[0,+\infty)\}\subset  H^1.
 $$
 Then $K$ is precompact in $H^1(\mathbb{R}^N)$.
\begin{proof} The proof is similar to that of Proposition $6.5$ in \cite{paper2} (replacing $3$ by $N$). So, we only give the main steps.
 
 Let $\{t_n\}\subseteq [0,+\infty )$ a sequence of times and $\phi_n=u_c(t_n)$ be a uniformly bounded sequence in $H^1(\mathbb{R}^N)$. We need to show that $u_c(t_n)$ has a subsequence converging in $H^1(\mathbb{R}^N)$. The result is clear if $\{t_n\}$ is bounded. Now assume that $t_n\rightarrow +\infty$. The linear profile expansion (Proposition \ref{LPD}) and the energy Pythagorean expansion (Proposition \ref{EPE}) yield the existence of profiles $\psi^j$ and a remainder $W_n^M$ such that
 $$
 u_c(t_n)=\sum_{j=1}^{M}U(-t_n^j)\psi^j+W_n^M
 $$	
and
\begin{equation}\label{ECCS}
\sum_{j=1}^{M}\lim_{n\rightarrow +\infty}E[U(-t_n^j)\psi^j]+\lim_{n\rightarrow +\infty} E[W_n^M]=E[u_c]=\delta_c,	
\end{equation}
which implies that\footnote{Since each energy in \eqref{ECCS} is nonnegative by Lemma \ref{LGS} (i).} $\lim_{n\rightarrow +\infty}E[U(-t_n^j)\psi^j]\leq \delta_c\;\;\;\;\forall\;j$.
Moreover, by \eqref{PDNHs} with $s=0$ we obtain
\begin{equation}\label{MCCS} 
\sum_{j=1}^{M}M[\psi^j]+\lim_{n\rightarrow +\infty}M[W_n^M]=M[u_c]=1,
\end{equation}
by Proposition \ref{ECS} (i).
 	
\ If more than one $\psi^j \neq 0$, similar to the proof in Proposition \ref{ECS}, we have a contradiction
with the fact that $\|u_c\|_{S(\dot{H}^{s_c})}=+\infty$. Thus, we address the case that only $\psi^j=0$ for all $j\geq 2$, and so
\begin{equation}\label{CCS1}
u_c(t_n)=U(-t_n^1)\psi^1+W_n^M.
\end{equation}
Also as in the proof of Proposition \ref{ECS}, we have
 \begin{equation}\label{MECS}
 M[\psi^1] =M[u_c]=1\;\;\;\textnormal{and}\;\;\;\lim_{n\rightarrow +\infty}E[U(-t_n^1)\psi^1]=\delta_c,	
 \end{equation}
 and using \eqref{ECCS}, \eqref{MCCS} together with \eqref{MECS}, we deduce that 
 \begin{equation}
 \lim_{n\rightarrow +\infty}M[W_n^M]=0\;\;\;\textnormal{and}\;\;\;\lim_{n\rightarrow +\infty}E[W_n^M]=0.
 \end{equation}	
By Lemma \ref{LGS} (i) we conclude that
\begin{equation}\label{ERCS}
 \lim_{n\rightarrow +\infty}\|W_n^M\|_{H^1}=0.
\end{equation} 
	
\ If $t^1_n$ converges to some finite $t^*$, it is easy to see that $u_c(t_n)$ converges in $H^1(\mathbb{R}^N)$, concluding the proof. 
  
\ Assume by contradiction that $|t^1_n|\rightarrow +\infty$, then we have two cases to consider. If $t^1_n\rightarrow -\infty$, by \eqref{CCS1}  
$$
\|U(t)u_c(t_n)\|_{S(\dot{H}^{s_c};[0,+\infty))}\leq\|U(t-t_n^1)\psi^1\|_{S(\dot{H}^{s_c};[0,+\infty))}+\|U(t)W_n^M\|_{S(\dot{H}^{s_c};[0,+\infty))}.
$$	
On the other hand, we also obtain
 \begin{equation*}\label{PPCS1}
  \|U(t-t_n^1)\psi^1\|_{S(\dot{H}^{s_c};[0,+\infty))}\leq \|U(t)\psi^1\|_{S(\dot{H}^{s_c};[-t_n^j,+\infty))}\leq \frac{1}{2}\delta,
\end{equation*}
and (given $\delta>0$ for $n, M$ large and using \eqref{SE2} \eqref{ERCS}) $\|U(t)W_n^M\|_{S(\dot{H}^{s_c})}\leq \frac{1}{2}\delta$. So
 $$
 \|U(t)u_c(t_n)\|_{S(\dot{H}^{s_c};[0,+\infty))}\leq \delta.
 $$ 
 Therefore, choosing $\delta>0$ sufficiently small, by the small data theory (Proposition \ref{GWPH1}) we get that $ \|u_c\|_{S(\dot{H}^{s_c})}\leq 2\delta$,
 which is a contradiction with Proposition \ref{ECS}(iv).
 
 \ Similarly, we have a contradiction when $t^1_n\rightarrow -\infty$.
\end{proof}
\end{proposition}

\section{Rigidity theorem}
The goal in this section is a rigidity theorem, which implies that the critical solution $u_c$ constructed in Section \ref{CCS} must be identically zero and so reaching a contradiction in view of Proposition \ref{ECS} (iv). To this end, we need the following results. 
\begin{proposition}\label{PFIUL} Let $u$ be a solution of \eqref{INLS} such that
$$
K=\{u(t)\;:\;t\in[0,+\infty)\}
$$
is precompact in $H^1(\mathbb{R}^N)$. Then for each $\varepsilon>0$, there exists $R>0$ so that
\begin{equation}\label{PFIUL1}
\int_{|x|>R}|\nabla u(t,x)|^2dx\leq \varepsilon,\;\textnormal{for all}\;0\leq t<+\infty.
\end{equation}
\end{proposition}
\begin{proposition}\label{VI}  {\bf (Virial identity)} Let $\phi\in C^\infty_0(\mathbb{R}^N)$, $\phi \geq 0$ and $T>0$. For $R>0$ and $t\in [0,T]$ define
\begin{equation*}\label{DFZ}
z_R(t)=\int_{\mathbb{R}^N} R^2 \phi\left(\frac{x}{R}\right)|u(t,x)|^2dx,
\end{equation*}  
where $u$ is a solution of \eqref{INLS}. Then we have
\begin{equation}\label{FD}
	z'_R(t)=2R Im\int_{\mathbb{R}^N} \nabla\phi\left(\frac{x}{R}\right)\cdot\nabla u \bar{u}dx
\end{equation}
and 
\begin{align}\label{SD}
z''_R(t)&= 4\sum_{j,k} Re\int \frac{\partial u}{\partial_{x_k}} \frac{\partial \bar{u}}{\partial_{x_j}}\frac{\partial^2\phi}{\partial x_k\partial x_j}\left(\frac{x}{R}\right)dx-\frac{1}{R^2}\int |u|^2\Delta^2 \phi\left(\frac{x}{R}\right) dx \nonumber \\
& -\frac{2\alpha}{\alpha+2}\int |x|^{-b}|u|^{\alpha+2}\Delta \phi\left(\frac{x}{R}\right) dx+\frac{4R}{\alpha+2}\int \nabla (|x|^{-b})\cdot\nabla \phi\left(\frac{x}{R}\right) |u|^{\alpha+2}dx.
\end{align}
\end{proposition}

\ The proof of Proposition \ref{PFIUL} is identical to the one in \cite[Lemma $5.6$]{HOLROU}, so we omit the details. On the other hand, Proposition \ref{VI} will proved at end of this section. 

\ Applying the previous results we now show the rigidity theorem.
\begin{theorem}\label{RT}{\bf (Rigidity)} Suppose $u_0\in H^1(\mathbb{R}^N)$ satisfying
$$
E[u_0]^{s_c}M[u_0]^{1-s_c} <E[Q]^{s_c}M[Q]^{1-s_c}
$$
and 
$$
\|  \nabla u_{0} \|_{L^2}^{s_c} \|u_{0}\|_{L^2}^{1-s_c}<\|\nabla Q \|_{L^2}^{s_c} \|Q\|_{L^2}^{1-s_c}.
$$
If the global $H^1(\mathbb{R}^N)$-solution $u$ with initial data $u_0$ satisfies
$$
K=\{u(t)\;:\;t\in[0,+\infty)\}\; \textnormal{is precompact in}\; H^1(\mathbb{R}^N)
$$
then $u_0$ must vanish, i.e., $u_0=0$.
	
\begin{proof} The proof follows similar ideas as in our paper \cite{paper2}. It follows from Theorem \ref{TG} that $u$ is global in $H^1(\mathbb{R}^N)$ and
\begin{equation}\label{TR1} 
\|  \nabla u(t) \|_{L^2_x}^{s_c} \|u(t)\|_{L^2_x}^{1-s_c}<\|\nabla Q \|_{L^2}^{s_c} \|Q\|_{L^2}^{1-s_c}.
\end{equation}
	
\ Set $\phi \in C_0^\infty$ be radial, with
$$
\phi(x)=\left\{\begin{array}{cl}
|x|^2&\textnormal{for}\;|x|\leq 1\\
0&\textnormal{for}\;|x|\geq 2.
\end{array}\right.
$$
The relation \eqref{FD}, the H\"older inequality and \eqref{TR1} imply that
\begin{eqnarray*}
|z'_R(t)| &\leq & cR\int_{|x|<2R}|\nabla u(t)||u(t)|dx\leq cR\|\nabla u(t)\|_{L^2}\|u(t)\|_{L^2}\lesssim cR.
\end{eqnarray*}
Thus,
\begin{eqnarray}\label{FDI}
|z'_R(t)-z'_R(0)|\leq 2cR,\;\;\textnormal{for all }\;t>0. 
\end{eqnarray}
  
\ The idea now is to get a lower bound for $z''_R(t)$ strictly greater than zero and reach a contradiction. Indeed, we deduce (using the local virial identity \eqref{SD})
\begin{align}\label{SDz} 
z''_R(t)&= 4\sum_{j,k} Re\int \partial_{x_k} u \partial_{x_j}\bar{u}\frac{\partial^2\phi}{\partial x_k\partial x_j}\left(\frac{x}{R}\right)dx-\frac{1}{R^2}\int |u|^2\Delta^2 \phi\left(\frac{x}{R}\right) dx   \nonumber   \\
&  -\frac{2\alpha}{\alpha+2}\int |x|^{-b}|u|^{\alpha+2}\Delta \phi\left(\frac{x}{R}\right) dx+\frac{4R}{\alpha+2}\int \nabla (|x|^{-b})\cdot\nabla \phi\left(\frac{x}{R}\right) |u|^{\alpha+2}dx \nonumber  \\
& = 8 \| \nabla u \|^2_{L^2_x}- \frac{4(N\alpha+2b)}{\alpha+2} \left\||x|^{-b}|u|^{\alpha+2}\right\|_{L^1_x}+R(u(t)),
\end{align}
where
\begin{align*}
R(u(t))&= 4 \sum\limits_{j}Re\int \left( \partial^2_{x_j} \phi\left(\frac{x}{R}\right)-2\right )|\partial_{x_j}u|^2
  + 4 \sum\limits_{j\neq k}Re\int \frac{\partial^2\phi}{\partial x_k\partial x_j}\left(\frac{x}{R}\right)\partial_{x_k}u\partial_{x_j}\bar{u}\\
  &- \frac{1}{R^2}\int |u|^2\Delta^2\phi\left(\frac{x}{R}\right)+\frac{4R}{\alpha+2}\int \nabla(|x|^{-b})\cdot\nabla \phi\left(\frac{x}{R}\right)|u|^{\alpha+2}\\
  &+\int \left(\frac{ -2\alpha(\Delta\phi\left(\frac{x}{R}\right)-2N)+8b}{\alpha+2}\right) |x|^{-b}|u|^{\alpha+2}.
\end{align*}
In view of $\phi(x)$ is radial and $\phi(x)=|x|^2$ if $|x|\leq 1$, the sum of all terms in the definition of $R(u(t))$ integrating over $|x|\leq R$ is zero. Indeed, by the definition of $\phi(x)$ it is clear for the first three terms. In the fourth term we have 
 $$
 \frac{8}{\alpha+2}\int_{|x|\leq R} \nabla(|x|^{-b})\cdot x|u|^{\alpha+2}dx=\frac{8}{\alpha+2}\int_{|x|\leq R} -b|x|^{-b}|u|^{\alpha+2}dx,
$$ 
and adding the last term also integrating over $|x|\leq R$, we have zero\footnote{Since $\Delta \phi =2N$, if $|x|\leq R$.}. Hence,
\begin{eqnarray*}
|R(u(t))|&\leq& c\int_{|x|>R} \left(|\nabla u(t)|^2+\frac{1}{R^2}|u(t)|^2+|x|^{-b}|u(t)|^{\alpha+2}\right)dx
\end{eqnarray*}
\begin{eqnarray}\label{RESTO}
&\leq & c\int_{|x|>R} \left(|\nabla u(t)|^2+\frac{1}{R^2}|u(t)|^2+\frac{1}{R^b}|u(t)|^{\alpha+2}\right)dx,
\end{eqnarray}
where we have used that all derivatives of $\phi$ are bounded and  $|R\partial_{x_j}(|x|^{-b})|\leq c|x|^{-b}$ if $|x|>R$. 
 
 \ Using the fact that $K$ is precompact in $H^1(\mathbb{R}^N)$. By Proposition \ref{PFIUL}, given $\varepsilon>0$ there exists $R_1>0$ such that $\int_{|x|>R_1} |\nabla u(t)|^2\leq\varepsilon$. Also, by mass conservation \eqref{mass}, there exists $R_2>0$ such that $\frac{1}{R^2_2}\int_{|x|>R_2} | u(t)|^2\leq \varepsilon$. Finally, by the Sobolev embedding $H^1\hookrightarrow L^{\alpha+2}$, there exists $R_3$ such that $\frac{1}{R_3^b}\int_{|x|>R_3} |u(t)|^{\alpha+2}\leq c\varepsilon$.\footnote{Recalling that $\|u(t)\|_{H^1_x}$ is uniformly bounded for all $t>0$ by \eqref{TR1} and Mass conservation \eqref{mass}.} Taking $R=\max\{R_1,R_2,R_3\}$ and by \eqref{RESTO} we conclude
 \begin{equation}\label{RESTO1}
 |R(u(t))|\leq c\int_{|x|>R} \left(|\nabla u(t)|^2+\frac{1}{R^2}|u(t)|^2+\frac{1}{R^b}|u(t)|^{\alpha+2}\right)dx\leq c\varepsilon.
 \end{equation} 
Furthermore, Lemma \ref{LGS} (iii), \eqref{SDz} and \eqref{RESTO1} imply that
\begin{equation*}\label{RESTO2}
z''_R(t)\geq 16A E[u]-|R(u(t))|\geq 16AE[u]-c\varepsilon,
\end{equation*} 
where $A=1-w^{\frac{\alpha}{2}}$ and $w=\frac{E[v]^{s_c}M[v]^{1-s_c}}{E[Q]^{s_c}M[Q]^{1-s_c}}$. \\
Choosing $\varepsilon=\frac{8A}{c}E[u]$, with $c$ as in \eqref{RESTO1} we have
$$
z''_R(t)\geq 8AE[u].
$$ 
Thus, integrating the last inequality from $0$ to $t$ it follows that
\begin{equation}\label{SDI}
z'_R(t)-z'_R(0)\geq 8AE[u]t.
\end{equation}

\ Taking $t$ large, we obtain a contradiction with \eqref{FDI}, which can be resolved only if $E[u]=0$. This implies by Lemma \ref{LGS} (i) that $u\equiv 0$. 
\end{proof}
\end{theorem}

\ We end this section by showing Proposition\ref{VI}.
\begin{proof}[\bf{Proof of Proposition \ref{VI}}] 
Observe that $
\partial_t|u|^2=2Re( u_t\bar{u})=2Im(iu_t\bar{u}). 
$
Since $u$ satisfies \eqref{INLS} and using integration by parts, we have
\begin{eqnarray*}
z'_R(t)&=&2Im\int  R^2\phi\left(\frac{x}{R}\right)iu_t \bar{u}dx\\
&=&-2Im\int R^2\phi\left(\frac{x}{R}\right)\left(\Delta u\bar{u}+|x|^{-b}|u|^{\alpha+2}\right)dx\\
&=&2RIm\int \nabla \phi\left(\frac{x}{R}\right)\cdot\nabla u \bar{u}dx.
\end{eqnarray*}
Again using integration by parts and the fact that $z-\bar{z}=2iIm z$, it follows that
\begin{eqnarray*}
z''_R(t)&=& 2RIm\int \nabla \phi\left(\frac{x}{R}\right)\cdot \left( \bar{u}_t\nabla u +\bar{u}\nabla  u_t\right)dx\\
&=& 2RIm\left\{\sum_{j}\int \bar{u}_t\partial_{x_j}u\partial_{x_j}\phi\left(\frac{x}{R}\right)dx- u_t\partial_{x_j}\left(\bar{u}\partial_{x_j}\phi\left(\frac{x}{R}\right)\right)  dx \right\}\\
&=& 2RIm\left\{\sum_{j}2i Im \int \bar{u}_t\partial_{x_j}u\partial_{x_j}\phi\left(\frac{x}{R}\right) dx-\int\frac{1}{R}u_t \bar{u}\partial^2_{x_j}\phi\left(\frac{x}{R}\right) dx \right\}\\
&=&4R I_1+2I_2,
\end{eqnarray*}
where
$$
I_1=Im \sum_{j}\int \bar{u}_t\partial_{x_j}u\partial_{x_j}\phi\left(\frac{x}{R}\right)\;\;\textnormal{and}\;\;I_2=-Im \sum_{j}\int u_t\bar{u} \partial^2_{x_j}\phi\left(\frac{x}{R}\right) dx.
$$

\ In view of $u$ is a solution of \eqref{INLS}, we deduce
\begin{eqnarray*}
I_2&=&-Im\left\{\sum_{j,k}  \int i \partial^2_{x_k}u \bar{u}\partial^2_{x_j}  \phi\left(\frac{x}{R}\right) dx  \right\}-\sum_{j}\int |x|^{-b}|u|^{\alpha+2}\partial^2_{x_j}\phi\left(\frac{x}{R}\right) dx\\
&= & Im \left\{\sum_{j,k}  \int i \left( |\partial_{x_k}u|^2 \partial^2_{x_j}  \phi\left(\frac{x}{R}\right)+ \frac{1}{R}\partial_{x_k}u\bar{u}\frac{\partial^3\phi}{\partial x_k \partial x^2_j}\left(\frac{x}{R}\right) \right)dx\right\} \\
&  &  - \int |x|^{-b}|u|^{\alpha+2}\Delta \phi\left(\frac{x}{R}\right) dx\\
&=&\int \left(|\nabla u|^2-|x|^{-b}|u|^{\alpha+2}\right) \Delta \phi\left(\frac{x}{R}\right) dx+\frac{1}{R} \sum_{j,k} Re\int \partial_{x_k}u\bar{u} \frac{\partial^3\phi}{\partial x_k \partial x^2_j}\left(\frac{x}{R}\right)dx.
\end{eqnarray*}
Another integration by parts yields 
\begin{eqnarray}\label{VI1} 
I_2&=&\int \left(|\nabla u|^2-|x|^{-b}|u|^{\alpha+2}\right) \Delta \phi\left(\frac{x}{R}\right) dx-\frac{1}{2R^2}\sum_{j,k}\int |u|^2\frac{\partial^4\phi}{\partial x^2_k\partial x^2_j}\left(\frac{x}{R}\right)dx \nonumber \\
& =&\int \left(|\nabla u|^2-|x|^{-b}|u|^{\alpha+2}\right) \Delta \phi\left(\frac{x}{R}\right) dx-\frac{1}{2R^2}\int |u|^2\Delta^2\phi\left(\frac{x}{R}\right) dx.
\end{eqnarray}

\ We now evaluate $I_1$. It follows from the equation \eqref{INLS} and $Im (z)=-Im (\bar{z})$ that\footnote{using $Im (iz)=Re (z)$ and $\partial_{x_j}(|u|^{\alpha+2})=(\alpha+2)|u|^\alpha Re (\partial_{x_j}\bar{u} u) $.}
\begin{align*}
I_1&=-Im \sum_{j}u_t \partial_{x_j}\bar{u}\partial_{x_j} \phi\left(\frac{x}{R}\right) dx\\
&=-Im i\sum_{j}\left\{\int \left(    \Delta u + |x|^{-b}|u|^{\alpha}u \right)\partial_{x_j}\bar{u}\partial_{x_j} \phi\left(\frac{x}{R}\right) dx   \right\}\\
&= -Re \sum_{j,k}\int \partial^2_{x_k} u\partial_{x_j} \bar{u} \partial_{x_j} \phi\left(\frac{x}{R}\right) dx-\sum_j \int |x|^{-b}\partial_{x_j} \phi\left(\frac{x}{R}\right)|u|^{\alpha}  Re (\partial_{x_j}\bar{u} u) dx\\
&=-Re \sum_{j,k}\int \partial^2_{x_k} u\partial_{x_j} \bar{u} \partial_{x_j} \phi\left(\frac{x}{R}\right)dx-\frac{1}{\alpha+2}\sum_j \int |x|^{-b}\partial_{x_j}\phi\left(\frac{x}{R}\right)\partial_{x_j}(|u|^{\alpha+2}) dx\\
&\equiv A+B.
\end{align*}
Since $\partial_{x_j}|\partial_{x_k}u|^2=2Re \left(\partial_{x_k}u \frac{\partial^2\bar{u}}{ \partial x_k \partial x_j}\right) $ and applying integration by parts twice, we obtain
\begin{eqnarray*}
A&=& Re \sum_{j,k} \left\{  \int \left(\partial_{x_j}\phi\left(\frac{x}{R}\right) \partial_{x_k}u \frac{\partial^2\bar{u}}{ \partial x_k \partial x_j}+ \frac{1}{R}  \partial_{x_k}u    \partial_{x_j}\bar{u}  \frac{  \partial^2\phi}{\partial x_j\partial x_k}\left(\frac{x}{R}\right)  \right)dx \right\}\\
&=&-\sum_{j,k}\frac{1}{2R} \int |\partial_{x_k}u|^2\partial^2_{x_j}\phi\left(\frac{x}{R}\right) dx+\frac{1}{R} \sum_{i,j} Re \int \partial_{x_k}u    \partial_{x_j}\bar{u}  \frac{  \partial^2\phi}{\partial x_j\partial x_k}\left(\frac{x}{R}\right) dx \\
&=&-\frac{1}{2R} \int |\nabla u|^2 \Delta \phi\left(\frac{x}{R}\right) dx+\frac{1}{R} \sum_{i,j} Re \int \partial_{x_k}u    \partial_{x_j}\bar{u}  \frac{  \partial^2\phi}{\partial x_j\partial x_k}\left(\frac{x}{R}\right) dx.
\end{eqnarray*}

\ Similarly, integrating by parts we have
\begin{align*}
B&=\frac{1}{\alpha+2} \sum_j \left( \int \partial_{x_j}\phi\left(\frac{x}{R}\right)\partial_{x_j}(|x|^{-b}) |u|^{\alpha+2} dx+ \frac{1}{R} \int \partial^2_{x_j}\phi\left(\frac{x}{R}\right)|x|^{-b} |u|^{\alpha+2} dx \right) \\
&=\frac{1}{\alpha+2}  \int \nabla \phi\left(\frac{x}{R}\right)\cdot \nabla (|x|^{-b}) |u|^{\alpha+2} dx+ \frac{1}{R(\alpha+2)}  \int \Delta \phi\left(\frac{x}{R}\right)|x|^{-b} |u|^{\alpha+2} dx. \\
\end{align*}
Therefore,
\begin{align}\label{VI2}
I_1&=-\frac{1}{2R} \int |\nabla u|^2 \Delta \phi\left(\frac{x}{R}\right) dx+\frac{1}{R} \sum_{i,j}Re  \int \partial_{x_k}u    \partial_{x_j}\bar{u}  \frac{  \partial^2\phi}{\partial x_j\partial x_k}\left(\frac{x}{R}\right) dx\nonumber \\
& +\frac{1}{\alpha+2}  \int \nabla \phi\left(\frac{x}{R}\right)\cdot \nabla (|x|^{-b}) |u|^{\alpha+2} dx+ \frac{1}{R(\alpha+2)}  \int \Delta \phi\left(\frac{x}{R}\right)|x|^{-b} |u|^{\alpha+2} dx.
\end{align}
Combining \eqref{VI1} and \eqref{VI2}, we deduce \eqref{SD}, which completes the proof. 
\end{proof} 

\section*{Acknowledgments} 
L.G.F. was supported by CNPq/Brazil and FAPEMIG/Brazil and C.M.G. was supported by CAPES/Brazil.

\section*{Appendix}
\ In this short Appendix we check the conditions of admissible pair used in Section $4$ and $6$.

\ {\bf A.1.} We claim $\frac{3\alpha}{2-b}=\frac{6}{3-2s_c}<p<6$, i.e., $p$ (see \eqref{paradmissivel1}) satisfies the condition \eqref{CPA2} (and therefore \eqref{L2Admissivel}, since $\frac{6}{3-2s_c}>2$) for $N=3$. Indeed, $\frac{3\alpha}{2-b}<p \Leftrightarrow (4-2b)(\alpha-\theta)+\alpha<(4-2b)(\alpha+1-\theta)\Leftrightarrow \alpha<4-2b$, which is true by hypothesis. Moreover, $p<6 \Leftrightarrow \alpha(\alpha+1-\theta)<(4-2b)(\alpha-\theta)+\alpha\Leftrightarrow \alpha(\alpha-\theta)<(4-2b)(\alpha-\theta)\Leftrightarrow \alpha<4-2b$. \\

\ {\bf A.2.} We notice that $r$ defined in \eqref{claim21} satisfies \eqref{CPA2}, that is $\frac{2N}{N-2s_c}<r< \frac{2N}{N-2}$. Indeed $r< \frac{2N}{N-2}$ is equivalent to $
\alpha(N^2-4)<2(4-2b)+\alpha N(N-2)\Leftrightarrow \alpha <\frac{4-2b}{N-2}.$
Moreover, $r>\frac{2N}{N-2s_c}=\frac{N\alpha}{2-b}$ is equivalent to 
$
(N+2)(4-2b)>2(4-2b)+\alpha N(N-2)\Leftrightarrow \alpha<\frac{4-2b}{N-2}.
$
\\

\ {\bf A.3.} Note that $\bar{r}$ defined in \eqref{C2N22} satisfies assumption \eqref{H-s} with $N=2$, that is $\frac{2}{1-s_c}=\frac{2\alpha}{2-b}<\bar{r}\leq \left((\frac{2}{1+s_c})^+\right)'$. The first inequality is equivalent to $\frac{2\alpha}{\varepsilon}>\frac{2\alpha}{2-b}$ and this holds since $2-b-\varepsilon>0$. On the other hand by the definition of $\left((\frac{2}{1+s_c})^+\right)'$ (see \eqref{a^+}) we conclude $\bar{r}=\frac{2\alpha}{\varepsilon}\leq \left((\frac{2}{1+s_c})^+\right)'$. Similarly, it easy to see that $r$ defined in \eqref{C2N21} satisfies \eqref{Hsdefinition}.

\bibliographystyle{abbrv}
\bibliography{bibguzman}

\end{document}